\documentclass[12pt]{article}

\usepackage{amsmath}
\usepackage{amsfonts}
\usepackage{amsthm}
\usepackage{amssymb}
\usepackage[usenames,dvipsnames]{color}

\usepackage{pgf,tikz,pgfplots}
\usepackage{verbatim}

\newtheorem{theo}{Theorem}[section]
\newtheorem{prop}[theo]{Proposition}
\newtheorem{lemm}[theo]{Lemma}

\newtheorem{coro}[theo]{Corollary}

\theoremstyle{remark}
\newcommand{\eqco}{\setcounter{equation}{0}}

\newcommand{\allco}{\eqco}
\renewcommand{\emptyset}{\varnothing}

\newcommand{\ER}{{Erd\"{o}s-R\'enyi~}}

\newcommand{\mus}{\mu^{\perp}}
\newcommand{\bigblock}{M}
\newcommand{\Po}{{\cal P}}

\newcommand{\Q}{{\cal Q}}
\newcommand{\cX}{{\cal X}}
\newcommand{\cY}{{\cal Y}}
\newcommand{\cZ}{{\cal Z}}
\newcommand{\cS}{{\cal S}}
\newcommand{\cI}{{\cal I}}
\newcommand{\cJ}{{\cal J}}
\newcommand{\cU}{{\cal U}}
\newcommand{\cV}{{\cal V}}

\def\N{\mathbb{N}}
\def\cN{{\cal N}}

\def\R{\mathbb{R}}
\def\E{\mathbb{E}}
\def\Pr{\mathbb{P}}

\def\0{{\bf 0}}

\def\G{{\cal G}}
\def\cA{{\cal A}}
\def\cD{{\cal D}}
\def\cB{{\cal B}}

\def\cH{{\cal H}}

\newcommand{\tilh}{\tilde{h}}
\newcommand{\tzeta}{\tilde{\zeta}}

\newcommand{\ogamma}{\overline{\gamma^\infty}}
\newcommand{\hatgamma}{\gamma_*}
\newcommand{\ozeta}{\overline{\zeta}}
\newcommand{\otheta}{\overline{\theta}}
\newcommand{\okappa}{\overline{\kappa}}
\newcommand{\oQ}{\overline{Q}}

\newcommand{\oalpha}{\overline{\alpha}}

\def\cC{{\cal C}}

\renewcommand{\E}{\mathbb E \,}

\newcommand{\toas}{\stackrel{{{\rm a.s.}}}{\longrightarrow}}
\newcommand{\toL}{\stackrel{{L^2}}{\longrightarrow}}
\newcommand{\tocc}{\stackrel{{c.c.}}{\longrightarrow}}
\newcommand{\toccL}{
 ~ \substack{ c.c. \\ \xrightarrow{\hspace*{0.7cm}}
\\ L^2}}

\newcommand{\wmax}{w_{{\rm max}}}
\newcommand{\Leb}{{\rm Leb}}

\newcommand{\X}{{\cal X}}
\newcommand{\LL}{{\cal L}}
\renewcommand{\H}{{\cal H}}

\newcommand{\A}{{\cal A}}

\newcommand{\Var}{{\rm Var}}

\newcommand{\Y}{{\cal Y}}

\newcommand{\U}{{\cal U}}

\newcommand{\cK}{{\cal K}}
\renewcommand{\U}{{\cal U}}
\newcommand{\F}{{\cal F}}
\renewcommand{\G}{{\cal G}}

\newcommand{\eps}{\varepsilon}
\def\bdm{\begin{displaymath}}
\newcommand{\edm}{\end{displaymath}}
\def\benu{\begin{enumerate}}
\def\eenu{\end{enumerate}}
\def\beqn{\begin{equation}}
\def\eeqn{\end{equation}}
\def\be{\begin{equation}}
\def\ee{\end{equation}}
\def\bea{\begin{eqnarray}}
\def\eea{\end{eqnarray}}
\newcommand{\bean}{\begin{eqnarray*}}
\newcommand{\eean}{\end{eqnarray*}}
\newcommand{\bear}{\begin{eqnarray}}
\newcommand{\eear}{\end{eqnarray}}

\renewcommand{\epsilon}{\varepsilon}

\def\R{\mathbb{R}}

\def\qed{\hfill\hbox{${\vcenter{\vbox{
    \hrule height 0.4pt\hbox{\vrule width 0.4pt height 6pt
    \kern5pt\vrule width 0.4pt}\hrule height 0.4pt}}}$}}

\begin{document}
\title{\bf Limit theory of combinatorial
 optimization for random geometric graphs}
%\title{\bf Subadditive limit theory of random geometric graphs}

\author{
Dieter Mitsche$^{1}$  and 
Mathew D. Penrose$^{2}$\\
%\author{Mathew D. Penrose , {\normalsize{\em University of Bath}} }
{\normalsize{\em Institut Camille Jordan, Univ. Jean Monnet and  University of Bath}} }

%\date{}
\maketitle

%{\em Department of
%Mathematical Sciences, University of Bath, Bath BA2 7AY, United
%Kingdom:} {\texttt{ m.d.penrose@bath.ac.uk}  } \\

 \footnotetext{ $~^1$ Institut Camille Jordan, Univ. Jean Monnet, Univ. St Etienne, Univ. Lyon, France: {\texttt dmitsche@univ-st-etienne.fr} \\ Dieter Mitsche has been supported by IDEXLYON of Universit\'{e} de Lyon (Programme Investissements d'Avenir ANR16-IDEX-0005) and by Labex MILYON/ANR-10-LABX-0070.}

 \footnotetext{ $~^2$ Department of
Mathematical Sciences, University of Bath, Bath BA2 7AY, United
Kingdom: {\texttt m.d.penrose@bath.ac.uk} }

%Running title: {\bf Limit theory for random geometric graphs} \\

\footnotetext{ AMS classifications: 
Primary 05C80, 60D05, 60F15. Secondary 90C27, 60G55.}

\footnotetext{
Keywords:
%Key words and phrases:
random geometric graph, thermodynamic limit,
dense limit, subadditivity, independence number,
domination number, clique-covering number, sphere packing, travelling salesman problem, minimum-weight matching.}

\begin{abstract}

In the random geometric graph $G(n,r_n)$, $n$ vertices are placed randomly in Euclidean $d$-space and edges are added between any pair of vertices distant at most $r_n$ from each other.  We establish strong laws of large numbers (LLNs) for a large class of graph parameters, evaluated for $G(n,r_n)$ in the thermodynamic limit with $nr_n^d =$ const., and also in the dense limit with $n r_n^d \to \infty$, $r_n \to 0$.  Examples include domination number, independence number, clique-covering number, eternal domination number and triangle packing number. The general theory is based on certain subadditivity and superadditivity properties, and also yields LLNs for other functionals such as the minimum weight for the travelling salesman, spanning tree, matching, bipartite matching and bipartite travelling salesman problems, for a general class of weight functions with at most polynomial growth of order $d-\eps$, under thermodynamic scaling of the distance parameter.

\end{abstract}

%\newpage

\section {Introduction}

Random geometric graphs (RGGs) are a well-known baseline stochastic model for
 combinatorial structures with spatial or
  multivariate content. Starting with the seminal paper of
 Gilbert~\cite{Gilbert}, random geometric graphs have in recent 
  decades received a lot of attention as a model for large 
communication networks such as sensor
networks, see~\cite{Akyildiz}. Network
 agents are represented by the vertices of the graph, and direct connectivity is represented by edges.  
Applications arise in other fields
 including 
%communications network engineering,
 theoretical computer science,
 geography, biology, topological
data analysis, network science and astronomy -
for more on applications of random geometric graphs, we refer to Chapter 3
of~\cite{Hekmat}, as well as to \cite{RGG}.

The classical random geometric graph  \cite{RGG}, also called
 the Gilbert graph, has its vertex set
 %$G(n,r)$ is
 given by taking  the points
of a sample of size $n$ from some specified probability distribution in Euclidean $d$-space, and edges between any two points distant at most
 $r$ from each other.
 %apart.
 In the terminology of \cite{RGG}, a
% {\em critical} or
 {\em thermodynamic} limiting regime
 involves taking $r=r(n) =
\Theta(n^{-1/d})$ as $n$ grows large, so that average degrees remain
bounded away from zero and infinity as $n \to \infty$, while
a {\em dense} limiting regime is one
 with $ 1 \gg r(n) \gg n^{-1/d}$, i.e. $n r(n)^d \to \infty$
%as $n \to \infty$ 
and $r(n) \to 0$ as $n \to \infty$.
%This is an important limiting regime, and is our main object of
%study here.

The first order (and in some cases, second order) limit theory, for RGGs in
the thermodynamic limit, of quantities such as number of edges,
 number of components, and number of isolated vertices,
%and chromatic 
is described in \cite{RGG}.  Loosely speaking, these enjoy
linear growth in $n$ because they are the sum of locally determined
quantities. The
order of the  largest component (also considered in \cite{RGG}) also
falls into this category if $n r^d$ is above a certain continuum percolation
threshold, in which case we say it is {\em supercritical}.

  In the present paper we derive laws of large numbers (LLNs),
 in thermodynamic or dense
limiting regimes, for certain other types of graph invariant. These 
include independence number, domination number, 
and many others (to be listed shortly).
Each of these quantities
is obtained as the solution of some constrained optimization problem in
the graph, and enjoys linear growth in $n$ in the thermodynamic regime.
 If $n r^d$ is subcritical,
they can be given as a sum of locally determined quantities,
 but if $n r^d$ is supercritical
they cannot.

 Our LLNs (Theorems \ref{th:therm2} and \ref{th:densegen}) come with
a rather simple set of conditions, applicable to a wide range of parameters 
of interest for RGGs, for which LLNs were previously available only in
the special case of domination number with $d=2$ and $nr(n)^d \gg \log n$
\cite{Lozier}, which we shall relax to $n r(n)^d \gg 1$. Moreover,
the dense 
limit theory relates these graph parameters to certain classic quantities of interest in deterministic combinatorial optimization such as sphere packing and sphere covering densities.

%The importance of this work lies in the presentation,
% in Theorems \ref{th:therm2} and \ref{th:densegen}, of general 
% laws of large  numbers 
%(LLNs) for both the thermodynamic and dense limits,
%with a rather simple set of conditions, applicable
%to a wide range of parameters of interest for RGGs,
% for which  LLNs in these limiting regimes
% were previously available  only in the special case
%of domination number  
% with $d=2$ and $nr(n)^d \gg \log n$ \cite{Lozier}, which we shall
%relax to $n r(n)^d \gg 1$. 
%Moreover, the dense limit theory relates 
%these graph parameters to certain classic quantities of interest in
%deterministic combinatorial optimization such as sphere packing
%and sphere covering densities.  

Furthermore, we shall use our methods to derive LLNs for
 weighted Travelling Salesman, Minimum Spanning Tree, 
 Minimum Matching  and Minimum Bipartite Matching problems
with edge weights determined by
 inter-point distances, via some arbitrary 
weight function that is either bounded, or grows at most
polynomially of order $d-\eps$, 
for some $\eps>0$, under thermodynamic scaling  
(for the unbounded case we require $\mu$ to have compact support).
These results (Theorems \ref{lem:TSP}, \ref{lem:MM1}, \ref{MBMthm}
and \ref{lem:MST1})
 generalize earlier work
on these problems, see e.g. \cite{BHH,Yukich},
 in which only power law weight functions were considered.

 As in the earlier work such as  
% on LLNs for Travelling Salesman type problems on random points such
%as
\cite{BHH,Yukich}, we use 
methods based on {\em subadditivity}. Here, however, 
we develop a method for using subadditivity
which does not require the spatial homogeneity assumptions
used in \cite{Yukich}; instead we use the thermodynamic
scaling of the extraneous distance parameter $r$.
This is what enables us to deal with the larger class
of weight functions than those considered in \cite{Yukich},
and moreover to apply subadditivity  to RGG functionals
of the
%
%Unlike the earlier works,  
%However, 
%has bewhich
%has previously been used to derive laws of large numbers (LLNs) for so-called
%{\em  Euclidean functionals} on random point samples;
% see \cite{Yukich}. Subadditivity
%does not appear to have been used before for RGG functionals of the
type considered here, for which subadditivity does not
appear to have been used before.
%which  are {\em not} Euclidean 
%in the sense of \cite{Yukich} because they involve an extraneous
%distance parameter $r$, whereas Euclidean functionals are completely
%determined by the sample.
%Also Euclidean functionals have a 
%homogeneity (scaling) property  which does not apply in the 
%present setting.
%However, by using the thermodynamic scaling of the parameter $r$ as a proxy
%for the homogeneity property, we are able to develop a limit theory
%analogous to that in \cite{Yukich} and apply it to the graph invariants
%already mentioned. 

Another methodology for deriving LLNs in the thermodynamic limit
was developed in \cite{PY}, based on {\em stabilization}. 
A point process functional is said to be stabilizing if it is a sum
of locally determined contributions from the points. 
Most of the graph parameters considered here do not seem to be 
stabilizing in the supercritical case.

We now define various graph parameters for an arbitrary finite graph 
$G=(V,E)$.  We are concerned here with the limit theory
 of these parameters for RGGs.
\begin{itemize}
\item{\emph{Independence number:}}
%Given a graph $G=(V, E)$, 
A set 
$A \subseteq V$ is said to be \emph{independent} (or {\em stable}), if for any $u,v \in A$,
 $uv \notin E$. The {\em independence number}
(or {\em stability number})
 of $G$,  here denoted $\alpha(G)$, is the maximum possible cardinality
 of an independent set $A \subseteq V$.

\item{\emph{Domination number:}}
%Given a finite graph $G$, its 
The {\em domination number} of $G$,
 here denoted $\gamma(G)$,  
 % $\gamma(G)$,  
is the minimum number $m$ such that there exists
a set $A$ of $m$ vertices with every
vertex of $G$ at graph distance at most
1 from $A$; such a set $A$ is known as a {\em dominating set} for $G$.

\item{\emph{Clique-covering number:}}
%Given a graph $G=(V,E)$, 
The clique-covering number of $G$, $\theta(G)$, is defined to be
 the chromatic number of the complement of $G$, that is, the minimum number of 
colours needed for colouring the vertices of $G$ in such a way that no two adjacent vertices in the complement of $G$ (i.e. non-adjacent vertices in $G$) obtain the same colour. In other words, $\theta(G)$ is the mimimum 
possible size (i.e. cardinality) of a {\em clique-partition}
of $V$, where by a clique-partition we mean a partition $\pi$ of $V$
such that for each set $W \in \pi$, the subgraph of $G$ induced by $W$
is the complete graph on $W$.

\item{\emph{Eternal domination number:}}
If $A \subset V$ is dominating, we may think of $A$ as representing
		the locations of a set of `guards', such that for any
		`attack' on an unguarded  vertex,
		there is a guard that can defend
		by moving to that vertex from an adjacent vertex.
		We say $A$ is {\em eternally
		dominating} if  guards, placed initially on the vertices 
		of $A$, can defend against any  finite or infinite 
		{\em sequence}
		of attacks (so after each attack, the defender can move a guard from an adjacent vertex to the attacked vertex leaving a new
		configuration which is also eternally dominating).
		Here we allow each attack in the sequence to be decided
		on the spot 
		%(rather than in advance) 
		by the attacker$^{1}$\footnotetext{$^{1}$
		There is also a variant, famous in computer science,
		where the whole attack sequence is given in advance, the 
		so called $k$-server problem.
		We will not consider this here.}.

%%Given a graph $G=(V,E)$, 
%A set $D \subseteq V$ is said to be 
%\emph{eternally dominating}
% if $D$ is a dominating set (we may think of $D$ as a subset of vertices such that one guard is located on each of them),
% and for any infinite sequence of attacks on unoccupied vertices (that is only decided on the spot
	%	$^{1}$\footnotetext{$^{1}$
	The eternal domination number of $G$,
 $\gamma^{\infty}(G)$, is the minimum number of vertices possible in an
 eternally dominating set of vertices $A$. See e.g. \cite{Klostermeyer2} for further discussion.
 %the initial set $D$. 

\end{itemize}
It is well known (see~\cite{Burger}) and easy to see that 
\bea
\gamma(G) \le
  \alpha(G) \le \gamma^{\infty}(G) \le \theta(G).
\label{graphineqs}
\eea 
We shall also consider some further graph parameters, namely
 vertex  cover  number, $H$-packing number, edge cover number, 
 number of components, and number of isolated vertices. We shall define
 these later on, in Section \ref{secfurther}.
%For two further graph parameters, namely {\em number of components}
%and {\em number of isolated vertices}, LLNs in the thermodynamic limit were
%previously obtained in \cite{RGG} when the underlying distribution has a density; here we shall extend these results to the case where this distribution
%has a singular part.

A second class of applications is to optimization problems over the
weighted complete graph on the sample of $n$ random points.
The classic example \cite{BHH} is the {\em Travelling salesman problem} (TSP):
find a tour (i.e., a Hamiltonian cycle)
  through the points of minimum total edge-length.
It is natural to consider  the analogous problem when the cost of
each edge $e$ is some {\em function} of the edge-length $\|e\|$,
denoted
$f(\|e\|)$ say, and one chooses the tour which minimizes
 the total cost. The case where $f(\|e\|) = \|e\|^{p}$ for some fixed constant
$p$ (i.e., power-weighted edges) 
has been much studied, see \cite{Yukich} and references therein,
 but many other
functions are available. One way to take the limit
is by thermodynamic scaling, i.e. consider a weight of the form
 $f(r_n^{-1} \|e\|)$ for a sample of size $n$, in the thermodynamic limit. 
Here we obtain LLNs under thermodynamic scaling
for a general  class
of weight functions $f$, not only for the TSP but also for
the minimum-weight matching (MM), minimum spanning tree (MST) and
minimum bipartite matching (BM) problems (for exact definitions of these problems, see Section \ref{secweighted}).
In the case of BM, and bounded $f$,
 our results go beyond what was previously available for power-weighted edges.

 We now describe some of the motivation and relevant past work.
On the one hand, the study of the parameters in a random setup is motivated by the fact that several of the decision problems studied here (that is, to decide whether a certain parameter is at least $k$), in particular independence number, domination number, clique-covering number, vertex cover and $H$-packing
 are well known to be $NP$-complete even for restricted graph classes, see~\cite{GareyJohnson, Kratochvil, Chan, Dor}. Therefore unless $P=NP$, one cannot
 expect a polynomial-time algorithm for such problems.
 %one can look at their typical sizes in a random setup as the current one.
 This motivates looking for polynomial-time algorithms
 which are `near-optimal 
 most of the time'. To quantify such terminology requires
 the study of these problems in a random setup (such as the current one).

We briefly discuss previous work on the corresponding problems for other
random graph models. Consider first  the \ER random graph $G(n,p)$ with $p=c/n$ for some constant $c$, which corresponds to our thermodynamic limit.
When $c<1$,
a weak LLN for any of the graph parameters considered here can
readily be obtained by computing the first two moments,  using
standard branching process approximation arguments for the exploration
process. In the case of independence number $\alpha(\cdot)$,
the value of the limit in the LLN is determined in \cite{Pittel},
where a central limit theorem is also provided. 
 For $c \geq 1$, a weak LLN
 for $\alpha(G(n,p))$
 was established in \cite{Bayati},
resolving a long-standing open problem. In fact, \cite{Bayati}
 gives this LLN for the `other' \ER model  $G(n,m)$
  with a deterministic number of edges $m$
and with $m \propto n$, but the result for $G(n,p)$ with $p= c/n$ can
then be readily derived from this. 
Similar LLNs for independence number have also been obtained for
the random $d$-regular graph \cite{Bayati} and the configuration
model \cite{Salez}. For $p = c/n$ with arbitrary $c$,
a LLN is known for the {\em maximum matching number},
which is a special case of the $H$-packing
number that we shall define later;
% which is the $H$-packing number with $H$ taken to be a complete
%graph on 2 vertices;
see \cite{KarpSipser}.

%, most results concern the subcritical case, that is, $c < 1$: in this regime, more precise results than ours (in fact, central limit theorems) are known for the parameters we study here~\cite{Bollobas,Pittel}. For $c \ge 1$, however, only for the number of components central limit theorems are known (in fact, since most of them are trees, they can be derived from results of~\cite{Pittel, Janson}). For other parameters studied here we are not aware of results similar to ours in that regime (excluding the number of isolated vertices, for which such results could easily be computed).

 On the other hand, the parameters have several practical applications.
Finding  dominating sets is important in
finding `central' or `important' sets of nodes
in a network, in contexts such as facility location~\cite{Haynes}, 
molecular biology~\cite{MMBP} for detecting significant proteins in 
protein-protein interaction networks, network controllability~\cite{Cowan}
and in wireless networks as centrality measure for efficient routing~\cite{SSZ}.
Dominating sets have attracted considerable attention in
the combinatorics literature; see the monograph \cite{Haynes}. 
% For a random geometric graph, it is clear that the domination number is
% not affected by how many times a repeated point is counted (though
% this might not be the
%case for a {\em soft} random geometric graph).
Domination numbers of random geometric graphs have been considered
in \cite{Lozier}. 

Finding large independent sets has different applications: for example vertices might represent intervals of tasks, and there is an edge between vertices if the corresponding intervals overlap; in job scheduling one likes to find a maximum number of jobs to be scheduled on one machine~\cite{Kolen}, corresponding to the independence number. Also, the problem of finding maximum independent sets in geometric graphs has been studied, for example, in the context of automatic label placement: given a set of locations in a map, find a maximum set of disjoint rectangular labels near these locations~\cite{Agarwal}.
% Even for special classes of graphs, the independence number is an NP-complete problem (see for example the introduction of~\cite{Chan}.}

Covering a graph with cliques can be seen as colouring the complement of the graph, and applications of colouring apply to the clique covering number. An iterative procedure of covering a graph with `almost' cliques was proposed in an influential work of~\cite{Song}, as a way of defining fractal dimension of networks.

 Regarding the eternal domination number,
its study was motivated by ancient problems in military defence (see for 
example~\cite{Ardilla}). It is known to be $NP$-hard~\cite{Klostermeyer}, but 
the decision problem of having an eternal domination number of at most $k$ is not known to be $NP$-complete since it is not known whether it belongs to
 $NP$:
 it is not clear how to confirm in polynomial time that a given initial configuration of guards can defend all possible sequences of attacks, see~\cite{Klostermeyer, Klostermeyer2}.

All of the classic optimization problems we consider here such
as the weighted TSP, MM and BM functionals, and 
 the weighted MST functional 
 have numerous applications in operations research (for example, in food delivery)
 and computer science. Among these, only the TSP
is NP-complete~\cite{GareyJohnson}; see \cite{Edmonds} and \cite{Leiserson} for polynomial-time
algorithms for the matching and MST problems respectively.

\textbf{Organization of the paper.} In Section~\ref{secstate} we state the
general results of this paper. In Section~\ref{secfurther} we then give 
applications of the general results to graph parameters; in
Section~\ref{secweighted} we show how to apply the results to classic
optimization problems such as TSP. Section~\ref{secproofgen} is devoted to the 
proof of the general results. Finally, Section~\ref{secdom} contains
the proofs of additional results about
the domination number in the dense regime.

\section{Statement of general results}
\label{secstate}

\subsection{Preliminaries}
\label{secprelim}

We now describe our setup in more detail.
Let $d \in \N$.
The $d$-dimensional  random geometric graph
$G(\X_n,r)$ is defined 
 as follows.
We take $\X_n$ to be a set of  $n$ independent identically
distributed (i.i.d.) random
 points in $\R^d$ with a specified common probability distribution
  $\mu$, and
assume $(r_n)_{n \geq 1}$ is a sequence of positive real numbers.
For any locally finite $\X \subset \R^d$ and distance parameter
 $r>0$, the graph
$G(\X,r)$ is defined to have vertex set $\X$, with
 two vertices  connected
by an edge if and only if their spatial
locations are at (Euclidean) distance at most $r$ from each other.

Throughout this paper  we assume the measure $\mu$ is diffuse (i.e.,
$\mu(\{x\}) =0$ for all $x \in \R^d$), so that the points of $\X_n$ are
almost surely distinct. (Actually,  our general results carry through
to the case  where $\mu$ is not diffuse,
provided we allow for $\X_n$ to be a {\em multiset},
that is we count each repeated point of $\X_n$
as many times as it arises, but we shall not pursue this further 
here.) Unlike in   \cite{RGG},
 we do {\em not} generally assume $\mu$ has a probability density 
function.

Let $\R_+ := (0,\infty)$. 
Given sequences $(a_n)_{n \geq 1}$ and
$(b_n)_{n \geq 1}$ taking
values in $\R_+$,
we write $a_n \ll b_n$, or $b_n \gg a_n$, or $a_n = o(b_n)$,
or $b_n = \omega(a_n)$,  to mean that
 $a_n/b_n \to 0$ as $n \to \infty$.
 Also we write
$a_n \sim b_n $ if $a_n/b_n \to 1$ as $n \to \infty$,
and $a_n = \Theta(b_n)$ if $a_n/b_n$
is bounded away from 0 and $\infty$.
For $x \in \R$ let $\lfloor x \rfloor$ denote the
integer part of $x$, i.e. the
largest integer not exceeding $x$.
Let $\lceil x \rceil$ denote the
smallest integer not less than $x$.
For any $x \in \R^d$ and $r>0$, let
 $B_r(x)$ denote the closed Euclidean ball (disk) centred
on $x$ of radius $r$. 
Let $o$ denote the origin in $\R^d$.
For $k \in \N$ set $[k]:= \{1,2,\ldots,k\}$.

For $s \in \R_+$, define $Q_s$,  
a half-open cube of side $s$ in $\R^d$ centred on the origin, by
$$
Q_s := [-s/2,s/2)^d. 
$$ 
%define the set $Q_r := [-r/2,r/2)^d$ 
%(a half-open cube of side $r$ in $\R^d$ centred on the origin).
Let $\Leb$ denote the Lebesgue measure on $\R^d$.
We set $\pi_d:= \Leb(B_1(o))$, the volume of the unit ball. 
%Define the set $Q_1:=[-\frac{1}{2},\frac{1}{2}]^d$ 
Let $\mu_U$ denote the uniform distribution on $Q_1$,
i.e. the restriction of Lebesgue measure
to $Q_1$.
% That is, set $\mu_U = \Leb(A \cap Q_1)$ for all Borel $A \subset \R^d$.
We write  $\mu = \mu_U$ in the case
where the common distribution $\mu$ of the points 
of $\X_n$ is this uniform distribution.

Given a sequence of random variables $\xi_n$ and a constant $c$,
we write $\xi_n \toas c$ (respectively $\xi_n \toL c$) 
if $\xi_n $ converges to $c$ almost surely
(respectively, in mean-square) as $n \to \infty$.
We write $\xi_n \tocc c$ if we have
{\em complete convergence} of $\xi_n$  to $c$, by which we
mean that for all $\eps >0$ we have
$\sum_{n=1}^\infty \Pr[|\xi_n - c| > \eps] < \infty$.
By the Borel-Cantelli lemma, if $\xi_n \tocc c$
then $\xi'_n \toas c$ for {\em any} sequence $(\xi'_n)_{n \geq 1}$
of  random variables on a common probability space with $\xi'_n$ having 
the same distribution as $\xi_n$ for each $n$.
For further discussion of complete convergence, see \cite{Yukich}.
If both $\xi_n \tocc c$ and $\xi_n \toL c$, we write $\xi_n \toccL c$.

Besides $\X_n$, we now define two further point processes, denoted
 $\Po_t$ and $\H_\lambda$. For $t >0$ we define the Poissonized point
 process $\Po_t$, coupled to $\X_n$ as follows.
Let $X_0,X_1,X_2,\ldots$ be a sequence of independent identically distributed
random $d$-vectors with common distribution $\mu$, and let $N_t$ be 
a Poisson random variable with mean $t$, independent of $(X_1,X_2,\ldots)$. 
Set $\X_n := \{X_1,\ldots,X_n\}$ and $\Po_t:= \{X_1,X_2,\ldots,X_{N_t}\}$.
Then $\Po_t$ is a Poisson process in $\R^d$ with intensity measure
$t \mu$ (see \cite[Proposition 3.5]{LP} or \cite[Proposition 1.5]{RGG}), 
coupled to $\X_n$.
Also, for $\lambda >0$,  let $\cH_\lambda$ be a homogeneous Poisson
 point process of intensity $\lambda$ on $\R^d$,  and for $s >0$
 set $\H_{\lambda,s} := \H_\lambda \cap Q_s$.
For any $\X \subset \R^d$ and $a >0$, we write $a \X$ for $\{ax: x \in \X\}$,
and for $y  \in \R^d$ set $y + \cX := \{y+x:x \in \cX\}$.
We write $|\cX|$ for the number
of elements of $\cX$.

\subsection{A class of functionals on point sets}

Let $\zeta(\cdot)$ be a 
non-negative real-valued function defined on the collection of 
all finite subsets of $\R^d$, with $\zeta(\emptyset)=0$, 
%(allowing multiplicities, as discussed earlier)
 and  having the following properties:
\begin{enumerate}
\item[P1]{\em Measurability:}
For all $k \in \N$,
the function $(x_1,\ldots,x_k) \mapsto \zeta(\{x_1,\ldots,x_k\})$
is measurable from $(\R^d)^k \setminus E$ to $\R$, where
$E$ denotes the set of $(x_1,\ldots,x_k) \in (\R^d)^k$ such
that $x_1,\ldots,x_k$ are not distinct.
\item[P2]
{\em Translation invariance}\label{propTI}: $\zeta(x+ \cX) = \zeta(\X)$ 
for all finite $\X \subset \R^d$ and  $x \in \R^d$.
\item[P3]
{\em Almost subadditivity:}\label{propsub} 
There exists a constant $c_1 \in [0,\infty)$,
% \in \mathbb{R}$, 
such that
$\zeta(\cY \cup \cZ) \leq
 \zeta(\cY) + \zeta(\cZ) +c_1$ for any two disjoint finite subsets
 $\cY,\cZ $  of $\R^d$. 

\item[P4]
{\em Superadditivity up to boundary:}\label{propsup} 
There exists a  constant $c_2 \in [0,\infty)$ such that
for any ordered pair $(\cY,\cZ)$  of disjoint finite sets in $\R^d$,
with  $\partial_{\cZ}(\cY)$ denoting the set of points of $\cY$ that lie
 at Euclidean distance at most $1$ from the set $\cZ$,
$$
\zeta(\cY \cup \cZ) \ge \zeta(\cY) + \zeta(\cZ) -
 %c \left( \partial_{\cY} ( \cZ)  + \partial_{\cZ} (\cY) \right)
 c_2 |\partial_{\cZ} (\cY) |.
$$
\end{enumerate}
Given the functional $\zeta$, and given a scaling  parameter $r>0$, we define
the scaled functional $\zeta_r(\cX)$ for all finite  
$\cX \subset \R^d$ by
$$
\zeta_{r}(\cX):=\zeta(r^{-1}\cX).
$$

Sometimes we shall assume $\zeta$ has one or more of the following
 further properties:

\begin{enumerate}
\setcounter{enumi}{3}
\item[P5] {\em Local sublinear growth:}
$\lim_{\delta \downarrow 0} ( \sup \{
\frac{\zeta(\X)}{|\X|}
: \X \subset B_\delta(o), \delta^{-1} \leq |\X|  \leq \infty \} ) =0.$
\item[P$5'$] {\em Local uniform boundedness:}\label{propLUB} It is the case that
$\sup_{\X \subset B_{1/2}(0), |\cX| < \infty} | \zeta(\X) | < \infty$.
\item[P6] {\em Upward monotonicity in $\cX$}: For all finite $\X \subset \R^d$
and $\Y \subset \R^d$ we have $ \zeta(\cX) \leq \zeta(\cX \cup \cY) $.
\item[P7] {\em Downward monotonicity in $r$}: For all finite $\X \subset \R^d$
we have  that   $ \zeta_r(\cX) $ is nonincreasing in $r$.
\item[P8] {\em RGG function:} $\zeta$ is a graph invariant of $G(\cX,1)$,
that is, $\zeta(\cX) = \zeta(\cY)$ whenever $G(\cX,1)$ is isomorphic
to $\zeta(\cY,1)$.
\end{enumerate}
Note that Property P$5'$ implies Property P5. Many of our examples
satisfy P$5'$.

For many of our examples, 
Property P8 holds, that is, $\zeta(\X)$ is some
specified graph invariant, evaluated on the geometric graph $G(\cX,1)$. Then
$\zeta_r (\X)$ is the same graph invariant evaluated for $G(\cX,r)$. 
Note that if P8 applies then both the 
 measurability condition P1 and the
translation invariance P2 
 follow automatically.

\subsection{General results for the thermodynamic limit}

Theorem \ref{th:therm2} below 
provides a law of large numbers, in the thermodynamic limit,
 for a general functional $\zeta$
 satisfying properties P1--P5,
 applied to $\X_n$.
The limit is expressed in terms 
of  a function $\rho: (0,\infty) \to [0,\infty)$, defined as follows for
all $\lambda >0$:
\bea
\rho(\lambda) :=
\lim_{s \to \infty}
%\lambda^{-1} s^{-d} 
 \E [ \zeta(\H_{\lambda,s}) /(\lambda s^d)].
% s^{-d} \E [ \zeta(\H_\lambda \cap Q_s) ].
% ~~~~~ \rho(\lambda) := \trho(\lambda)/\lambda.
\label{meanconv}
\eea

%Now we state our general thermodynamic limit theorem.
Given $(r_n)_{n \geq 1}$, we use $r_n$ as the scaling parameter
 when applying $\zeta$ to $\X_n$.
 Let $f_\mu$ denote the
density function (with respect to $\Leb$)
of the absolutely continuous part of $\mu$, and
 let  $\mus$  denote the  singular part of  $\mu$,
 in its Lebesgue decomposition.

\begin{theo}[General LLN in the thermodynamic limit]
\label{th:therm2}
Let $\zeta(\cdot)$ be a functional
 satisfying Properties  P1--P5. Then,
for all $\lambda > 0$, the limit in (\ref{meanconv}) exists in $\R$.
 %If also P5 holds} then
 Moreover, if  $nr_n^d \to t$ as $n \to \infty$, for some 
$t \in (0,\infty)$, we have 
\bea
n^{-1} \zeta_{r_n}(\X_n) \toccL \int_{\R^d} \rho(t f_\mu(x)) f_\mu(x) dx 
~~~{\rm as}~~ n \to \infty,
\label{0101c}
\eea
where for $\lambda >0$, we define $\rho(\lambda) $ by (\ref{meanconv}),
and moreover
\bea
n^{-1}  \zeta_{r_n}(\Po_n) \toccL \int_{\R^d} \rho(t f_\mu(x)) f_\mu(x) dx 
~~~{\rm as}~~ n \to \infty.
\label{0302a}
\eea
\end{theo}
{\bf Remark.}
 If $\zeta(\cdot)$ satisfies only P1--P4, then provided
 $\mu$ is absolutely continuous, the conclusion of Theorem \ref{th:therm2}
still holds. See Proposition \ref{propsofrho}(i) below, and
 the remark just before Lemma \ref{lemblock}.

In some examples, the functional $\zeta$ does not satisfy P3 and P4 because
the inequalities are the wrong way around, but can be obtained as a linear
combination of functionals which do satisfy P3 and P4, possibly including
the counting functional $\cX \mapsto |\cX|$. The following corollary deals
with these cases.

\begin{coro}\label{CorollaryMain2}
Suppose for some $c_3 \ge 0$ that we can write
 $\zeta(\cdot) = c_3 |\cdot| - \zeta'(\cdot) + \zeta''(\cdot)$, 
where $\zeta'(\cdot), \zeta''(\cdot)$ are  non-negative functionals both
satisfying Properties P1--P5.
Then, for all $\lambda >0$, the limit 
$\rho(\lambda)$ defined by (\ref{meanconv}) exists in $\R$.
Also, if  $nr_n^d \to t$ as $n \to \infty$, for some 
$t \in (0,\infty)$,
we have
\bea
n^{-1} \zeta_{r_n}(\X_n) \toccL \int_{\R^d} \rho(t f_\mu(x)) f_\mu(x) dx 
+ c_3 \mus
(\R^d) ~~~{\rm as}~~ n \to \infty. 
\label{0113b}
\eea
\end{coro}
\begin{proof}
%[Proof of Corollary~\ref{CorollaryMain2}]
For $\lambda >0$,
let $\rho'(\lambda), \rho''(\lambda)$ be defined analogously to
$\rho(\lambda)$ at (\ref{meanconv}), viz., 
as the large-$s$ limit of
$ 
	\E [ \zeta'(\H_{\lambda, s})/(\lambda s^d) ]$, 
 respectively $
	\E [ \zeta''(\H_{\lambda,s})/(\lambda s^d) ]$;
these limits exist by applying Theorem  \ref{th:therm2}
to $\zeta'$ and $\zeta''$, respectively. Hence
the limit $\rho(\lambda)$ defined by (\ref{meanconv}) also exists,
with  $\rho(\lambda)= c_3- \rho'(\lambda) + \rho''(\lambda)$.
Applying  
 Theorem~\ref{th:therm2} to  $\zeta'$ and $\zeta''$, 
 and using the fact that for any two sequences of random variables $(\xi_n)_n$ and $(\xi'_n)_n$, we have that $\xi_n \tocc c$ and $\xi'_n \tocc d$
 together imply
 $\xi_n+\xi'_n \tocc c+d$,
we obtain 
%under either hypothesis (a) or (b) 
the complete and $L^2$  convergence 
\bean
\lim_{n \to \infty} n^{-1} \zeta_{r_n}(\X_n) & =  &
\lim_{n \to \infty} (c_3 -  n^{-1} \zeta'_{r_n}(\X_n) + n^{-1} \zeta''_{r_n}(\X_n) ) 
\\
& = & c_3 + \int_{\R^d}
 [- \rho'(t f_\mu(x) ) +  \rho''(t f_\mu(x) ) ] f_\mu(x) dx,
\eean
which is equal to the right hand side of (\ref{0113b}),
as required. 
% The proof of (b) follows along the same lines.
\end{proof}

Our next result gives
some of the properties of the function $\rho(\cdot)$.

\begin{prop}[Properties of $\rho(\cdot)$]
\label{propsofrho}
Suppose $\zeta$ satisfies Properties P1--P4. Then:

(i) 
For all $\lambda > 0$ the limit in (\ref{meanconv}) exists in $\R$.
In fact,
if $((s_n,\lambda_n))_{n \geq 1}$
is any sequence in $\R_+^2$ satisfying $s_n \to \infty$
and $\lambda_n \to \lambda$ as $n \to \infty$, we have
\bea
\lim_{n \to \infty} s_n^{-d} \E[ \zeta(\cH_{\lambda_n} \cap Q_{s_n}) ] 
= \lambda \rho(\lambda).
\label{meanconvstrong}
\eea

(ii) The function $\lambda \mapsto \lambda \rho(\lambda)$ 
is Lipschitz continuous on $(0,\infty)$ with Lipschitz constant
at most $K := \max(c_1 + \zeta(\{o\}),c_2 - \zeta(\{o\}) )$, where
$c_1,c_2$ are the constants in P3 and P4 respectively. Hence
the function $\lambda \mapsto \rho(\lambda)$ is also continuous.

(iii) 
%$\trho(\lambda) \leq (c_1 + \zeta(\{o\})) \lambda$ ( so
 For all $\lambda >0$ we have 
 $\rho(\lambda) \leq c_1 + \zeta(\{o\})$.

(iv) If  $\zeta$ satisfies P6, then the function
  $\lambda \mapsto \lambda \rho(\lambda)$ is nondecreasing.

(v) If $\zeta$ satisfies P7, then the function
  $\lambda \mapsto \rho(\lambda)$ is nonincreasing.

(vi) If we can take $c_1=0$ in P3, 
then $\rho(\lambda) \to \zeta(\{o\})$ as $\lambda \downarrow 0.$

(vii) If $\zeta(\{o\}) >0$, then 
 $\lambda \mapsto \lambda \rho(\lambda)$ is strictly increasing for
$0 < \lambda < \zeta(\{o\}) /(c_2 \pi_d).$
\end{prop}

Let $\lambda_c := \lambda_c(d)$ be the critical value for
 continuum percolation,
namely the supremum of all $\lambda$ such that
all components of
$G(\cH_\lambda,1)$ are finite.
It is known that $\lambda_c(d) \in (0,\infty)$ for $d \geq 2$ and
$\lambda_c(1)= +\infty$. See \cite{RGG} for further discussion.
%Our next result provides an alternative characterization
%of $\rho(\lambda)$ when $\lambda < \lambda_c$.

\begin{prop}[Alternative characterization of
	$\rho(\lambda)$ for subcritical $\lambda$]
\label{thsubcrit}
Suppose $\lambda \in (0,\lambda_c)$. Suppose $\zeta$ satisfies Properties
P1--P4, with $c_1=0$ in P3. 
Then the constant $\rho(\lambda)$ given at (\ref{meanconv})
satisfies
\bea
\rho(\lambda) = \E [ \zeta(\cC_o(\lambda)) /|\cC_o(\lambda)| ],
\label{rhosubcrit}
\eea
where $\cC_o(\lambda)$ denotes the set of vertices of the component containing
$o$ of the graph $G(\cH_\lambda \cup \{o\},1)$.
\end{prop}

 It follows from Proposition \ref{thsubcrit}
 that if the hypotheses of Theorem \ref{th:therm2} 
 hold with $c_1=0$, and also
 $t f_\mu(x) < \lambda_c$ for Lebesgue-almost all $x \in \R^d$ then the
limit in (\ref{0101c}) is equal to
$\int_{\R^d} \E[\zeta(\cC_o(tf_\mu(x))) /|\cC_o(t f_\mu(x))|]f_\mu(x)dx$.
In the case where $\mu$ is absolutely continuous,
this can also be proved using the methods of \cite{PY}.

If the hypotheses of Theorem \ref{th:therm2} 
 hold with $c_1=0$, and also
 $\sup_{x \in \R^d} t f_\mu(x) < \lambda_c$, then
we expect that one could
prove a Gaussian limit for $\zeta_{r_n}(\X_n)$ (suitably scaled and
centred),
using the methods of \cite{PenEJP} for example. However, this is beyond
the remit of the present paper.

\subsection{General results for dense limiting regimes}
In Theorem
\ref{th:lambig} below, we 
  provide information about the limiting behaviour of 
$\rho(\lambda)$ for large $\lambda$, under the extra  conditions P$5'$ (local uniform boundedness),
  P6 (monotonicity) and P8 ($\zeta$ is a RGG functional).
 This complements the information provided 
for small $\lambda$ in parts (vi) and (vii) of Proposition \ref{propsofrho}.
 Before stating this result, we first need
to describe the deterministic limiting behaviour of $\zeta^*(Q_s)$, where
for  $A \subset \R^d$, 
we set
%Next we give a deterministic consequence of subadditivity,
\bea
\zeta^*(A) : = \sup \{\zeta(\cX): \cX \subset A, |\cX| < \infty\}.
\label{gamstardef}
\eea
\begin{lemm}
\label{lemdet}
 Suppose $\zeta(\cdot)$ satisfies P1--P4 and P$5'$.
Then the limit
\bea
  \overline{\zeta} :=
\lim_{s \to \infty} 
 \left( s^{-d} \zeta^*(Q_s) \right) 
\label{wellknown2}
\eea
exists in $[0,\infty)$.
Also for all $\lambda >0$, we have
\bea
\lambda \rho(\lambda) \leq \overline{\zeta}. 
\label{0922a}
\eea
\end{lemm}

\begin{theo}
\label{th:lambig}
Suppose $\zeta$ satisfies Properties P8, P3, P4, P$5'$ and P6.
 Then
\bea
\lim_{\lambda \to \infty} (\lambda \rho(\lambda)) = \ozeta.
\label{eq:lambig}
\eea
\end{theo}
Our main result in this subsection gives the limiting behaviour 
of $\zeta_{r_n}( \cX_n) $ in  the
 {\em dense limit} with $n r_n^d \to \infty$ and $r_n \to 0$.
 %In this case we consider only the case with $\mu = \mu_U$. 
Here we consider only the case where $\mu$ has
 a density $f_\mu$, and where moreover
$f_\mu^{-1} ((0,\infty))$ is Riemann measurable. 
Recall that we say a set $A \subset \R^d$ is {\em Riemann 
measurable} if ${\bf 1}_A$ is Riemann integrable,
i.e. if $A$ bounded and $\partial A$ has Lebesgue measure zero.
\begin{theo}[General LLN in the dense limit]
\label{th:densegen}
Suppose $\zeta$ satisfies Properties P3--P6 and P8.
Suppose $\mu$ is absolutely continuous, and $f_\mu$
can be chosen in such a way that $f_\mu^{-1}((0,\infty))$ is
Riemann measurable.
Suppose $r_n \to 0$ and 
$nr_n^d \to \infty$ as $n \to \infty$. Then, almost surely 
$$
\lim_{n \to \infty} r_n^d \zeta_{r_n}(\cX_n) = \ozeta
\;
\Leb(f_\mu^{-1}((0,\infty))).
$$
\end{theo}

We shall prove the general results stated so far
 in Section \ref{secproofgen}.

\section{Applications to graph parameters}
\label{secfurther}
\allco

In this section, we seek to apply  Theorems \ref{th:therm2},
\ref{th:lambig} and \ref{th:densegen} to the four graph parameters 
described in the Introduction,
and also to the following graph parameters:
\begin{itemize}
\item{\emph{Vertex cover number:}}
A set $W \subset V$ is called a {\em vertex cover} of $G$ if
every edge of $G$ is incident to at least one vertex in $W$. The
{\em vertex cover number} of $G$ is the smallest possible number of
vertices required for a vertex cover of $G$.

\item{\emph{$H$-packing number:}}
Let $H$ be a fixed  connected graph  with $k$ vertices, $2 \leq k < \infty$.
Given $m \in \N$,  
we refer to any collection of $m$
 vertex-disjoint $H$-subgraphs of $G$ as an $H$-{\em packing} of size $m$.
 Here an $H$-subgraph means a subgraph isomorphic to $H$
(it does not need to be an induced subgraph).
 The  {\em maximum $H$-packing number} $\psi_H(G)$ is defined to be the largest
possible size of $H$-packing in $G$.
In a remark just after Theorem \ref{lem:Packing}, we shall
describe a generalization that allows for more than one $H$.

\item{\emph{Edge cover number:}}
If $G$ has no isolated vertex,
 the edge cover number $\eta(G)$ is the
smallest number of edges that can be selected
 such that every vertex is incident to
at least one of the selected edges. 
For general $G$, we let $\eta(G)$ denote the edge cover
number of the graph obtained from $G$ by removing all isolated
vertices.
		\item {\em Number of connected components}, and 
			{\em number of isolated vertices}. 
			The definitions of these graph
			parameters are well known.
\end{itemize}

%two further
%graph parameters, namely the number of connected components and the
%number of isolated vertices. 
We consider in turn each of these graph parameters,
evaluated on the geometric graph $G(\cX,1)$.
%We also consider 
Clearly, in each case this gives a functional $\zeta(\cdot)$
%all nine of these graph parameters
satisfying Property P8, and hence both P1 and P2.
% are translation-invariant, so Property P2
%\ref{propTI}
% are satisfied for all of them.
 Also, each of these graph parameters, %all of them
 except vertex cover number
and edge cover number, satisfies Property P7, because it  can only be
reduced by adding an edge.

\subsection{Independence number}
\label{subsecindep}

For finite $\cX \subset \R^d$, let $\alpha(\cX) :=\alpha(G(\cX, 1))$, the
independence number of the geometric graph $G(\X,1)$. 
Thus if we take $\zeta(\cX) = \alpha(\cX)$, then $\zeta_r(\X) =
\alpha(G(\X,r))$ for all $r>0$.   
Let $\cB$ denote the class of bounded Borel subsets of $\R^d$.
Given $A \in \cB $, set
%We now discuss a
%notion of {\em packing density} which is
% in some sense dual to that of covering density, and is
%also adapted from \cite{Rogers}. For bounded $A \subset \R^d$,  define
\bea
\alpha^*(A)
 & := & \sup\{\alpha(\X): \X \subset A, |\cX | < \infty \}
\nonumber \\
 & = & \sup\{|\X|: \X \subset A, |\cX | < \infty, E(G(\X,1)) = \emptyset \},
\label{packdensfin}
\eea
where $E(\cdot)$ denotes the set of edges of the graph in question.
Thus $\alpha^*(A)$ is the maximum number of disjoint closed
 balls of radius $1/2$
which can be packed into $A$ (or at least, with their  centres in $A$),
which is clearly finite.
\begin{theo}[LLNs for independence number of RGG] 
\label{th:ind}
Theorems~\ref{th:therm2}, \ref{th:lambig} and \ref{th:densegen} all apply
 when choosing $\zeta(\cX) := \alpha(\X)$.
% to be the independence number $\alpha(G(\cX, 1))$.
\end{theo}
{\bf Remarks.} (a) A Poissonized version of the case $\mu= \mu_U$ of
Theorem \ref{th:therm2} for the present 
choice of $\zeta(\cdot)$ was given in \cite[Theorem 2.7]{PY}, but
that result required the limit of $n r_n^d$ to be
 below the percolation threshold
$\lambda_c$.

(b) With the current choice of $\zeta$, the quantity 
$\lambda \rho(\lambda)$ defined
at (\ref{meanconv}) is the
  intensity of a maximum hard-core thinning of the restriction of
 $\cH_\lambda$ to a box, in the limit of a large box.
Here, by `hard-core thinning' we mean a sub-point process with all
inter-point distances greater than 1. 
A related quantity is the maximum possible intensity of a 
{\em hard-core  stationary thinning}
of the whole of $\cH_\lambda$; we denote this quantity here by
 $\hat{\rho}(\lambda)$.
See \cite{HL}  for further details on hard-core stationary thinnings.
 It is not hard to see that $
\hat{\rho}(\lambda) \leq \lambda \rho(\lambda) $,
 and we conjecture that in fact equality holds here. 

(c) By Lemma \ref{lemdet}, the limit 
\bea
\oalpha := \lim_{s \to \infty} (s^{-d} \alpha^*(Q_s)), 
\label{opsidef}
\eea
exists in $\R$. 
It  is the optimal (i.e., maximal) {\em packing density} of balls of radius 1/2 in $\R^d$, where
we measure density here by `number of packed balls per unit volume', not `volume
of packed balls per unit volume'. 
For the present choice of $\zeta$,
the limit $\ozeta$ appearing in the statement
of Theorems \ref{th:lambig} and \ref{th:densegen} is equal to  $\oalpha$.

It can be seen that $\oalpha =1$ for $d=1$,
and $\oalpha = \sqrt{4/3}$ for $d=2$ by Thue's theorem on disk packing. 
As a consequence of the Kepler conjecture \cite{Hales},
 for $d=3$ we have $\oalpha = \sqrt{2}$.
For further discussion of packing densities, see \cite{Rogers,Conway}. 

\begin{proof}[Proof of Theorem \ref{th:ind}]
Given disjoint finite $\cY, \cZ \subset \mathbb{R}^d$, and   given
  independent $\cA \subset \cY \cup \cZ$ with
$|\cA| = \alpha(\cY \cup \cZ)$, the set $\cA \cap \cY$ is independent in 
$G(\cY,1)$ and $\cA \cap \cZ$ is independent in $G(\cZ,1)$. 
Therefore
 $\alpha(\cY \cup \cZ) \le \alpha(\cY)+\alpha(\cZ)$, so
 Property P3 (almost subadditivity) holds with $c_1=0$. 

To check Property P4,
%\ref{propsup},
%let  $\cY, \cZ,$ be finite subsets of $\R^d$. 
let $\cI \subset \cY$
be an  independent set of $G(\cY,1)$ with 
 $|\cI | = \alpha(\cY )$, 
 and let $  \cJ \subset  \cZ$
be an independent set of $G(\cZ,1) $ with
 $|\cJ | = \alpha(\cZ )$. 
Let $\cI' := \cI \setminus \partial_{\cZ} \cY$. 
% Then the number of points so removed is at most $|\partial_{\cZ}\cY|$. Also
Then  $\cI' \cup \cJ$ is an independent set in $G(\cY \cup \cZ,1)$. Hence
\bean
\alpha(\cY \cup \cZ)  & \geq  & | \cI' \cup \cJ| = |\cI'| + | \cJ|
\\
  & \geq  & |\cI| -  
|\partial_{\cZ} \cY| + |\cJ| = \alpha(\cY) + \alpha(\cZ) - |\partial _{\cZ} \cY|, 
\eean
 which gives us Property P4 with $c_2=1$.
%
% Property  (nonnegativity) is obvious.
Also
for any finite 
 nonempty 
$\cX \subset B_{1/2}(o)$,
$G(\X,1)$ is a complete graph so
 $\alpha(\X)  \leq 1$. 
% Q_1$, $\alpha(\X) \leq \alpha^*(Q_1)$ by
%$\cX \subset Q_1$, $\alpha(\X) \leq \alpha^*(Q_1)$ by
%(\ref{packdensfin}),  and $\alpha^*(Q_1) < \infty$, 
Hence  P$5'$
 (locally uniform boundedness), and hence also P5, holds here.
 Thus  Theorem \ref{th:therm2} applies.  

Property P6 %and P7 
clearly holds here, since
 the independence number can only be increased by adding vertices.
%and can only be decreased by adding edges.
%Since P6 holds, 
Therefore Theorems \ref{th:lambig} and \ref{th:densegen} also apply.
 \end{proof}

\subsection{Domination number}
\label{subsecdom}

For finite $\X \subset \R^d$, set
 $\gamma(\cX) :=\gamma(G(\cX,1))$,
 the domination number 
of $G(\X,1)$.
Then taking $\zeta(\cdot) := \gamma(\cdot)$, we have
 $\zeta_r(\cX)= \gamma(G(\X,r))$ for all $r >0$.

For this choice of $\zeta$, it is shown in  \cite[Theorem 2(c)]{Lozier} that  
when $d=2$ and $\mu= \mu_U$, 
if $n r_n^2$ converges to a positive finite limit,
then $ n^{-1} \zeta_{r_n}(\X_n)  = \Theta(1) $ in probability.
 Our next result improves on this by 
showing almost sure convergence to a limit, and
allowing for general $d$ and $\mu$.

\begin{theo}[LLN for domination number of RGG in the thermodynamic limit]
\label{th:dom}
Theorem~\ref{th:therm2} holds when choosing
 $\zeta(\cX)$ as $\gamma(\cX)$.
\end{theo}
\begin{proof}
Setting $\zeta(\cX) :=\gamma(\cX)$,
we need to check Properties P3--P5.
For P3 (almost subadditivity), let
 $\cY, \cZ$
be  disjoint finite
subsets
 of $\R^d$. Taking a minimum dominating
 set for the graph induced by $\cY$ together with a minimum dominating set for the graph induced by $\cZ$ yields a dominating set for the graph induced by $\cY \cup \cZ$, and hence $\gamma(\cdot)$ satisfies Property P3
 with $c_1=0$. 
Also by Lemma \ref{lemdomsup}, which we give at the end of this subsection,
$\gamma(\cdot)$ satisfies
Property P4 (superadditivity up to boundary). 

Finally 
$\gamma(\X) \leq 1$ for all $\X \subset B_{1/2}(o)$, so
P$5'$ (and hence  P5) holds for $\gamma(\cdot)$. 
Thus, $\gamma(\cdot)$ satisfies all the conditions
for Theorem \ref{th:therm2}.
\end{proof}

  Property P6 (upward
monotonicity in $\cX$) does {\em not} hold here,
 since adding vertices might make
 the domination number smaller. 
Therefore
 Theorems \ref{th:lambig} and 
\ref{th:densegen}, concerning
high density limits, are not applicable here. Nevertheless, we
are able to provide certain similar results for
domination number too.
First we need to recall the notion of {\em covering density}, 
 adapted from \cite{Rogers}.

For $A \in \cB$
(the class of bounded Borel subsets of $\R^d$),
 let $\kappa (A)$ denote the smallest possible 
number of  unit radius balls required to cover $A$, i.e.
\bea
\kappa (A)= \inf \{|{\cX} |: \cX \subset \R^d, |\cX| < \infty,
 A \subset \cup_{x \in \X} B_1(x)\}.  
\label{kappadef}
\eea
Then $\kappa (A \cup A')
\leq \kappa(A) + \kappa(A')$ for any $A,A' \in \cB$
with $A \cap A' = \emptyset$, since for any covering of $A$ and
any covering of $A'$, they can be combined to cover $A \cup A'$.
 Also $\kappa(A) \leq \kappa(Q_1)< \infty$
for all Borel $A \subset Q_1$.
 Hence by a simple deterministic version of the subadditive limit theorem
 (Lemma \ref{lemsubadd} below), 
\bea
 \lim_{s \to \infty} 
 \left(s^{-d}  
\kappa(Q_s)
 \right) 
= \inf_{s \geq 1} \left( s^{-d} \kappa(Q_s) \right)
=:  \okappa.
\label{psidef}
\eea
The quantity $\okappa$ is the optimal (i.e., minimal) covering density of unit balls in $\R^d$.
 We can now state our results for domination number that
are analogous to Theorems \ref{th:lambig} and \ref{th:densegen}.

\begin{theo}
\label{th:onemore}
Choosing $\zeta(\X)$ as $\gamma(\cX)$,
and defining $\rho(\lambda)$ by (\ref{meanconv}), we have
that
 $\lim_{\lambda \to \infty} (\lambda \rho(\lambda)) = \okappa$.
\end{theo}
\begin{theo}[LLN for domination number of RGG in the dense limit] 
\label{thsubcond}
Suppose $\mu = \mu_U$ and  $n^{-1/d} \ll r_n \ll 1$.
Then 
\bea
r_n^d \gamma(G(\X_n,r_n)) \toas \okappa ~~~ 
{\rm as } ~ n \to \infty.
\label{1221h}
\eea
%where $B$ is the Euclidean unit ball in $d $ dimensions.
\end{theo}
When $d=2$,
it is known that $\okappa =   \sqrt{4/27} $.
 See \cite{Kershner}, or 
\cite[page 16]{Rogers}. Therefore, Theorem \ref{thsubcond} extends
\cite[Theorem 2(b)]{Lozier}, which is concerned only with convergence
in probability in the case
where $d=2$ and 
$r_n \gg ((\log n)/n)^{1/2}$. 

We defer the proofs of Theorems \ref{th:onemore} and \ref{thsubcond}
to Section \ref{secsubd}.
It is likely that Theorem \ref{thsubcond} can be generalized from
the special case given here with $\mu=\mu_U$, to general absolutely
continuous $\mu$ such that $f_{\mu}$   
can be taken to be  Lebesgue-almost
 everywhere continuous with $f_\mu^{-1}((0,\infty))$ 
 Riemann measurable. In this generality the limit at
(\ref{1221h}) should be $\okappa \; \Leb
 (f_\mu^{-1}((0,\infty)))$. 

We now give some relations between $\okappa$, and $\ozeta$
defined at (\ref{wellknown2}), 
and the packing density $\oalpha$ from the previous
subsection.

 \begin{prop}
\label{corocov}
Choosing $\zeta(\cX)$ as $ \gamma(\X)$,
we have
 $\okappa \le \oalpha = \overline{\zeta}$.
 \end{prop}
 \begin{proof}
Given $A \in \cB$, we have $\kappa(A) \leq \alpha^*(A)$. Indeed,
for any $\X \subset A$ with
$|\X|= \alpha^*(A)$
and
 $E(G(\X,1)) = \emptyset$, we have  $A \subset \cup_{x \in \X} B_1(x)$
(else we could find a point $y\in A \setminus \cX$ 
%that could be added to $\X$
with $E(G(\X \cup \{y\},1)) = \emptyset$).
Hence by (\ref{psidef}) and (\ref{opsidef}), $\okappa \leq \oalpha$. 
%This provides another way
%to prove Corollary \ref{corocov}.

Since the domination number of a graph with no edges
equals the number of vertices,  with our current choice of $\zeta$
we have $\alpha^*(A) \leq \zeta^*(A) $. But also 
by (\ref{graphineqs})
%since $\gamma(G) \leq \alpha(G)$ for any finite graph $G$,
we have $\zeta(\cX) \leq \alpha(\X)$ for any finite
$\X \subset \R^d$, and hence
$\zeta^*(A) \leq \alpha^*(A)$.
% since since for any
%finite $\X \subset A$ and $\cY \subset \X$ with
%$\cY$ independent in $G(\X,1)$ and
%$|\cY| = \zeta(\X)$, we have $E(G(\cY,1))= \emptyset$ so
%that $\alpha^*(A) \geq |\cY|$.
 Thus $\alpha^*(A) = \zeta^*(A)$. 
 Hence by
(\ref{wellknown2}) and (\ref{opsidef}),  $\oalpha = \ozeta$.
%
%
%By Lemma~\ref{lemdet}, in particular~\eqref{0922a}, we have 
%$ \lambda \rho(\lambda) \leq \overline{\zeta}$ for all $\lambda > 0$. 
%Taking the limit $\lambda \to \infty$,
% and applying Theorem~\ref{th:onemore},
%gives us the result.
%% $\lim_{\lambda \to \infty} \lambda \rho(\lambda) =\kappa^*_B$. 
 \end{proof}

The inequality in  Proposition \ref{corocov} is strict,
at least for low dimensions.
Indeed, for $d=1$ it can be seen that $\okappa=1/2$, while $\oalpha=1$. 
For $d=2$,  $\okappa =  \sqrt{4/27}$ 
 while $\oalpha = \sqrt{4/3}$ 
as mentioned already.
 For $d=3$ we have $\oalpha = \sqrt{2}$
as mentioned already,  while 
%$\okappa \leq 1/2$
%by consideration of a face-centred cubic array and 
$\okappa \leq \sqrt{125/1024}$ by consideration
of a body-centred cubic array; see for example \cite{Bambah, Few}.
Thus, at least in low dimensions we have
$\okappa < \ozeta$, so that the limit in Theorems
 \ref{th:onemore} and \ref{thsubcond}
 (which equals $\okappa$) is not equal to $\ozeta$. 
This distinguishes the domination number from a number of 
other graph parameters (see Subsections \ref{subsecindep},
 \ref{subsecCC} and \ref{subsecED})
 for which we 
{\em can} apply Theorems \ref{th:lambig} and \ref{th:densegen} to obtain
limiting results analogous to  Theorems \ref{th:onemore} and
\ref{thsubcond}, but there with the limit equal to the relevant
$\ozeta$.

\begin{lemm}
\label{lemdomsup}
 Let  $\cY \subseteq \R^d$
and
 $\cZ \subseteq \R^d$ be finite and disjoint.
Then $\gamma( \cY \cup \cZ) \geq \gamma (\cY ) + \gamma(\cZ) 
- (1+ \kappa(B_2(o)) ) | \partial_{\cZ}\cY|$.
\end{lemm}
\begin{proof}
Let $\cS \subset \cY \cup \cZ$ be a dominating set in $\cY \cup \cZ$
with $|\cS| = \gamma( \cY \cup \cZ)$. Enumerate the
points of $ \cY$ distant at most 1 from $\cZ$ as $y_1,\ldots,y_m$,
where $m= |\partial_{\cZ}\cY|$.
Then $(\cS \cap \cY)  \cup \{y_1,\ldots,y_m\} $ is a dominating set 
in $\cY$.

Now set $k:=\kappa(B_2(o))$, the number of closed balls of radius 1 required
to cover a closed ball of radius 2 centred on the origin. 
For $1 \leq i \leq m$, let $B_{i,1} ,\ldots B_{i,k}$ 
be closed balls of radius $1/2$ that cover $B_1(y_i)$. For $1 \leq j \leq k$,
if $\cZ \cap B_{i,j} \cap B_1(y_i) \neq \emptyset$ then 
pick one element of $\cZ  \cap B_{i,j} \cap B_1(y_i)$ and denote it $z_{i,j}$. Then 
for each 
$ z \in B_1(y_i) \cap \cZ$, there exists $j$ such that $z \in B_{i,j}$,
and then $z_{i,j}$ is defined and $\|z - z_{i,j} \| \leq 1$. Hence
the set
$$
{\cal T} :=
\cup_{i=1}^m
\{z_{i,j}:
j \in \{1,\ldots,k\},  \cZ \cap B_{i,j} \cap B_1(y_i) \neq \emptyset 
\} 
$$
is a dominating set for  $ \cZ \cap (\cup_{i=1}^m B_1(y_i))$ 
(which equals $\partial_\cY \cZ$), and hence
$
(\cS \cap \cZ) \cup {\cal T}
$
is a dominating set for $\cZ$. Also $|{\cal T} | \leq mk $.  Then
\bean
\gamma(\cY) + \gamma(\cZ) \leq |(\cS \cap \cY) \cup \{y_1,\ldots,y_m\}| 
+ | (\cS \cap \cZ) \cup {\cal T}|
\\
\leq (|\cS \cap \cY | + m ) + (|\cS \cap \cZ| + km)
\\
= |\cS| + (1+k) m  
\\
= \gamma(\cY \cup \cZ) + (1+k) |\partial_{\cZ} \cY|,
\eean
which gives us the result. 
\end{proof}

Although Proposition \ref{propsofrho}(iv) does not apply with
the current choice of $\zeta$ (because P6 fails),
we conjecture that $\lambda \rho(\lambda)$ is
 nondecreasing in $\lambda$ for all values of $\lambda > 0$.
By 
Proposition \ref{propsofrho}(vii) and Lemma \ref{lemdomsup}, we
know at least
  that $\lambda \mapsto \lambda \rho(\lambda)$ is increasing on
 $(0,1/(\pi_d (1+ \kappa(B_2(o))))$.

\subsection{Clique-covering number}
\label{subsecCC}
For finite $\cX \subset \R^d$, let $\theta(\cX):=\theta(G(\cX, 1))$,
the clique-covering number of $G(\X,1)$.
\begin{theo}[LLNs for clique-covering number of RGG]
\label{lemclicov}
Theorems~\ref{th:therm2}, \ref{th:lambig} and \ref{th:densegen} all apply
 when choosing $\zeta(\cX):= \theta(\cX)$.
% to be the clique-covering number $\theta(G(\cX, 1))$.
\end{theo}
Denote by $\otheta$ the quantity $\ozeta$ appearing in Theorems
\ref{th:lambig} and \ref{th:densegen} when we take
 $\zeta(\cdot)= \theta(\cdot)$.
Then $\otheta$ is the minimum density (number of sets per unit volume)
of a partition of $\R^d$ into sets of Euclidean diameter at most 1.
It can be seen that $\otheta =1$ for $d=1$. For $d=2$, by
 partitioning the plane into regular hexagons of Euclidean diameter 1,
it can be shown that $\otheta \leq \sqrt{64/27}$.
We conjecture that in fact equality holds here, i.e. the most efficient
packing of the plane by sets of Euclidean diameter 1 uses regular
hexagons. For $d=3$, an upper bound for $\otheta$
can be obtained for example by tiling the three-dimensional space into 
 trapezo-rhombic doedecahedra of diameter 1.
\begin{proof}[Proof of Theorem \ref{lemclicov}]
Given  finite disjoint $\cY, \cZÂ \subset \mathbb{R}^d$,
we have
 $\theta(\cY \cup \cZ) \le \theta(\cY)+\theta(\cZ)$: indeed, colour
 all vertices in $\cY$ with $\theta(\cY)$ many colours, and all vertices
 in $\cZ \setminus \cY$ with 
 $\theta(\cZ \setminus \cY) \le \theta(\cZ)$
fresh
 colours. The colouring obtained is clearly an upper bound on $\theta(\cY \cup \cZ)$, and Property P3 (with $c_1=0$) follows. 

 For Property P4, let 
$\pi$ be
a clique-partition of $\cY \cup \cZ$ with $|\pi | = \theta (\cY \cup \cZ)$. 
Set 
$$
\pi_1 := \{S \cap \cY: S \in \pi, S \cap \cY \neq \emptyset\}; ~~~~~~
\pi_2 := \{S \cap \cZ: S \in \pi, S \cap \cZ \neq \emptyset\}.
$$
Then $\pi_1, \pi_2$ are clique-partitions of $\cY, \cZ$ respectively, with
$$
|\pi_1| + |\pi_2| - |\pi| = |\{S \in \pi: S \cap \cY \neq \emptyset
~{\rm and} ~
 S \cap \cZ \neq \emptyset\} |
\leq | \partial_{\cZ} \cY|,
$$
and rearranging this shows that P4
 holds with $c_2=1$. Also P$5'$ holds, since 
 $\theta(\X) \leq 1$
if
$\X \subset B_{1/2}(o)$ (so that $G(\X,1)$ is complete).
Therefore Theorem \ref{th:therm2} is applicable.

Clearly P6 and P7 also hold, since the clique-covering number can only be increased by adding vertices and can only be decreased by adding edges.
Since P6 holds, Theorems \ref{th:lambig} and \ref{th:densegen} also apply.
\end{proof}

\subsection{Eternal domination number}
\label{subsecED}
For finite $\cX \subset \R^d$, set 
$\gamma^\infty(\cX) := \gamma^\infty(G(\X,1))$,
 the eternal domination number.
% of $G(\X,1)$.
\begin{theo}[LLNs for eternal domination number of RGG]
\label{th:ED}
Theorems~\ref{th:therm2}, \ref{th:lambig} and \ref{th:densegen}
all  apply when choosing $\zeta(\cX) := \gamma^{\infty}(\cX)$.
 %$\gamma^{\infty}(G(\cX, 1))$.
\end{theo}

Denote by $\ogamma$ the quantity $\ozeta$ appearing in Theorems
\ref{th:lambig} and \ref{th:densegen}
when we take  $\zeta(\X) = \gamma^\infty(\X)$.
Then $\ogamma$ is the minimum density (number of guards per unit volume
in the large-$s$ limit)
of a  finite set of guards placed in $Q_s$
 that can defend against any sequence of attacks on locations in $Q_s$
with moves of Euclidean size at most 1 in each step.
By (\ref{graphineqs}) and
%  $\alpha(G) \leq \gamma^\infty(G) \leq \theta(G)$, we have
%from
 (\ref{wellknown2}) we have
 $\oalpha \leq \ogamma \leq \otheta$ in all dimensions.
Hence $\ogamma =1$ for $d=1$, and for $d=2$ we have
$\sqrt{4/3} \leq \ogamma \leq \sqrt{64/27}$.  It would be of interest
to find sharper upper and lower bounds when $d=2$.

\begin{proof}[Proof of Theorem \ref{th:ED}]
%As before, for $\cX \subset \R^d$, let
% $\gamma^{\infty}(\cX)=\gamma^{\infty}(G(\cX, 1))$. 
Given finite disjoint $\cY, \cZ \subset \mathbb{R}^d$, if
$\cA \subset \cY$ is eternally dominating in $\cY$, and
$\cA' \subset \cZ$ is eternally dominating in $\cZ$, then
$\cA \cup \cA'$ is eternally dominating in $\cY \cup \cZ$, since
one can defend all attacks on  $\cY$ using guards from $\cA$ and
 all attacks on  $\cZ$ using guards from $\cA'$.
Therefore
 $\gamma^{\infty}(\cY \cup \cZ) \le \gamma^{\infty}(\cY)+\gamma^{\infty}(\cZ)$,
so
 %indeed, consider a sequence of attacks to vertices of $\cY \cup \cZ$: if a vertex in $\cY$ is attacked, then using  $\gamma^{\infty}(\cY)$ guards,
% there exists a strategy to protect such vertices forever, and the set remains dominating, and similarly if a vertex in $\cZ$ is attacked, and 
Property P3 
%\ref{propsub}
 holds  with $c_1=0$.
We check P4 in Lemma \ref{lemEDP4} below.

Property P$5'$ (local uniform boundedness) follows from the fact that
$\gamma^\infty(\X) \leq \theta(\X)$ by (\ref{graphineqs}),
and the uniform boundedness of $\theta(\cdot)$ checked earlier.
Property P6 holds by Lemma \ref{lem:ED2} below.
 Therefore Theorems \ref{th:therm2},
 \ref{th:lambig} and \ref{th:densegen} all apply.
\end{proof}

Before checking P4 we give a further lemma.
 Given a finite graph $G= (V,E) $, consider a function 
$\phi: V \to \N \cup \{0\}$ representing an  assignment of finitely many
guards to the vertices,  now with
  with {\em more than one guard allowed} on each
vertex. 
Let us say that $\phi$ is 
eternally dominating if starting from  $\phi$, one can defend against
any sequence of attacks on unoccupied vertices,
 moving one guard at a time.
In the special case where $\phi$ takes values only in $\{0,1\}$ this reduces
to our earlier definition of eternal domination.
   Let $\hatgamma^\infty(G)$ be the minimum number of guards (i.e.,
the minimal value of $\sum_{v \in V} \phi(v)$) required
for an eternally dominating $\phi$.   

\begin{lemm}
\label{EDlem}
Let $G = (V,E)$ be a finite graph. Then $\hatgamma^\infty(G) = \gamma^\infty(G)$. 
\end{lemm}
\begin{proof}
Essentially this is shown in~\cite{Burger}, but for completeness
we include a more
detailed explanation than is given there.
% it is shown that $\gamma^{\infty}(G)$ does not depend on
% whether at most one guard is allowed to be located at any vertex at any time or not.
%% Then
%
The inequality $\hatgamma^\infty(G) \leq \gamma^\infty(G)$ is obvious.

Let $\phi: V  \to  \N \cup \{0\}$ be eternally dominating
with $\sum_{v \in V} \phi(v)  = \hatgamma^\infty(G)$,
and also with maximum total {\em coverage}
out of all such functions,
where the  coverage of $\phi$  is defined to be the total number of
 $v \in V$ for which $\phi(v) > 0$.
Suppose for some  $ v_0  \in V$ that $ \phi(v_0) > 1 $.
We shall derive a contradiction.

For any finite sequence of attacks,  the assignment of guards
after defending against these attacks must still be eternally dominating.
Also the coverage never goes down as one defends against the successive
attacks, because one moves a guard to an empty vertex each time.

The coverage never goes up as one defends against successive
attacks, because we assumed the original $\phi$ had maximum coverage
out of all eternally dominating assignments of guards.
Therefore in defending against any sequence of attacks, we never move
any guard off from $v_0$, because if we did then the coverage would go up.

But this shows we could still defend ourselves if we reduced $\phi(v_0)$ to 1
(leaving $\phi(v)$ unchanged for all other $v \in V$) and never moving
 the (single)
guard at $v_0$. This change reduces the total number of guards, contradicting
the earlier assertion that  originally $\sum_{v \in V} \phi(V) = 
\hatgamma^\infty(G)$.
Thus by contradiction, $\phi(v) \leq 1 $ for all $v \in V$, so
$\gamma^\infty(G) \leq \sum_{v \in V} \phi(v) = \hatgamma^\infty(G)$.
\end{proof}

\begin{lemm}
\label{lemEDP4}
 Property P4 (superadditivity up to boundary) holds for $\gamma^\infty(\cdot)$.
\end{lemm}
\begin{proof}
let  $\cY, \cZ $ be disjoint finite subsets of $\R^d$, and
let $\cA$ be an eternally dominating set in $\cY \cup \cZ$ with
$|\cA| = \gamma^\infty(\cY \cup \cZ)$.

Define $\phi: \cY \to \N \cup \{0\}$ by
$\phi:= {\bf 1}_{\cA \cap \cY} + {\bf 1}_{\partial_{\cZ}\cY}$.
%Let $\cA'= \cA'_1 \cup \cA'_2$
That is, $\phi$ is an assignment of guards to $\cY$ obtained by
adding to the guards of $\cA \cap \cY$,
 a guard at each vertex in
$\partial_{\cZ} \cY$ (in addition to any guards that were already there).

Then $\phi $ is eternally dominating for $\cY$. Indeed, since
$\cA$ is eternally dominating for $\cY \cup \cZ$, the guards  of
$\cA$ can defend against any sequence of attacks  on 
vertices of $\cY \setminus \partial_{\cZ} \cY$, and for any such
sequence of attacks the guards on $\cA \cap \cZ$ are unable
to help with the defence, so therefore the guards of $\cA \cap \cY$
are able to defend against any such sequence.  
On the other hand, for any sequence of attacks on vertices  in
$\partial_{\cZ}\cY$, the added guards in $\partial_{\cZ} \cY$ are able to
defend without moving at all. So the combined guards of
 ${\bf 1}_{\A \cap \cY} + {\bf 1}_{\partial_{\cZ}\cY}$ 
can defend against any sequence of attacks on $\cY$.
Together with Lemma \ref{EDlem},
this shows that 
\bea
\gamma^\infty( \cY) \leq |\cA \cap \cY|
+ |\partial_{\cZ} \cY| .
\label{0902a}
\eea

For each $y \in \partial_\cZ \cY$, let $\pi_y$ be a clique-partition
of $\cZ \cap B_1(y)$, of size $k_y$ with $k_y \leq k := \kappa(B_2(o))$.
This can be found by a similar argument to the one in the proof of 
Lemma \ref{lemdomsup}.
Enumerate the sets of $\phi_y$ as $\cS_{y,1},\ldots,\cS_{y,k_y}$.
For each $j=1,\ldots,k_y,$  pick an element $z_{y,j} $ of $\cS_{y,j}$.
Then define $\phi: \cZ \to \N \cup \{0\}$ by 
$$
\phi' = {\bf 1}_{\cA \cap \cZ} + \sum_{y \in \partial_{\cZ} \cY } 
\sum_{j=1}^{k_y} {\bf 1}_{z_{y,j}}.
$$
In other words we take the set of guards in $\cA \cap \cZ$, and
for each $y \in \partial_\cZ \cY$, and each set $\cS_{y,j}$ in the
 clique-partition  $\pi_y$, we add one guard at a
vertex  in $\cS_{y,j}$.

Then we claim $\phi'$ is eternally dominating for 
$\cZ$. Indeed, 
similarly to before, the guards of $\cA \cap \cZ$
 can defend against any sequence of attacks on vertices in 
$\cZ \setminus \partial_{\cY} \cZ$. 
%
%Let $\cA''_1 $ be the set of guards in $\cA \cap \cZ$, and let
%$\cA''_2$ be a set of guards in $\partial_{\cY} \cZ$ obtained
%as follows:
%for each $y \in \partial_\cZ \cY$, and for each set $\cS$ in the
% clique-partition  $\pi_y$
% pick an element  of $\cS$ and place a guard there.  
%Then, similarly to before, the guards of $\cA''_1$ can defend
%against any sequence of attacks on vertices in 
%$\cZ \setminus \partial_{\cY} \cZ$.
 Also for each $y \in \partial_{\cZ} \cY$, and each
set $\cS_{y,j}$ of the partition $\pi_y$, the added guard placed
at  $z_{y,j}$ can defend against  
 any sequence of attacks on vertices in
$\cS_{y,j}$,
%$\partial_{\cY} \cZ$, since for each $y \in \partial_{\cZ} \cY$,
%and for each set $\cS$ in $\pi_y$,
since the subgraph induced by $\cS_{y,j}$ is complete.
%the added guard placed in $\cS$ can be used to defend against all attacks
%within $\cS$. 

%Therefore 
%$\cA'' := \cA''_1 \cup \cA''_2$ can defend against

%any sequence of attacks on $\cZ$. 
Since $k_y \leq k$ for each $y \in \partial_\cZ \cY$,  we have
 $\sum_{z \in \cZ}
\phi'(z) \leq |\cA \cap \cZ| + k  |\partial_{\cZ} \cY|$. 
Hence, by Lemma \ref{EDlem},  
 $\gamma^\infty( \cZ) \leq |\cA \cap \cZ|
+ k |\partial_{\cZ} \cY|$. Combining this with (\ref{0902a})
yields
$$
\gamma^\infty(\cY) + \gamma^\infty(\cZ) \leq |\cA| + (1+k) |\partial_{\cZ} \cY|
= \gamma^\infty(\cY \cup \cZ) + (1+k) | \partial_\cZ \cY|,
$$
and rearranging this gives us Property P4 with $c_2 = 1+ \kappa(B_2(o))$.
\end{proof}

\begin{lemm}
\label{lem:ED2}
Property P6 (upward continuity in $\cX$) holds for $\gamma^\infty (\cdot)$.
\end{lemm}
\begin{proof}
Given finite $\cX \subset \R^d$ and 
 $x \in \R^d \setminus \cX$, let $\cA \subset \cX
\cup \{x\}$ be 
 eternally dominating for $\cX \cup \{x\}$ (i.e., for $G(\X \cup \{x\},1)$)
 with $|\cA| = \gamma^\infty(\cX \cup \{x\})$. 

If $x \notin \cA$, then starting from $\cA$ we can defend against
any sequence of attacks on vertices in $\cX$, so $\cA$ is
eternally dominating for $\cX$, and hence 
% and each move we use in doing so is still allowed if 
 $\gamma^\infty(\cX) \leq |\cA| .$
%= \gamma^\infty(\cX \cup \{x\})$.

Suppose $x \in \cA$. If any sequence of attacks on vertices in $\cX$
can be defended without moving the guard at $x$, then $\cA \setminus \{x\}$
is eternally dominating for $\cX$, so that $\gamma^\infty(\cX) \leq |\cA| -1$.
Otherwise, there exists some finite sequence of attacks on vertices in $\cX$, 
such that the eternally dominating defence against these attacks ends with
 moving the guard at $x$. After defending against  this sequence of attacks we
are left with a configuration
 $\cA' \subset \cX$  that is  
 eternally dominating (for $\cX \cup \{x\}$)
and satisfies $|\cA'| = |\cA|$. Then, as for the first case considered above,
$\cA'$ is eternally dominating for $\cX$, so $\gamma^\infty(\cX)
\leq |\cA'| = |\cA|$. 

Thus in all cases we have $\gamma^\infty(\cX) \leq |\cA| 
= \gamma^\infty(\cX \cup \{x\})$, which gives us P6.
\end{proof}

\subsection{Vertex cover number}
\begin{theo}[LLN for vertex cover number of RGG]
Corollary~\ref{CorollaryMain2} (with $c_3=1$) applies when choosing 
$\zeta(\cX)$ to be the vertex cover number of $G(\cX, 1)$.
\end{theo}
\begin{proof}
Let $\zeta(\X)$ here denote the vertex cover number of $G(\X,1)$.
Suppose $\cY$ and $\cZ$ are disjoint finite subsets of $\R^d$.
If $\cD$ is a vertex cover of $G(\cY \cup \cZ,1)$ 
with $|\cD|= \zeta(\cY \cup \cZ)$, then 
$\cD \cap \cY$
and $\cD \cap \cZ$
are vertex covers of $G(\cY,1) $ and of $G(\cZ,1)$ respectively. 
Therefore $\zeta(\cY) + \zeta(\cZ) \leq \zeta(\cY \cup \cZ)$.

Let $\cA$, $\cA'$ be vertex covers  of $G(\cY,1) $ and  $G(\cZ,1)$
respectively, with $|\cA|= \zeta(\cY)$ and $|\cA'|= \zeta(\cZ)$.
Then the  set $\cA \cup \cA' \cup \partial_{\cZ}(\cY)$ is a vertex
cover for $G(\cY \cup \cZ,1)$. Indeed, the set $\cA$ covers all 
edges of $G(\cY \cup \cZ,1)$ having both endpoints in $\cY$, 
while $\cA'$ covers
all edges of $G(\cY \cup \cZ,1)$ having both endpoints in $\cZ$, 
and $\partial_{\cZ}(\cY)$ covers all 
 edges of $G(\cY \cup \cZ,1)$ having one endpoint
in $ \cY $ and one endpoint  in $\cZ$. Therefore $\zeta(\cY \cup \cZ)
\leq \zeta(\cY) + \zeta(\cZ) + |\partial_{\cZ}(\cY)|$.

Thus we have Properties P3, P4 for $\zeta$ but with the inequalities
the wrong way round. Therefore taking $\zeta'(\X) :=  
|\X| - \zeta(\X)$ gives us a functional satisfying
Properties P3 and P4 with the inequalities the right way round
and with $c_1=0$, $c_2 = 1$. Moreover,
clearly $\zeta(\X) \leq |\cX|$ so $\zeta'(\cdot)$ is
 nonnegative. Also, if $G(\X,1)$ is a complete graph
then $\zeta(\X) = |\X| -1 $ and
 $\zeta'(\X) =1$, so 
$\zeta'(\cdot)$ satisfies P$5'$. Therefore, taking $\zeta'' (\cdot) \equiv 0$,
we can apply Corollary~\ref{CorollaryMain2} 
 to $\zeta$, with $c_3=1$. 
\end{proof}

\subsection{Number of connected components}
When we take $\zeta(\cX)$ to be the number of components
of $G(\X,1)$, and $\mu$ is absolutely continuous,
 a result along the lines of Theorem~\ref{th:therm2} was already known to be true (see \cite[Theorem 13.25]{RGG}). 
  Here we extend that result to cases where 
$\mu$ has a singular part.
\begin{theo}[LLN for number of components of RGG]
\label{th:numcon}
Theorem~\ref{th:therm2} applies when choosing 
$\zeta(\cX)$ to be the number of connected components of $G(\cX, 1)$.
\end{theo}
As we shall discuss in Section \ref{subsecMST},
this result can be obtained as
a special case of Theorem \ref{lem:MST1}
below. 
Alternatively, one can check Properties P1--P4 and P$5'$ directly.

In this case P6 does not apply.
However, P7 does hold here,
so by Proposition \ref{propsofrho}(v), 
 $\rho(\lambda)$ is nonincreasing in
$\lambda$. In fact by comparing (\ref{0101c}) with
 \cite[Theorem 13.25]{RGG}, one can see that 
$\rho(\lambda) = \E[ | \cC_o(\lambda)|^{-1}]$
in this case. Therefore, here we have $\lim_{\lambda \to \infty} \rho(\lambda)
=0$.

\subsection{Number of isolated vertices}
\label{subseciso}
Let $\sigma(G)$ denote the number of isolated vertices 
 of a graph $G$, and for finite $\X \subset \R^d$ set
$\sigma(\X):= \sigma(G(\cX,1))$. 

When $\mu$ is absolutely continuous, it was already known from
\cite[Theorem 3.15]{RGG} that
$n^{-1} \sigma(G(\X_n,r_n))$ converges in the thermodynamic limit. The following
result adds to that result by allowing $\mu$ to have a singular part.
\begin{theo}[LLN for number of isolated vertices of RGG]
\label{lemiso}
Theorem~\ref{th:therm2} applies when choosing $\zeta(\cX) $ to
be $\sigma (\cX)$. 
\end{theo}
\begin{proof}
Suppose $\cY,\cZ$ are finite disjoint subsets of $\R^d$. Then every isolated
point in $\cY \cup \cZ$ is  either an isolated point in $\cY$ or
an isolated point in $\cZ$, so $\sigma (\cY \cup \cZ) \leq
\sigma(\cY) + \sigma(\cZ)$; thus P3 holds with $c_1 =0$. 
Moreover, if $\cI$ denotes the set of isolated vertices 
in $\cY \cup \cZ$, then every isolated vertex in $\cY$ 
is either  in $\cI$ or lies within unit distance  of
$\cZ$, so that
\bea
 \sigma(\cY) \leq |\cI \cap \cY | + |\partial_{\cZ} \cY|.
\label{0924b}
\eea
Similarly, each isolated vertex in $\cZ$
is either in $\cI$ 
 or lies within unit distance of $\cY$.
Moreover, for each $y \in \Y$ the number of isolated vertices of $\cZ$
within unit distance of $y$ is bounded by  $3^d$,
since  this is an upper bound for the number of disjoint balls of
radius 1/2 which can be fitted inside a ball of radius 3/2.
Therefore the number of isolated vertices of $\cZ$
within unit distance of $\cY$ is bounded by $3^d |\partial_{\cZ} \cY|$. 
Hence,
$$
 \sigma(\cZ) \leq |\cI \cap \cZ| + 3^d |\partial_{\cZ} \cY|.
$$
Combined with (\ref{0924b}) this shows that
$$
\sigma(\cY \cup \cZ) = |\cI| \geq \sigma (\cY) + \sigma(\cZ)
- (1+ 3^d) |\partial_\cZ \cY|,
$$
so P4 holds with $c_2 = 1+ 3^d$. 
Also P$5'$ holds since 
$\sigma(\cX) \leq 1 $ 
 for
any $\cX \subset B_{1/2}(o)$. 
%is bounded
%by the maximum number of disjoint balls of radius $1/2$ that can be packed
% inside $Q_2$, which is finite. 
\end{proof}
Note that P7 holds here, and hence 
% holds as well, since the number of isolated vertices can only decrease
% by adding edges. Hence,
 Proposition \ref{propsofrho}(v) is applicable here.
However, P6 does not hold here.

\subsection{$H$-packing number}
Let $H$ be a fixed  connected graph  with $h$ 
 vertices, $2 \leq h < \infty$.
 Recall that the   $H$-packing number $\psi_H(G)$ is defined to be the largest
possible size of $H$-packing in $G$.
For finite $\X \subset \R^d$ set $\psi_H(\X):= \psi_H(G(\X,1))$. 

\begin{theo}[LLN for $H$-packing number of RGG]
	\label{lem:Packing}
Corollary~\ref{CorollaryMain2} applies,  with
$c_3= 1/h$, when choosing 
$\zeta(\cX)$ to be $\psi_H(\X)$.
\end{theo}
{\bf Remark.} 
A generalization of $H$-packing allows for more than
one $H$. Suppose we are given a collection
of pairs $\{(H_1,v_1),\ldots,(H_m,v_m)\}$, where for each 
$i \in \{1,\ldots,m\}$,
$H_i$ is a finite connected graph and $v_i > 0$,
with $H_1,\ldots,H_m$ pairwise non-isomorphic. We call $v_i$ the `value'
of $H_i$. Given a graph $G= (V,E)$,
 the aim now is to pack vertex-disjoint 
copies of the sets $H_1,\ldots,H_m$  
 into $G$ (allowing repetitions)
 with maximum possible total value. That is, a packing of
$G$ is a collection $\pi = \{J_i, i \in {\cal I}\}$ of  
vertex-disjoint subgraphs, each of which is isomorphic to
one of $H_1,\ldots,H_m$, and the total value $v(\pi)$  of the packing
equals $\sum_{i \in {\cal I}} w(J_i)$, where
we set $w(J_i)$ to be $v(H_j)$ if $J_i$ is isomorphic to $H_j$. Then
define $MP_{(H_1,v_1),\ldots,(H_m,v_m)}(G)$ to
be  the maximum of $v(\pi)$, over all packings
$\pi$.
 Theorem \ref{lem:Packing} can be generalized to this setting by
minor modifications of the proof.
 In this generality we need to take 
$c_3= \max_{i \in \{1,\ldots,m\}} (v_i/h_i)$, where  $h_i$ denotes
the number of vertices of $H_i$. 

\begin{proof}[Proof of Theorem \ref{lem:Packing}] 
Let $\cY,\cZ \subset \R^d$ be disjoint and finite. 
Then any  $H$-packing of $\cY$
and $H$-packing of $\cZ$ can be put together to give an $H$-packing
of $\cY \cup \cZ$; hence $\psi_H( \cY \cup \cZ) \geq \psi_H(\cY) + \psi_H(\cZ)$.

%Next, observe that 
For any $H$-packing of $\cY \cup \cZ$, if we remove
all $H$-subgraphs in the packing that contain
 at least one vertex in $\cY$ and at least
one vertex in $\cZ$, we are left with an $H$-packing of $\cY$ and
an $H$-packing of $\cZ$. The number of removed $H$-subgraphs is
at most $|\partial_{\cZ}\cY|$, and hence
$$
\psi_H(\cY \cup \cZ) - | \partial_{\cZ}\cY| \leq \psi_H(\cY) + \psi_H(\cZ).
$$
Thus, taking
 $\zeta(\cdot)  =  \psi_H(\cdot)$, Properties P3
and P4 hold with the inequalities the wrong way round.
Therefore,   the functional  $\zeta'$, defined by
$
\zeta'(\cX) := ( |\cX|/h) -  \zeta(\cX),
$
%which
%gives us a functional
 satisfies P3 and P4 with $c_1=0$, $c_2 =1$.
Moreover clearly  $\zeta' (\cX) \geq 0$ for all $\cX$.

If $G(\cX,1)$ is a complete graph, then 
$h \zeta(\X) \geq |\X | -h$ so $\zeta'(\cX) \leq 1$;
hence  P$5'$ holds for $\zeta'$.
 Thus, taking $\zeta''(\cdot) \equiv 0$, we can apply
 Corollary~\ref{CorollaryMain2}  with $c_3=1/h$.
\end{proof}

\subsection{Edge cover number}

Recall that  $\eta(G)$ denotes the edge cover
number of the graph obtained from $G$ by removing all isolated
vertices.
 For $\cX \subset \R^d$,  write  $\eta(\cX)$
for $\eta(G(\cX,1))$. 

\begin{theo}[LLN for edge cover number of RGG]
Corollary~\ref{CorollaryMain2} applies, with $c_3 = 1/2$, when choosing 
$\zeta(\cX) := \eta(\cX)$.
\end{theo}
\begin{proof}
For any finite graph $G=(V,E)$, let $\theta_2(G)$ denote the minimum
size of partition
 $\pi$
of $V$ such that
every set $W \in \pi$ is either a single vertex or two vertices
connected by an edge.
We shall refer to any such partition as an {\em edge-partition}
of $V$.
 For $\cX \subset \R^d$ set  $\theta_2(\X)
:= \theta_2(G(\cX,1))$.

Recalling from Subsection \ref{subseciso}
the notation $\sigma(G)$ 
 for the number
 of isolated vertices of a graph $G$, we assert that
 $\eta(G) = \theta_2(G) - \sigma(G)$ for any finite graph $G$
(and hence $\eta(\cX) = \theta_2(\cX) - \sigma(\cX)$ for
 any finite $\cX \subset \R^d$). It clearly suffices to check this  
assertion
in the case where $\sigma(G)=0$. In this case
$\eta(G) \leq \theta_2(G)$ since when $\sigma(G)=0$,
starting with any edge-partition of $V$ 
we can replace each singleton in the partition
by an edge incident to it, to get an edge cover of the same cardinality.
But also given any edge cover, we can think of each
edge as a set of two vertices and  go through these sets
one by one,
at each step
 removing any vertex that was already counted,
(i.e. possibly sometimes changing the pair of vertices to a singleton,
or even removing the pair altogether if the original edge cover
 was of non-minimum cardinality),
to end up with an edge-partition
of $V$ consisting of at most the same number of sets as there were
 edges in the original edge cover.
Hence $\theta_2(G) \leq \eta(G)$ so 
 $\theta_2(G) = \eta(G)$ as asserted.

Now let $H$ be a complete graph on 2 vertices, 
i.e. a graph with two vertices and one edge.
For all finite $\cX \subset \R^d$,
we claim that 
\bea
 \eta(\X)= |\X|- \psi_H(\X) - \sigma(\X).
\label{1230a}
\eea
This is because any minimum edge-partition of
$\cX$
consists of a maximum $H$-packing in $\cX$, together with
$|\cX| - 2\psi_H(\cX)$ singletons, making a total of $\psi_H(\cX) +
|\cX| - 2 \psi_H(\cX)$ sets in the partition.
Thus we have
$\theta_2(\cX) = |\cX| - \psi_H(\X)$, and 
(\ref{1230a}) follows.

Define $\zeta(\cdot):= \eta(\cdot)$,
$\zeta'(\cdot) = \sigma(\cdot)$ and  
 $\zeta''(\cdot)= (|\cdot|/2) -  \psi_H(\cdot)$.
 By Theorem \ref{lemiso}, $\zeta'(\cdot)$ satisfies
Properties P1--P5, 
 and by the proof of Theorem \ref{lem:Packing},
so does $\zeta''(\cdot)$.
By (\ref{1230a}), for all $\cX$
we have $\zeta(\cX) = (|\X|/2) - \zeta'(\X) + \zeta''(\X)$. 
We may thus apply Corollary~\ref{CorollaryMain2} with $c_3=1/2$, as asserted.
\end{proof}

\section{Applications to weighted complete graphs}\label{secweighted}
\allco

In this section we consider the limit theory of four classic combinatorial
optimization problems (travelling salesman, minimum matching,
 minimum bipartite matching, minimum spanning tree), defined
 on a weighted complete  graph on vertex set $\X_n$, 
where the edge weights are determined 
by the inter-point displacements via  a specified measurable 
{\em weight function} $w: \R^d \to [0,\infty)$.
Given $w (\cdot)$, define  
 $\wmax \in [0,\infty]$ by
$\wmax := \sup \{w(x): x \in \R^d\}$.
We shall require  $w(\cdot)$ to satisfy 
the following conditions (but allow it to be otherwise arbitrary).
\begin{enumerate}
\item[W1] ({\em symmetry}):  $w(x) = w(-x)$ for all $x \in \R^d$.
\item[W2]  ($w$ is small near the origin):
 $w(o)=0$, and
$w(x) \to 0$ as $x \to o$.
 \item[W3] ($w $  is large far from the origin):  $w(x) \to \wmax 
%\sup_{y \in \R^d} w(y)
$  
as 
 $\|x \| \to \infty $.
\item[W4] (Polynomial growth bound): There exist $c_4 \in (0,\infty) $ and
$p \in (0,d)$ such that $w(x) \leq c_4 \max(\|x\|^p,1)$ for all 
$x \in \R^d$.
\end{enumerate}

In some lemmas
%Sometimes
we shall require the following
stronger versions of W3, W4:
\begin{enumerate}
\item[W5] ($w$ achieves its maximum far from the origin). There is a
 constant $c_5 \in (0,\infty)$ such that $w(x) = \wmax < \infty$ for all
$x \in \R^d$  with $\|x\| > c_5$.
\item[W6] There exists $p \in (0,d)$, $c_6 \in (0,\infty)$ such that
$w(x) \leq c_6 \|x\|^p$ for all $x \in \R^d$.

\end{enumerate}

Given $a \in [0,\wmax)$, define the truncated weight function
$w_a(x) = \min(w(x),a) $, $x \in \R^d$. Note that if $w(\cdot)$ satisfies
W3, then $w_a(\cdot)$ satisfies W5.
Moreover, if $w(\cdot)$ satisfies W4, then, for any $\delta >0$,
$w(\cdot)(1- {\bf 1}_{B_\delta(o)}(\cdot))$ satisfies W6.
Also, if $\wmax < \infty$, then W4 follows automatically.

As usual we are given a sequence of scaling parameters $r_n$. 
We shall  take the complete graph on vertex set $\X_n$, with 
 each edge  $e=\{x,y\}$  having weight  $w_n(e)  := w(r_n^{-1}(y-x))$ (the unscaled weight function is $w(e) := w(y-x)$.)
One example of interest is to take 
$w(x) =1 -  {\bf 1}_{B_1(o)}(x)$, so that $w_n(e) $ is zero if $e$ is
an edge of $G(\X_n,r_n)$, and otherwise is 1.

Our results in this section extend results given in~\cite{Yukich},
where attention is restricted to
weight functions of the form $w(x) = \|x\|^\alpha$.
Here we allow  for a much more general class of weight functions $w(\cdot)$; for example, if $d \geq 2$ one could take
%Here are some examples that are not covered in \cite{Yukich}:
%
$w(x) = \log ( |x| +1 )$, or 
$w(x) = |x| ( 2 + \sin |x| )$, or
$w(x) =  |x|^{3/2} + |x|^{1/2}$, using some arbitrary
norm $|\cdot|$.
	Also we consider the bipartite matching problem, which is
not addressed in \cite{Yukich}.

Recall that for $(\lambda,s)\in (0,\infty)^2$,
   we write $\H_{\lambda,s} $ for $\H_\lambda \cap Q_s$.
   We say that $\mu$ has {\em bounded support} if 
 $\mu(Q_s)= 1 $ for some $s \in (0,\infty)$.
   The next two lemmas are messy to state, but will be repeatedly
   useful for proving LLNs for  certain functionals via
   approximating functionals. 
\begin{lemm}
	\label{lemlong}
	%Let $p \in [0,d)$.
	Suppose that for each $k \in \N $, 
	$\zeta^{(k)}(\X)$ is a measurable nonnegative real-valued
	functional defined on all finite $\X \subset \R^d$.
	Assume for each $\X$ that $\zeta^{(k)}(\X)$ is nondecreasing in $k$,
	and set $\zeta(\X) := \lim_{k \to \infty} \zeta^{(k)}(\X)$.
	Suppose for each $k \in \N, \lambda >0$ that the limit 
	$\rho_k(\lambda)
	:= \lim_{s \to \infty} \E[ \zeta^{(k)}(\H_{\lambda,s})/(\lambda
	s^d)]$ exists and is finite, and that if $n r_n^d \to t \in (0,\infty)$
	as $n \to \infty$, then 
%	$\zeta_a(r_n^{-1} \X_n) \toccL \int \rho_a (t f_\mu(x)) f_\mu(x) dx$.
\bea
	\zeta^{(k)}(r_n^{-1} \X_n) \toccL 
\int \rho_{k} (t f_\mu(x)) f_\mu( x) d x. 
\label{0316e2}
\eea
	Let $(h(k), k \in \N)$ be a real-valued sequence
	with $h(k) \to 0$ as $k \to \infty$. Suppose moreover that either (a)  
	$\mu$ has bounded support and there exist constants
	$p \in (0,d), c_7 \in (0,\infty)$ such that
	for all $s\geq 1$, $k \in \N$ and all finite $\X \subset Q_s$,  we have 
	\bea
	\label{0321b}
	%\zeta^{(k)}(\X) & \leq & c_7 s^p |\X|^{(d-p)/d};
	\zeta^{(k)}(\X) & \leq & c_7 s^p |\X|^{1-(p/d)};
\\
	\label{0321c}
	\zeta(\X) - \zeta^{(k)}(\X) & \leq
	%& s^p h(k) (\zeta^{(k)}(\X))^{(d-p)/d},
	& s^p h(k) (\zeta^{(k)}(\X))^{1-(p/d)},
	\eea
	or (b) 
	$\zeta(\X) - \zeta^{(k)}(\X) \leq h(k)|\X|$
	for all finite
	$\X \subset \R^d$, $k \in \N$.  
	%we have
%
	Then for each $ \lambda >0$ the limit $\rho(\lambda)
	:= \lim_{s \to \infty} \E[ \zeta(\H_{\lambda,s})/(\lambda
	s^d)]$ exists and is finite, and if $n r_n^d \to t \in (0,\infty)$
	as $n \to \infty$, then $\zeta(r_n^{-1} \X_n) \toccL \int
	\rho (t f_\mu(x)) f_\mu(x) dx$, i.e. (\ref{0101c}) holds.
\end{lemm}

\begin{proof}
	First we assume condition (a) holds.  Using
(\ref{0321b}), taking expectations
and using Jensen's inequality yields 
for all $k,\lambda,s >0$ that
$$
	\E[\zeta^{(k)}(\H_{\lambda,s})] \leq  c_7  s^p \E [ |\H_{\lambda,s}|^{(d-p)/d} ] 
	\leq   c_7s^p (\lambda s^d)^{1-(p/d)} =  c_7\lambda^{1-(p/d)} s^d.
$$
Hence using (\ref{0321c}) and Jensen's inequality again yields
	that
\bea
	\E[\zeta(\H_{\lambda,s}) - \zeta^{(k)}(\H_{\lambda,s})  ] & \leq & 
	 s^p h(k) ( \E[\zeta^{(k)} (\H_{\lambda,s})])^{1-(p/d)} 
\nonumber \\
	& \leq & 
	c_7^{1-(p/d)} s^p h(k) \lambda^{(1-(p/d))^2} s^{d-p}
\nonumber \\
	& = & 
	c_7^{1-(p/d)}  h(k) \lambda^{(1-(p/d))^2} s^{d}.
\label{0316c2}
\eea

	Let $\lambda >0$.
	By (\ref{0316c2}), given $\eps >0$ we can choose 
	$k_0 \in \N$ such that %for all for all $s >0$ we have 
\bean
	0 \leq \E [ \zeta(\H_{\lambda,s} ) -  \zeta^{(k)}(\H_{\lambda,s} ) ]/(\lambda
s^d) \leq \eps, ~~~ \forall~ s \geq 1, k \geq k_0.
\eean 
	Note that this is also true under condition (b).
	Taking the large-$s$ limit yields for $k \geq k_0$ that
$$
	\rho_k(\lambda) \leq \liminf_{s \to \infty}
	\E[\zeta(\H_{\lambda,s}) /(\lambda s^d)]  
	\leq \limsup_{s \to \infty} \E[\zeta(\H_{\lambda,s}) /(\lambda s^d)]  
\leq \rho_k(\lambda) 
+ \eps.
$$
Therefore since $\rho_k(\lambda)$ is nondecreasing in $k$, the limit
$\rho(\lambda) := \lim_{k \to \infty} \rho_k(\lambda) $ exists, is
finite, and is equal to
 $\lim_{s \to \infty} \E[\zeta(\H_{\lambda,s}) /(\lambda s^d)]$
	under either condition (a) or (b).
By monotone convergence,
\bea
\lim_{k \to \infty} \int \rho_{k} (t f_\mu(x)) f_\mu( x) d x 
=  \int \rho (t f_\mu(x)) f_\mu( x) d x.
\label{frommon2}
\eea 

	Under condition (a), by the assumption of bounded support we
can and do choose $C = C(\mu) \in (0,\infty)$ such that $r_n^{-1} \X_n \subset
	Q_{r_n^{-1} C}$ almost surely, for all $n \in \N$.
 Then  by
 (\ref{0321b}), %and the 
	%(depending on $\mu$) such that 
	$\zeta^{(k)}(r_n^{-1} \X_n) \leq c_7 C^p r_n^{-p} n^{1-(p/d)}$,
 almost surely.
Using (\ref{0321c}), therefore we have, almost surely,

\bean
\zeta(r_n^{-1}\X_n) - \zeta^{(k)}(r_n^{-1}\X_n) 
\leq C^p  r_n^{-p} h(k)  (c_7C^p r_n^{-p} n^{1-(p/d)} )^{1-(p/d)}. 
\eean
Since we are taking the thermodynamic limit $nr_n^d \to t$,
 we have $r_n^p n^{p/d} \to t^{p/d}$, so there are constants $C', C''$
 (depending only on $\mu$ and $t$) such that, almost surely,
\bean
\zeta(r_n^{-1}\X_n) - \zeta^{(k)}(r_n^{-1}\X_n)  \leq
C' r_n^{-p} h(k) n^{1-(p/d)} 
\sim C'' h(k)  n.
\eean
Hence given $\eps >0$, we can choose $k_0, n_0 \in \N$ such that for all
$k \geq k_0, n \geq n_0$,
almost surely
$ \zeta(r_n^{-1}\X_n) - \zeta^{(k)}(r_n^{-1}\X_n)
\leq \eps n ;$ this is also true under condition (b).
Hence using (\ref{0316e2})
and (\ref{frommon2}) we can obtain the
desired result (\ref{0101c}) here.
\end{proof}

\begin{lemm}
\label{lemlong2}
	Suppose $\mu$ has bounded support.
Suppose for each $\ell,m \in \N$ that
	$\tzeta^{(\ell,m)}(\X)$ is a measurable 
non-negative real-valued functional defined for all  
all finite  $\X \subset \R^d$,
	%with $\zeta(\emptyset)=0$, 
	such that for all
	$\lambda >0$ the limit
	$\rho^{(\ell,m)} (\lambda) := \lim_{s \to \infty}
	\E[\tzeta^{(\ell,m)} (\H_{\lambda,s})/(\lambda s^d)] $ 
	exists and is finite, and if $n r_n^d \to t \in (0,\infty)$ as
	$n \to \infty$, then $n^{-1} \tzeta^{(\ell,m)}(r_n^{-1} \X_n)
	\toccL $ $\int \rho^{(\ell,m)}( t f_\mu(x)) f_\mu(x)dx$.
	Assume for all $\X$ that
	$\zeta^{(\ell,m)}(\X)$ is nondecreasing both in $\ell$ and
	in $m$, and set $
	\tzeta^{(\ell)}(\X) := \lim_{ m\to \infty}
	\tzeta^{(\ell,m)}(\X)  
	$
	%as $m \uparrow \infty$, for all $\ell \in \N$, and
	and
	$
	\tzeta(\X) := \lim_{\ell \to \infty} \tzeta^{(\ell)}(\X) 
	$.

	Let $p \in (0,d)$.
	Suppose for each $\ell \in \N$ that
	there is a constant $C(\ell)$ 
	and a real-valued sequence $(h_\ell(m), m \in \N)$
	with $h_\ell(m) \to 0$ as $m \to \infty$,
	such that for all $s \geq 1$
	and all finite $\X \subset Q_s$, we have $\tzeta^{(\ell,m)}(\X)
	\leq C(\ell) s^p |\X|^{1-(p/d)}$ 
	and
	$$
	\tzeta^{(\ell)}(\X) - \tzeta^{(\ell,m)}(\X) \leq
	h_{\ell}(m) s^p  (\tzeta^{(\ell,m)}(\X))^{1-(p/d)}.
	$$
	Suppose moreover that there is a sequence $(\tilh(\ell),\ell \in \N)$
	such that $\tilh(\ell) \to 0$ as $\ell \to \infty$, and
	such that $\tzeta(\X) - \tzeta^{(\ell)}(\X) \leq \tilh(\ell) |\X| $
	for all finite $\X \subset \R^d$.
	Then for all $\lambda >0$ the limit
	$\rho (\lambda) := \lim_{s \to \infty}
	\E[\tzeta (\H_{\lambda,s})/(\lambda s^d)] $ 
	exists and is finite, and if $n r_n^d \to t \in (0,\infty)$ as
	$n \to \infty$, then $n^{-1} \tzeta(r_n^{-1} \X_n)
	\toccL \int \rho( t f_\mu(x)) f_\mu(x)dx$.
	\end{lemm}

\begin{proof}
	Let $\ell \in \N$. Then taking 
	$\zeta(\cdot):= \tzeta^{(\ell)}(\cdot)$ 
	and $\zeta^{(k)}(\cdot):= \tzeta^{(\ell,k)}(\cdot)$ 
	gives us functionals which 
	together 
	satisfy (\ref{0321b}) and (\ref{0321c}) if we take 
	$C= C(\ell)$ and $h(k): = h_\ell(k)$.
	%Hence 
	By
	Lemma \ref{lemlong} (a),
for all $\lambda >0$
the limit 
$
	\rho^{(\ell)}(\lambda) := \lim_{s \to \infty}
	\E[\tzeta^{(\ell)} (\H_{\lambda,s})/ (\lambda s^d)]
$
exists and is finite, and if $nr_n^d \to t \in (0,\infty)$ 
	%as $n \to \infty$, 
then
$
	\tzeta^{(\ell)}( r_n^{-1} \X_n ) \toccL 
	%$ $
	\int \rho^{(\ell)}
	(t f_\mu(x)) f_\mu(x)dx.
$

	Now taking instead
	$\zeta(\cdot) := \tzeta(\cdot)$ and $\zeta^{(k)}(\cdot) :=
	\tzeta^{(k)}(\cdot)$ 
	 for $k \in \N$,
%	these functionals together satisfy the hypothesis of
	 we can apply
	Lemma \ref{lemlong} (b), now with $\tilh(k)$ playing the
	role of $h(k)$, to 
	obtain the desired conclusion.
\end{proof}

\subsection{Weighted travelling salesman problem}
Given finite $\cX \subset \R^d$,
 the (weighted) travelling
 salesman problem (TSP) is to find a tour of minimum
 weight (using weight function $w(\cdot)$)
 that passes exactly once through all points and returning to the starting
 vertex.
Formally, a   {\em tour} of $\X$ is defined to be a Hamilton cycle in the
complete graph  on vertex set $\X$. For each tour
$\tau$ of $\X$, 
define $w(\tau) = \sum_{e \in E(\tau)}
w(e)$, where $E(\tau)$ denotes the set of edges of $\tau$.
Define $TSP_w(\cX) := \min w(\tau) $, where the minimum is over
all  tours of $\X$.

\begin{theo}[LLN for weighted TSP]
	\label{lem:TSP}
	Suppose  $w(\cdot)$ satisfies conditions W1--W4,
	and either (a) $\wmax < \infty$ or (b) $\mu$ has bounded support,
	$\wmax = \infty$.
	  Then 
when choosing $\zeta(\cX)$ as $TSP_w(\cX)$, the limit $\rho(\lambda)$ given by
(\ref{meanconv}) exists and is finite for all $\lambda \in (0,\infty)$, 
and if $n r_n^d \to t \in (0,\infty)$ as $n \to \infty$,
	then (\ref{0101c}) holds.
% Theorem~\ref{th:therm2}  applies.
\end{theo}

%When $\wmax = +\infty$, we consider only the case
%where $\mu$ has bounded support.
%We also require the polynomial growth bound W4.
%
%\begin{theo}\label{lem:TSPunbound}
%Suppose  $w(\cdot)$ satisfies conditions W1--W4, with
%	$\wmax = + \infty$.  Suppose $\mu$ has
%bounded support.
%Then when
% choosing $\zeta(\cX)$ as $TSP_w(\cX)$, the limit $\rho(\lambda)$ given by
%(\ref{meanconv}) exists and is finite for all $\lambda \in (0,\infty)$, 
%and if $n r_n^d \to t \in (0,\infty)$ as $n \to \infty$,  
%then (\ref{0101c}) holds.
%\end{theo}
% Theorem~\ref{th:therm2}  applies.

{\bf Remark.} One could also consider a more general {\em directed}
version of the TSP, dropping condition W1 so that
the weight of an edge depends on the direction of travel.
Theorem~\ref{lem:TSP}(a) and its proof apply without change
to this directed  setting 
(we do not
know if the directed version of Theorem~\ref{lem:TSP}(b)
holds).
 For example,
in $d=2$ one could take $w(\cdot) = 1- {\bf 1}_S$ with $S$ being
the right half of $B_1(o)$ (a semi-circle); in this case some
minor modifications to the proof are needed because W2 fails too, but
the directed version of Theorem \ref{lem:TSP} still holds.
A different spatial directed MST problem  
has been considered in \cite{Steele}.

%The next lemma provides the first step towards the preceding theorems.
We first prove Theorem \ref{lem:TSP} under the extra condition W5.
\begin{lemm}
\label{LemTSPW5}
Suppose $w(\cdot)$ satisfies W1, W2 and W5. 
Then when choosing 
$\zeta(\cX)$ as $TSP_w(\cX)$, the limit $\rho(\lambda)$ given by
(\ref{meanconv}) exists and is finite for all $\lambda \in (0,\infty)$, 
and if $n r_n^d \to t \in (0,\infty)$ as $n \to \infty$, then (\ref{0101c}) holds.
\end{lemm}

\begin{proof}
Assume first that W5 holds with $c_5 =1$ (as well as
W1 and W2). 
We will show that Theorem \ref{th:therm2}  holds for $\zeta(\cdot) :=
 TSP_{w}(\cdot)$ in this case.
Clearly, $TSP_w(\cdot)$ satisfies Properties P1 and P2.
 Next,
 for any disjoint finite $\cY, \cZ \subset \R^d$, 
 let $\tau'$ be a tour of $\cY$ and $\tau''$ a tour of $\cZ$ with
$w(\tau')= MST_w(\cY)$ and 
$w(\tau'')= MST_w(\cZ)$.
 Pick an edge of $\tau'$, and denote its endpoints $y_1,y_2$. 
 Pick an edge of $\tau''$, and denote its endpoints $z_1,z_2$. 
Then $(E(\tau') \cup E(\tau'')
 \cup \{y_1z_1, y_2 z_2\}) \setminus \{y_1 y_2, z_1z_2\}$
is the edge-set of a tour of $\cY \cup \cZ$, and hence
$$
TSP_w(\cY \cup \cZ) \le TSP_w(\cY) +TSP_w(\cZ)+2 \wmax.
$$
%
%	since the optimal solution on $\cY \cup \cZ$ is clearly bounded
% above by an optimal solution on $\cY$, followed by an edge from the last
% point visited in $\cY$ (right before returning to the starting vertex in
% $\cY$) to $\cZ$ of weight at most $\wmax$, an optimal solution on $\cZ$ and
% another edge of weight at most $\wmax$ from the last point visited in $\cZ$ 
%(before returning to the first vertex in $\cZ$) returning to the starting 
%vertex in $\cY$.
 Thus, Property P3 holds with $c_1=2 \wmax$.

For Property P4,
 let
% $\cY \subset \R^d$ and
%$
% \cZ \subset \R^d$ be disjoint finite subsets of $\R^d$.
%Let
 $\tau$ be
% an optimal tour of $\cY \cup \cZ$, i.e.
 a tour of $\cY \cup \cZ$ with $w(\tau) =
TSP_w(\cY \cup \cZ)$.
Suppose without loss of generality
 that the tour $\tau$ starts at a vertex in $\cY$.
 Let $\tau_1$ be the tour through $\cY$ obtained by taking the vertices 
of $\cY$ in the same order as they arise in $\tau$,
 while leaving out from
 $\tau$ all vertices belonging to $\cZ$, 
thereby introducing some new edges within $\cY$. The number
of new edges equals half the number of edges of $\tau$ that go
from $\cY$ to $\cZ$.
Similarly let $\tau_2$ be the tour though $\cZ$  obtained from 
$\tau$ by going through the vertices of $\cZ$ in the order they
	arise in $\tau$.
	%
	%starting at the first vertex of $\cZ$ visited,
 %leaving out all vertices belonging to $\cY$, and going back
 %to this first vertex of $\cZ$.

	Now writing just
	$\tau $ for the set of edges of $\tau$ (i.e., identifying
a tour with its set of edges), we have
 $|\tau_1 \cup \tau_2| = |\cY | + |\cZ| =
|\tau|$, and hence
\bea
|\tau \setminus (\tau_1 \cup \tau_2)| =
|(\tau_1 \cup \tau_2) \setminus \tau|.
\label{0830c}
\eea
 Also, 
%By definition
$  TSP_w(\cY) \leq w(\tau_1)$ and $ TSP_w(\cZ) \leq w(\tau_2) $.
Hence
\bea
TSP_w(\cY) + TSP_w(\cZ) - TSP_w(\cY \cup \cZ) \leq
w(\tau_1) + w(\tau_2) - w(\tau)
\nonumber \\
	= \left( \sum_{e \in (\tau_1 \cup \tau_2) \setminus \tau } w(e) 
	\right)
- \sum_{e \in \tau\setminus (\tau_1 \cup \tau_2)    } w(e) 
\nonumber \\
\leq  \sum_{e \in \tau \setminus (\tau_1 \cup \tau_2) }
(\wmax- w(e) ),
\label{0830b}
\eea
where for the second inequality we have used both (\ref{0830c})
and the assumption that all edge weights are at most $\wmax$.

Let us say that an edge $uv$ of $\tau$ is {\em short}
if $\|u-v \| \leq 1$, {\em long} if $\|u-v\|>1$.
The last sum in (\ref{0830b}) is at most $\wmax$ times
 the number of {\em short} edges 
in $\tau$ that go from $\cY$ to $\cZ$ (because all long edges
have weight $\wmax$ so do not contribute
to the sum), and 
this is bounded above by $2 \wmax |\partial_{\cZ}\cY|$.
Therefore
 (\ref{0830b}) gives us
$$
TSP_w(\cY) + TSP_w(\cZ) - TSP_w(\cY \cup \cZ) 
\leq 2 \wmax | \partial_{\cZ} \cY|,
$$
and hence Property P4 with $c_2=2 \wmax$.

To check Property P5 (local sublinear growth),
 note that given any $\delta >0, n \in \N$
and  finite $\X \subset B_\delta(o)$,
for {\em any} tour $\tau$ of
$\X$
 we have 
   $w(\tau) \leq |\X| \sup_{x \in B_{2\delta}(o)} w(x)$. Therefore
$\sup\{ TSP_w(\X)/|\X|: \X \subset B_{\delta}(o) \}
\leq   
\sup_{x \in B_{2\delta}(o)} w(x)$, which tends to 0 as $\delta \downarrow 0$
by W2, and P5 follows. 

Thus,
% Lemma \ref{propsofrho} and
Theorem \ref{th:therm2} is applicable for
$\zeta(\cdot) := TSP_w(\cdot)$,
in the special case where $c_5=1$ in W5. 
In the general case where  $c_5 \in (0,\infty)$,
consider the
rescaled weight function $w'(x) :=  w(c_5 x)$, so that
$w'(x) = \wmax $ for $\|x\| > 1$.
By the special case already considered,
%Lemma \ref{propsofrho}
Theorem \ref{th:therm2}
 applies to $\zeta (\cdot) = TSP_{w'}(\cdot)$.  Thus
 for all $\lambda>0$ 
 %   writing $\H_{\lambda,s} $ for $\H_\lambda \cap Q_s$, we have that
the limit $\rho'(\lambda) := \lim_{s \to \infty}
\E[ TSP_{w'}(\H_{\lambda,s}) / (\lambda s^d) ] $ 
exists. Since $TSP_{w}(\X) =  TSP_{w'}(c_5^{-1} \X)$ for all $\cX$,
by the Mapping theorem
for Poisson processes \cite{LP}, 
\bean
\E[ TSP_w(\H_{\lambda,s} )]  /(\lambda s^d) 
=  \E[ TSP_{w'} (c_5^{-1} \H_{\lambda,s} )]  / (\lambda s^d) 
\\
=  
\E[ TSP_{w'} (\H_{c_5^d \lambda,c_5^{-1} s} )]  /( \lambda s^d), 
\eean
which converges to $ \rho'(c_5^d \lambda) =:
 \rho(\lambda)$ as $s \to \infty$.

Now suppose $n r_n^d \to t \in (0,\infty)$.
Then as $n \to \infty$ we have $n (c_5 r_n)^d \to c_5^d t$, so
by Theorem \ref{th:therm2} applied to $\zeta(\cdot) = TSP_{w'}(\cdot)$,
\bean
n^{-1} TSP_w(r_n^{-1} \X_n) =
n^{-1}  TSP_{w'}((c_5 r_n)^{-1}\X_n) 
\toccL   \int \rho'(c_5^dt f_{\mu}(x) ) f_{\mu}
(x) dx
\\
= \int \rho( t   f_\mu( x) )  f_\mu( x) dx, 
\eean
completing the proof in the general case.
\end{proof}

\begin{proof}[Proof of Theorem \ref{lem:TSP}(a)] 
Assume  that $w(\cdot)$ satisfies W1--W3 with $0 < \wmax < \infty$
(and hence also W4).
	Then
% since $w(\cdot)$ satisfies W3,
 the function $w_a (\cdot):= \min(w(\cdot),a)$
%	has a maximum value of
%$a$, and $w_a(x) = a$ for all $x$ sufficiently far from $o$; that is,
%$w_a$
	satisfies condition W5. Therefore we can apply Lemma \ref{LemTSPW5}
to the weight function $w_a$. In particular, for all $\lambda \in (0,\infty)$
 the limit $\rho_a(\lambda) := \lim_{s \to \infty}
\E[TSP_{w_a}(\H_{\lambda,s})/(\lambda s^d)]$ exists, and if
	$n r_n^d \to t \in (0,\infty)$ then (\ref{0101c}) holds
	for $\zeta(\cdot) = TSP_{w_a}(\cdot)$, i.e.
\bea
n^{-1} TSP_{w_a} (r_n^{-1} \X_n) \toccL \int \rho_a( t f_\mu(x)) f_\mu(x) dx.
\label{0319b}
\eea

%It is easy to see that $\rho_a(\lambda)$ is nondecreasing in $a$,
%and so the limit $\rho(\lambda) := \lim_{a \uparrow \wmax} \rho_a(\lambda)$
%exists, and that $\rho(\lambda) \leq \wmax$.
%Moreover, 
	For any finite $\X \subset \R^d$ and $a \in (0,\wmax)$
 we have $0 \leq TSP_{w}(\X) - TSP_{w_a}(\X) \leq (\wmax -a)|\X|$,
since for each edge $e$ of any tour of $\X$, we have
 $0 \leq w(e)- w_a(e) \leq \wmax - a$.
	Therefore setting  $a(k):= \max(\wmax - 1/k,0)$ for each
	$k \in \N$ 
	%with $a(k) \uparrow \wmax$ as $k \to \infty$,
	and setting
	$\zeta^{(k)}(\cdot) := TSP_{w_{a(k)}}(\cdot)$, gives us a
	sequence of functionals satisfying the conditions
	of Lemma \ref{lemlong}(b).    By that result we obtain
	that
	$\E[TSP_w(\H_{\lambda,s} ) /(\lambda s^d)]$
	converges to $\rho(\lambda)$ as $s \to \infty$
	for all $\lambda \in(0,\infty)$, and 
	if $n r_n^d \to t \in (0,\infty)$, then
	(\ref{0319b}) holds with $w_a(\cdot)$ replaced by $w(\cdot)$, $\rho_a$
	replaced by $\rho$,
as required.
\end{proof}

Now we aim to prove
%consider the case with $\wmax = \infty$, and
%work towards proving 
Theorem \ref{lem:TSP}(b), so we assume 
W1--W4, with $\wmax = +\infty$.
Initially we assume W6 too.
Again set $w_a(\cdot) := \min(w(\cdot),a)$.

\begin{lemm} \label{lemTSPdet}
	Suppose that  $w(\cdot)$ satisfies W1--W4 and W6
	with $\wmax = +\infty$.
	There is a constant $C$ such that
	for all $(a,s) \in (0,\infty)^2$,
	for all finite $\X \subset Q_s$ we have
\bea
TSP_{w_a}(\X) \leq 
	%TSP_{w} ( \X ) \leq
	C s^p |\X|^{(d-p)/d} , 
%~~~~ \X \subset Q_s, \X ~{\rm finite},
\label{0316a3}
\eea  
and moreover
\bea
	TSP_w(\X) - TSP_{w_a}(\X) \leq C s^p a^{(p/d) - 1} 
	(TSP_{w_a} (\X))^{(d-p)/d}.
	%~~~~ \X \subset Q_s, \X ~{\rm finite},
	\label{0321a}
\eea
\end{lemm}
\begin{proof}
By \cite[eqn (3.6) or (3.7)]{Yukich},
there is a constant $C \in (0,\infty)$ such that for any 
 finite $\cY \subset Q_1$, we have 
$TSP_{c_6 \|\cdot\|^p}(\cY) \leq C |\cY |^{(d-p)/d}$, where $c_6$ is the constant in assumption W6. 
 Hence  
	for any $s >0$,
	by
	%by assumption W6 and
	scaling
	%and any finite $\cY \subset Q_1$,
$ TSP_w(s \cY ) \leq TSP_{c_6\|\cdot\|^p} (s\cY) \leq C s^p |\cY|^{(d-p)/d} $.
	Since $TSP_{w_a}(\X) \leq TSP_w(\X)$ for all $\X$, this 
implies
%Hence for all $s \in (0,\infty)$,
(\ref{0316a3}).

Given $s, a >0$ and given finite $\X \subset Q_s$,
let $\tau_a(\X)$ be (the set of edges in)  a
minimum-weight tour  of $\X$, using weight function $w_a$ (and using some
	deterministic rule to choose if there exist several tours of
minimum weight). That is, let $\tau_a(\X)$ be a  tour of $\X$
with $w_a(\tau_a(\X)) = TSP_{w_a}(\X)$. Then  
let $\tau_{(=a)} (\X) = \{e \in \tau_a(\X): w_{a}(e) =a\}$, 
and  let $N_a (\X) := |\tau_{(=a)}(\X)|$.
% and let $W^{(a)} (\X) := \sum_{e \in \tau_{(a)}(\X)} w(e) $.
 Thus $N_a(\X) $ is the  total number of edges
of weight $a$ in the minimum $w_a$-weight tour of
$\X$ (in the top figure of Figure~\ref{FigureTSP} all edges of weight $a$ are drawn as thick edges).
 Let $\cY_a(\X)$ be the set of those vertices in $\X$ that
are incident to at least one edge in $\tau_{(=a)}(\X)$,
	and set $r= |\Y_a(\X)|$.   Note that
	$ N_a(\X) \leq r \leq 2N_a(\X)$.

We now remove the edges of $\tau_{(=a)}(\X)$ from $\tau_a(\X)$,
and add further edges to make a new tour of $\X$, for
which we can bound the total weight of the added edges.
To do this, we use the {\em space-filling curve heuristic}
(see \cite{Yukich}).
Let $\phi:[0,1] \to [-1/2,1/2]^d$ be
a continuous surjection. Assume moreover
that $\phi$ is  
Lipschitz of order $1/d$, that is, for any $t,t' \in [0,1]$, we have 
$\|\phi(t')-\phi(t)\| \le c_8 
	|t'-t|^{1/d}$ for some absolute constant $c_8 > 0$. 
Such a $\phi$  is well known to exist; see for example~\cite{Buckley}.
Consider a tour of $s^{-1} \cY_a(\X)$,
 obtained by visiting the points in order of this curve. Let 
$y_1, \ldots, y_{r}$ be the corresponding  points of $\cY_a(\X)$
taken in this order, and let $t_1 < \ldots < t_{r} \in [0,1]$
 be such that $s \phi(t_i)=y_i$; to be definite, if $\phi^{-1} (\{s^{-1} y_i\})$
has more than one element, let $t_i$ be the minimal element. 

We now create a multigraph
	on vertex set $[r] := \{1,\ldots, r\}$ as follows. 
Starting with the empty graph on this vertex set,
add  an edge from $i$ to $j$ for each $\{i,j\} \subset [r]$
 such that there exists a path from $y_i$ to $y_j$
	in the graph with vertex set $\X$ and edge set
	$\tau_a \setminus \tau_{(=a)}$, i.e. a path
	from $y_i$ to $y_j$ in the original tour $\tau_a$
	with all edge weights less than $a$ along this path.
	%$\{y_i,y_j\}$ is an edge of  
%$\tau_{(=a)}(\X)$, and
We shall refer to the  the edges added so far as \emph{red} edges.
At this stage, each vertex has degree $0$ or $1$, and the number
of red edges equals $r-N_a(\X)$.
% add `blue arcs' from $i$ to $i+1$, for each $i \in \{ 1,ldots,r-1\}$. 
%
%Also, letting $y_\alpha$ be the first
%of the $y_i$'s to appear in the original tour $\tau$,
%and $y_\omega$ to be the last, add a red  edge from $0$ to $\alpha$
%and a red edge from $\omega$ to $r+1$. 
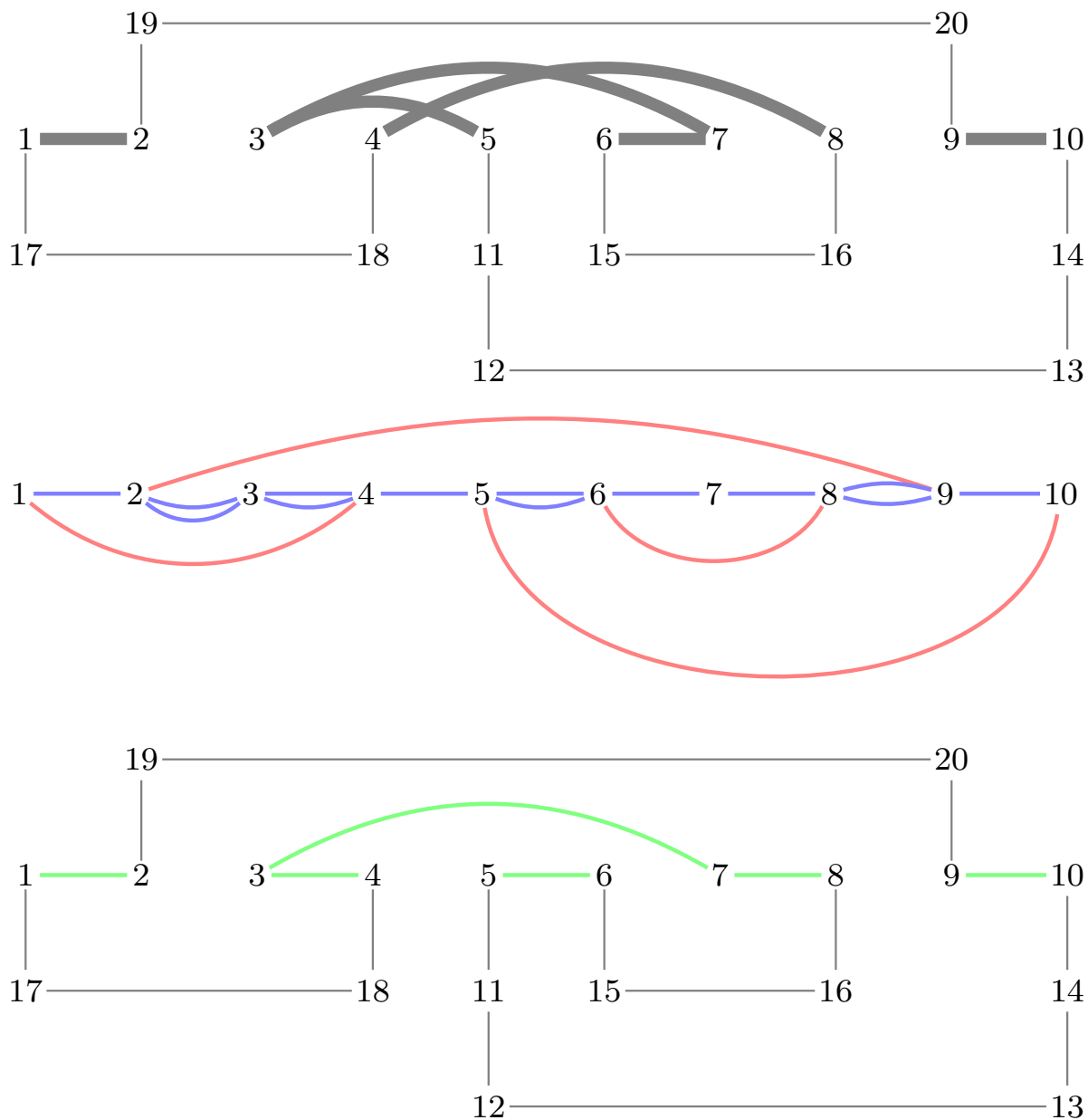
\begin{figure}
 \centering
 \scalebox{1.7}
 {
\begin{tikzpicture}
%[line cap=round,line join=round,x=1cm,y=1cm, scale = 0.5]
%\clip(-2.2,-8.5) rectangle (31,3.1);
\begin{scriptsize}
	 \tikzstyle{vertex}=[circle, minimum size=0.1pt, inner sep=0pt]
	 \tikzstyle{selected vertex} = [vertex, fill=red!24]
	 \tikzstyle{selected edge} = [draw,line width=1pt,-,blue!50]
	 \tikzstyle{selected edge2} = [draw,line width=3pt,-,black!50]
	 \tikzstyle{edge} = [draw,line width=0.5pt,-,black!50]
	 \node[vertex] (v1) at (1,3) {$1$};
	 \node[vertex] (v2) at (2,3) {$2$};
	 \node[vertex] (v3) at (3,3) {$3$};
	 \node[vertex] (v4) at (4,3) {$4$};
	 \node[vertex] (v5) at (5,3) {$5$};
	 \node[vertex] (v6) at (6,3) {$6$};
	 \node[vertex] (v7) at (7,3) {$7$};
	 \node[vertex] (v8) at (8,3) {$8$};
	 \node[vertex] (v9) at (9,3) {$9$};
	 \node[vertex] (v10) at (10,3) {$10$};
	 \node[vertex] (v11) at (5,2) {$11$};
	 \node[vertex] (v12) at (5,1) {$12$};
	 \node[vertex] (v13) at (10,1) {$13$};
	 \node[vertex] (v14) at (10,2) {$14$};
	 \node[vertex] (v15) at (6,2) {$15$};
	 \node[vertex] (v16) at (8,2) {$16$};
	 \node[vertex] (v17) at (1,2) {$17$};
	 \node[vertex] (v18) at (4,2) {$18$};
	 \node[vertex] (v19) at (2,4) {$19$};
	 \node[vertex] (v20) at (9,4) {$20$};
	 \draw[selected edge2] (v1) -- (v2);
	 \draw[selected edge2] (v9) -- (v10);
	   \draw[selected edge2] (v3) to [out=30,in=150] (v5);
	 	 \draw[selected edge2] (v3) to [out=30, in=150] (v7);
	 \draw[selected edge2] (v4) to [out=30, in=150] (v8);
	    \draw[selected edge2] (v6) -- (v7);
	  \draw[edge] (v1) -- (v17) -- (v18) -- (v4);
	  \draw[edge] (v5) -- (v11) -- (v12) -- (v13) -- (v14) -- (v10);
	  \draw[edge] (v6) -- (v15) -- (v16) -- (v8);
	  \draw[edge] (v2) -- (v19) -- (v20) -- (v9);
	  \end{scriptsize}
 \end{tikzpicture}
 }
 \scalebox{1.7}
 {
 \begin{tikzpicture}
 \begin{scriptsize}
	 \tikzstyle{vertex}=[circle, minimum size=0.1pt, inner sep=0pt]
	 \tikzstyle{selected vertex} = [vertex, fill=red!24]
	 \tikzstyle{blue edge} = [draw,line width=1pt,-,blue!50]
	  \tikzstyle{red edge} = [draw,line width=1pt,-,red!50]
	 \tikzstyle{selected edge2} = [draw,line width=3pt,-,black!50]
	 \tikzstyle{edge} = [draw,line width=0.5pt,-,black!50]
	 \node[vertex] (v1) at (1,3+8) {$1$};
	 \node[vertex] (v2) at (2,3+8) {$2$};
	 \node[vertex] (v3) at (3,3+8) {$3$};
	 \node[vertex] (v4) at (4,3+8) {$4$};
	 \node[vertex] (v5) at (5,3+8) {$5$};
	 \node[vertex] (v6) at (6,3+8) {$6$};
	 \node[vertex] (v7) at (7,3+8) {$7$};
	 \node[vertex] (v8) at (8,3+8) {$8$};
	 \node[vertex] (v9) at (9,3+8) {$9$};
	 \node[vertex] (v10) at (10,3+8) {$10$};
	  \draw[red edge] (v1) to [out=320,in=220] (v4);
	   \draw[red edge] (v6) to [out=300,in=240] (v8);
	   \draw[red edge] (v5) to [out=280,in=260] (v10);
	    \draw[red edge] (v2) to [out=18,in=162] (v9);
	   \draw[blue edge] (v1) -- (v2);
	    \draw[blue edge] (v2) to [out=320,in=220] (v3);
	      \draw[blue edge] (v2) to [out=340,in=200] (v3);
	       \draw[blue edge] (v3) -- (v4);
	       \draw[blue edge] (v3) to  [out=340,in=200] (v4);
	        \draw[blue edge] (v4) -- (v5);
	            \draw[blue edge] (v5) -- (v6);
	       \draw[blue edge] (v5) to  [out=340,in=200] (v6);
	      \draw[blue edge] (v7) -- (v6);
	      \draw[blue edge] (v8) -- (v7);
	         \draw[blue edge] (v8) to [out=-15,in=195] (v9);
	      \draw[blue edge] (v8) to [out=15,in=165] (v9);
	       \draw[blue edge] (v9)--(v10);
	  \end{scriptsize}
 \end{tikzpicture}
 }
   \scalebox{1.7}
 {
\begin{tikzpicture}
%[line cap=round,line join=round,x=1cm,y=1cm, scale = 0.5]
%\clip(-2.2,-8.5) rectangle (31,3.1);
\begin{scriptsize}
	 \tikzstyle{vertex}=[circle, minimum size=0.1pt, inner sep=0pt]
	 \tikzstyle{selected vertex} = [vertex, fill=red!24]
	 \tikzstyle{selected edge} = [draw,line width=1pt,-,blue!50]
	 \tikzstyle{green edge} = [draw,line width=1pt,-,green!50]
	 \tikzstyle{edge} = [draw,line width=0.5pt,-,black!50]
	 \node[vertex] (v1) at (1,3+13) {$1$};
	 \node[vertex] (v2) at (2,3+13) {$2$};
	 \node[vertex] (v3) at (3,3+13) {$3$};
	 \node[vertex] (v4) at (4,3+13) {$4$};
	 \node[vertex] (v5) at (5,3+13) {$5$};
	 \node[vertex] (v6) at (6,3+13) {$6$};
	 \node[vertex] (v7) at (7,3+13) {$7$};
	 \node[vertex] (v8) at (8,3+13) {$8$};
	 \node[vertex] (v9) at (9,3+13) {$9$};
	 \node[vertex] (v10) at (10,3+13) {$10$};
	 \node[vertex] (v11) at (5,2+13) {$11$};
	 \node[vertex] (v12) at (5,1+13) {$12$};
	 \node[vertex] (v13) at (10,1+13) {$13$};
	 \node[vertex] (v14) at (10,2+13) {$14$};
	 \node[vertex] (v15) at (6,2+13) {$15$};
	 \node[vertex] (v16) at (8,2+13) {$16$};
	 \node[vertex] (v17) at (1,2+13) {$17$};
	 \node[vertex] (v18) at (4,2+13) {$18$};
	 \node[vertex] (v19) at (2,4+13) {$19$};
	 \node[vertex] (v20) at (9,4+13) {$20$};
	  \draw[edge] (v1) -- (v17) -- (v18) -- (v4);
	  \draw[edge] (v5) -- (v11) -- (v12) -- (v13) -- (v14) -- (v10);
	  \draw[edge] (v6) -- (v15) -- (v16) -- (v8);
	  \draw[edge] (v2) -- (v19) -- (v20) -- (v9);
	  \draw[green edge] (v1) -- (v2);
	    \draw[green edge] (v9) -- (v10);
	      \draw[green edge] (v5) -- (v6);
	        \draw[green edge] (v8) -- (v7);
	         \draw[green edge] (v3) -- (v4);
	         \draw[green edge] (v3) to  [out=30,in=150] (v7);
	  \end{scriptsize}
 \end{tikzpicture}
 }
 \caption{Top: Example of an optimal tour $\tau_a$ with respect to $w_a$ (on $|\X|=20$ vertices); thick edges have weight $a$ ($r=10$ here, the space-filling curve visits vertices $1,\ldots, 10$ in this order). Middle: construction of red and blue edges of $\{1,\ldots, 10\}$. Bottom: an Eulerian circuit visits vertices in this order:
 1,2,9,10,5,6,8,9,8,7,6,5,4,3,2,3,4,1. New list of vertices: 1,2,9,10,5,6,8,7,3,4,1. In green: not-red edges corresponding to the new list added to the new tour.}
 \label{FigureTSP}
 \end{figure}

Next, if $1$ has odd degree, add a \emph{blue} edge from $1$ to $2$;
otherwise add two blue edges from $1$ to $2$.
Then if $2$ has odd degree, add a blue edge from $2$ to $3$; otherwise,
add two blue edges from $2$ to $3$.
Then if $3$ has an odd degree, add a blue edge  from $3$
to $4$, and otherwise add two such edges, and so on.

Continuing in this way up to vertex $r-1$,
 we add one or two blue edges between $i$ and $i+1$
for each $i \in \{1,\ldots,r-1\}$ in turn (see the middle figure of
 Figure~\ref{FigureTSP} for an example).  We end up with
a multigraph which we denote $\G$, in which each
 $i \in \{1,\ldots, r-1\}$ has even degree, so
   vertex $r$ must also have even degree
since the sum of degrees is even. Thus,
$\G$ is connected and all of its degrees are even; in fact,
each vertex has degree $2$ or $4$. We define the {\em red-degree}
of a vertex in $\G$ to be the number of
red vertices incident to that vertex (which is either $0$ or $1$). 

By Euler's ``K\"onigsberg bridge'' theorem, we can
and do choose an 
Eulerian circuit  through $\G$, starting and ending
at vertex $1$. Choose a direction of travel through this circuit.
We denote this directed Eulerian circuit by $\gamma$. 
All  edges of $\G$ now become arcs (i.e., directed edges),
with the direction of an edge determined by the direction of travel 
through the edge in $\gamma$.

List the vertices and arcs of $\G$ in the order they appear in the circuit
$\gamma$. This list is an alternating sequence of vertices and arcs, starting and ending at vertex $1$. Vertices of degree $4$ on $\G$  appear twice in the list
(as does vertex $1$); other vertices appear once. Now reduce this list to a list of vertices only,
in which each vertex other than vertex $1$ appears once, as follows.

Suppose vertex $i >1$ appears twice in the old list. If $i$ has red-degree
$0$ then omit the second appearance of $i$ in the old list. If $i$ has
red-degree $1$, then retain the appearance of $i$ in the old list
which adjacent in that list to the red arc incident to $i$, 
and omit the other appearance of vertex $i$.

Write the vertices in the order of this new list as $i(1)=1, i(2), \ldots,
i(r), i(r+1)=1$. Let $\cK$ be the set of indices $k \in [r]$ such that
$(i(k), i(k+1)) $ is {\em not} a red arc in $\G$.

We now create our new tour of $\X$, by taking 
the  set of  edges 
$$
\{\{y_{i(k)},y_{i(k+1)}\},  k \in \cK \} \cup
(\tau \setminus \tau_{(=a)})
$$
(see the bottom figure of Figure~\ref{FigureTSP} for the edges added corresponding to the Eulerian circuit chosen in the text below the figure).
 In other words, we take the edges in the
 tour $(y_{i(1)}, y_{i(2)}, \ldots, y_{i(r)},
 y_{i(1)})$ of $\cY_a(\X)$, and then replace each red edge of
 this tour (i.e., a step $\{y_{i(k)}, y_{i(k+1)}\}$
 for which $\{i(k), i(k+1)\}$ is a red
 edge of $\G$) with  the portion of the original tour $\tau_{(a)}$
 that goes from $y_{i(k)} $ to $y_{i(k+1)}$ and consists entirely
 of edges of weight less than $a$.

The total $w$-weight of
 this new tour of $\X$ is an upper bound for $TSP_w(\X)$, so
$TSP_w(\X) - TSP_{w_a}(\X)$ is bounded above by the total
weight of edges in this tour that are not in $\tau$. 
%the added edges.
Thus using W6, we have
 \bea
 TSP_w(\X) - TSP_{w_a}(\X)  
 %&\le &
 \leq \sum_{k\in \cK} 
w(y_{i(k+1)}  - y_{i(k)}) 
  %\nonumber \\ &\le&
  \le c_6 \sum_{k \in \cK}
\|y_{i(k+1)}  -y_{i(k)}\|^p.  
\nonumber 
\eea
Using 
the Lipschitz  continuity of $\phi$, and then
H\"older's inequality, we obtain that
\bea
 TSP_w(\X) - TSP_{w_a}(\X)  
 %&=& 1+ C' \sum_{i=1}^{r-1} \| \phi(t_i) - \phi(t_{i-1}) \|^d \\
 %&\le& C'' \sum_{i=1}^{r-1} \| t_i-t_{i-1} \|^d \\
 &\le&   c_6c_8^p s^p \sum_{k\in \cK}
 | t_{i(k+1)} -t_{i(k)} |^{p/d}
\nonumber \\
 &\le&   c_6 c_8^p  s^p
 %(r-N_a(\X))
 |\cK|^{1-(p/d)}
 \sum_{k \in \cK} | t_{i(k+1)}  -t_{i(k)} |.
\label{0318a}
\eea

In the multigraph $\G$, each blue edge goes between two neighbouring
vertices in the list $\{1,\ldots,r\}$. For $1 \leq \ell \leq r-1$,
let each blue edge from $\ell $ to $\ell+1$
 be given weight $|t_{\ell+1} -t_\ell|$.
Then  for $ k \in \cK$, the value of
$|t_{i(k+1)} - t_{i(k)}|$ is bounded above by the
sum of the weights of the edges making up the portion
of the  Eulerian circuit  $\gamma$ that goes from $i(k)$
to $i(k+1)$ (all of which are blue).
Since this circuit is Eulerian, no blue edge is traversed more
than once, so the sum in the last line of (\ref{0318a})
is at most  the total weight of all the blue 
edges, which is at most 2 because there are at most two
blue edges from $i $ to $i+1$, $1 \leq i \leq r-1$. Also
%$|\cK| \leq r- N_a(\X) \leq N_a(\X)$, 
$|\cK| = N_a(\X)$, 
so
\bean
TSP_w(\X) - TSP_{w_a}(\X) \leq
2 c_6 c_8^p s^p (  N_a(\X))^{1-(p/d)},   ~~~~~ \X \subset Q_s, \X {\rm ~
finite.} 
\eean
Since $a N_a(\X) \leq TSP_{w_a}(\X)$, this yields 
(\ref{0321a}).
\end{proof}

%\begin{proof}[Proof of Theorem \ref{lem:TSPunbound}]
	\begin{proof}[Proof of Theorem \ref{lem:TSP}(b)]
		Assume that $w(\cdot)$ satisfies W1--W4 with
		$\wmax= +\infty$,
		and that $\mu$ has bounded support.
Given $\delta > 0$,
	define the new weight function $w'_\delta(x) := w(x)
	(1- {\bf 1}_{B_\delta(o)}(x))$.
	For $m \in \N $ set $w'_{\delta,m}(x) := \min(w_\delta(x),m)$.

	Then $w'_{\delta,m}(\cdot)$ satisfies W2 and W5. Hence by 
	Lemma \ref{LemTSPW5}, for all $\lambda >0$ the limit
	$\rho'_{\delta, m} (\lambda) := \lim_{s \to \infty}
	\E[TSP_{w'_{\delta,m}} (\H_{\lambda,s})/(\lambda s^d)] $ 
	exists, and if $n r_n^d \to t \in (0,\infty)$ as
	$n \to \infty$, then $n^{-1} TSP_{w'_{\delta,m}}(r_n^{-1} \X_n)
	\toccL \int \rho'_{\delta, m}( t f_\mu(x)) f_\mu(x)dx$.

	Since we assume $w(\cdot)$ satisfies W4, for each
	$\delta >0 $ the function 
	$w'_\delta(\cdot)$ satisfies W6. Therefore by Lemma
	\ref{lemTSPdet}, 
	there is a constant $C(\delta)$ such that for all $s >0$
	and all finite $\X \in Q_s$, we have $TSP_{w'_{\delta,m}}(\X)
	\leq C(\delta) s^p |\X|^{(d-p)/d}$ and 
	\bea
	TSP_{w'_\delta}(\X) - TSP_{w'_{\delta,m}}(\X) \leq
	C(\delta) s^p m^{(p/d)-1} (TSP_{w'_{\delta,m}}(\X))^{(d-p)/d}.
	\label{0407a}
	\eea

	 Now for $\ell,m \in \N $,  
	set $\tzeta^{(\ell,m)}(\X)
	:= TSP_{w'_{1/\ell,m}}(\X)$, which is nondecreasing both in 
	$\ell$ and in $m$.
	Set $\tzeta^{(\ell)}(\X):=
	\lim_{m \to \infty} \tzeta^{(\ell,m)}(\X) 
	%$ is equal to $
	= TSP_{w'_{1/\ell}}(\X)$,
	and set  $\tzeta(\X):= \lim_{\ell \to \infty} \tzeta^{(\ell)}(\X)
	= TSP_w(\X)$.
	Set $\tilh(\ell):= \sup_{x  \in B_{1/\ell}(o)} w(x)$.
	Then $\tilh(\ell) \to 0$ as $ \ell \to \infty$,
	by  W2, and
	\bea
	%0 \leq
	\tzeta(\X) - \tzeta^{(\ell)}(\X) 
	= TSP_{w}(\X) - TSP_{w'_{1/\ell}}(\X) \leq  \tilh(\ell) |\X|.
	\label{0407b}
	\eea
	By (\ref{0407a}) and (\ref{0407b}) we can 
	  apply Lemma \ref{lemlong2} to the functionals
	 $\tzeta^{(\ell,m)}(\cdot)$, $\ell,m \in \N$,
	 %taking $h_{\ell}(m) = C(1/\ell)m^{-(d-p)/d}$. Hence
	 taking $h_{\ell}(m) = C(1/\ell)m^{(p/d)-1}$. Hence
	the limit
	$ 
	\lim_{s \to \infty}
\E[TSP_{w}(\H_{\lambda,s})/(\lambda s^d) ]
=:	\rho(\lambda)
$
	exists and is finite, 
	for all $\lambda >0$, 
	and if $n r_n^d \to t \in (0, \infty) $ as $n \to \infty$,  then
$
TSP_{w}( r_n^{-1} \X_n ) \toccL \int \rho(t f_\mu(x)) f_\mu(x)dx,
$
as required. 
\end{proof}

\subsection{Minimum-weight matching}
 Given finite  $\cX \subset \R^d$,
 %$\cX = \{X_1, \ldots, X_n\}$, 
the minimum-weight matching (MM) problem is to find a near-perfect matching 
(that is, a perfect matching if $|\X|$ is even and a matching excluding one 
vertex if $|\X|$ is odd) of minimum weight. More precisely, 
we define a (near-perfect) {\em matching} of $\cX$ to be a collection $\tau$
of $\lfloor |\cX|/2\rfloor$ pairwise disjoint edges of the complete graph on
vertex set $\cX$,
where two  edges $uv$ and $xy$ are said to be disjoint if $u,v,x,y$
are all distinct. Given such a matching $\tau$, we define $w(\tau):= \sum_{e \in
\tau} w(e)$. We define $MM_w(\X) := \min_\tau w(\tau)$, where the minimum
is taken over all matchings $\tau$ of $\cX$.  

\begin{theo}[LLN for minimum-weight matching]
		\label{lem:MM1}
%Suppose  $w(\cdot)$ satisfies conditions W1--W3, with $\wmax < \infty$.  
	Suppose  $w(\cdot)$ satisfies conditions W1--W4, and either (a)
	 $\wmax < \infty$ or (b) $\mu$ has bounded support, $\wmax = \infty$.
	 Then 
when choosing $\zeta(\cX)$ as $MM_w(\cX)$, the limit $\rho(\lambda)$ given by
(\ref{meanconv}) exists and is finite for all $\lambda \in (0,\infty)$, 
and if $n r_n^d \to t \in (0,\infty)$ as $n \to \infty$,
	then (\ref{0101c}) holds.
% Theorem~\ref{th:therm2}  applies.
\end{theo}

%In the case where $\wmax = +\infty$, we consider only the case
%where $\mu$ has bounded support,
% i.e. $\mu(Q_s)= 1 $ for some $s \in (0,\infty)$.
%We also require W4.

%\begin{theo}\label{lem:MMunbound}
%Suppose  $w(\cdot)$ satisfies conditions W1--W4, with
%	$\wmax = + \infty$.  Suppose $\mu$ has
%bounded support.
%Then when
% choosing $\zeta(\cX)$ as $MM_w(\cX)$, the limit $\rho(\lambda)$ given by
%(\ref{meanconv}) exists and is finite for all $\lambda \in (0,\infty)$, 
%and if $n r_n^d \to t \in (0,\infty)$ as $n \to \infty$,
%	then (\ref{0101c}) holds.
%\end{theo}
% Theorem~\ref{th:therm2}  applies.

We first prove Theorem \ref{lem:MM1} under the extra condition W5.

%The next lemma provides the first step towards the preceding theorems.
\begin{lemm}
\label{LemMMW5}
Suppose $w(\cdot)$ satisfies W1, W2 and W5. 
Then when choosing 
$\zeta(\cX)$ as $MM_w(\cX)$, the limit $\rho(\lambda)$ given by
(\ref{meanconv}) exists and is finite for all $\lambda \in (0,\infty)$, 
and if $n r_n^d \to t \in (0,\infty)$ as $n \to \infty$,
	then (\ref{0101c}) holds.
\end{lemm}

\begin{proof}
Assume first that W5 holds with $c_5 =1$ (as well as
	W1 and W2). We check that Theorem \ref{th:therm2}
	is applicable when taking $\zeta(\cdot) := MM_w(\cdot)$
	in this case.

	Properties P1 and P2 are clear. To check P3, let
 $\cY, \cZ $ be disjoint finite subsets of $ \R^d$. 
If $\tau$ is an optimal matching of $\cY$ (i.e., a matching
of $\cY$ with $w(\tau)=MM_w(\cY)$), and $\tau'$ is an optimal
matching of $\cZ$, then we can create a matching  $\tau'' $ 
of $\cY \cup \cZ$ by taking all edges of $\tau \cup \tau'$, together with
possibly one more edge.
Then 
$$
MM_w(\cY \cup \cZ) \leq
w(\tau'') \leq w(\tau) + w(\tau') + \wmax = 
 MM_w(\cY)+MM_w(\cZ)+ \wmax,
$$ 
so that P3 holds with $c_1= \wmax$ here.

Now we check P4. Suppose $\tau$ is an optimal matching of $\cY \cup \cZ$, and
enumerate
the  edges in $\tau$ that have one endpoint in $\cY$
and one endpoint in $\cZ$, as 
$e_1,\ldots e_k$. For $ 1 \leq i \leq k$ 
 denote the endpoints of $e_i$ by $y_i, z_i$ with $y_i \in \cY$
and $z_i \in \cZ$.

Now create a matching $\tau'$ on $\cY$ as follows.
Choose a matching
$\tau_1$ on $\{y_1,\ldots,y_k\}$, and let $\tau'$ consist of
all edges of $\tau$ that have both endpoints in $\cY$, together with
the edges of $\tau_1$, together with one further edge if $k$ is odd 
and $\cY$ has an unmatched vertex in $\tau$.

Similarly, create 
 a matching $\tau''$ on $\cZ$ as follows. Choose a matching
$\tau_2$ on $\{z_1,\ldots,z_k\}$, and let $\tau''$ consist of
all edges of $\tau$ that have both endpoints in $\cZ$, together with
the edges of $\tau_2$, together with one further edge if $k$ is odd 
and $\cZ$ has an unmatched vertex in $\tau$.

Since $\tau' \cup \tau''$ is either a matching of $\cY \cup \cZ$, or is
 contained in such a matching, $|\tau' \cup \tau''| \leq |\tau|$,
so that $|(\tau' \cup \tau'') \setminus \tau| \leq
|\tau \setminus (\tau'  \cup \tau'') | = k$.
Then
\bean
MM_w(\cY) + MM_w(\cZ) 
& \leq &  w(\tau') + w(\tau'')
\\
& = & w(\tau) + \left( \sum_{e \in (\tau' \cup \tau'') \setminus \tau} w(e) \right) 
- \sum_{i=1}^k w(e_i)
\\
& \leq
& w(\tau)
 + \sum_{i=1}^k (\wmax -w(e_i)).
\eean
As before, we say that an edge $uv$ of $\tau$ is {\em short}
if $\|u-v \| \leq 1$, {\em long} if $\|u-v\|>1$.
In the last sum, the $i$th term is zero whenever
$e_i$ is a 
 long edge.
The number of  short edges
$e_i$ is at most $|\partial_{\cZ} \cY|$,
and therefore
$MM_w(\cY) + MM_w(\cZ) \leq MM_w(\cY \cup \cZ) + \wmax  | \partial_{\cZ} \cY|$.
This gives us P4 with $c_2 =\wmax$.

Since one can obtain a matching from a tour by removing alternate
edges, $MM_w(\X) \leq TSP_w(\X)$ for all $\X$. Hence
Property P5 
for
 $MM_w$ holds as a consequence of
 the corresponding property for $TSP_w$, which we checked earlier.

 Thus Theorem \ref{th:therm2} is directly applicable in the
 case where $c_5=1$, giving us the desired conclusion in that case.
 We can then extend the result to the general case $c_5 \in (0,\infty)$
 in the same manner as we did in the proof of Lemma
 \ref{LemTSPW5}.
\end{proof} 

\begin{proof}[Proof of Theorem \ref{lem:MM1}(a)]
Assume that $w(\cdot)$ satisfies W1--W3 with $0 < \wmax < \infty$
(and hence  also W4). Then
	%Since we assume W3 with $\wmax <\infty$, 
	 for $a \in (0,\wmax)$ we have
	that $w_a(\cdot):= \min(w(\cdot),a)$ satisfies W5
	(as well as W1 and W2).
	So by Lemma \ref{LemMMW5},
	the limit
	$\rho_a (\lambda) := \lim_{s \to \infty} \E[ MM_{w_a}(\H_{\lambda,s}) 
	/(\lambda s^d) ]$ exists and is finite for all $\lambda >0$,
	and if $n r_n^d \to t \in (0,\infty) $, 
	%as $n \to \infty$,
	then $ n^{-1} MM_{w_a}(r_n^{-1} \X_n) \toccL \int \rho_a(t f_\mu(x))
	f_\mu(x) dx$.

	For any finite $\X \subset \R^d$,  $0 \leq MM_w(\X) - MM_{w_a}(\X)
	\leq |\X|(\wmax -a)$. Therefore, setting $a(k) := \max(\wmax - 1/k,0)$
	for each $k \in \N$,  we can apply Lemma 
	\ref{lemlong}(b), taking $\zeta^{(k)}(\cdot) := MM_{w_{a(k)}}(\cdot)$
	and $\zeta(\cdot) := MM_w(\cdot)$, to deduce the desired conclusion.
\end{proof}

	%Now we aim to prove Theorem \ref{lem:MMunbound},
	Now we aim to prove Theorem \ref{lem:MM1}(b),
	so we assume
	%$w(\cdot)$ satisfies
	W1--W4 with $\wmax = +\infty$.
	Initially we shall also assume W6. 
	Again set $w_a(\cdot):= \min(w(\cdot),a)$.

\begin{lemm}
	\label{lemMMdet}
	Suppose $w(\cdot)$ satisfies W1--W4 and W6 with $\wmax =+\infty$.
	Then there exists $C \in (0,\infty)$ such that for all
	$a \in (0,\infty)$, 
	$s \geq 1$ and all finite $\X \subset Q_s$ we have
\bea
	MM_{w_a}(\X) \leq
	MM_{w}(\X) \leq
	C s^p |\X|^{(d-p)/d}
\label{0316b}
\eea  
	and
	\bea
	0 \leq MM_w(\X) -  MM_{w_a}(\X) \leq C s^p a^{(p/d)-1}
	MM_{w_a}(\X)^{(d-p)/d}.
\label{0316a}
\eea
\end{lemm}
\begin{proof}
	Let $c_6$ be as in W6.
	By \cite[(3.6) or (3.7)]{Yukich},
	there is a constant $C>0$ such that 
	for any finite $\cY \subset Q_1$, 
$MM_{c_6\|\cdot\|^p}(\cY) \leq C |\cY |^{(d-p)/d}$.
 Hence by scaling, 
	%for finite $\cY \subset Q_1$,
$ MM_w(s \cY ) \leq MM_{c_6 \|\cdot\|^p} (s\cY) \leq C s^p |\cY|^{(d-p)/d} $.
	%for all $a \in (0,\infty)$, 
	Then % second inequality of
	(\ref{0316b}) follows.
	%(and the first inequality of (\ref{0316b}) is trivial).

Now let $\tau_a(\X)$ be (the set of edges in)  a
minimum-weight matching of $\X$, using weight function $w_a$ (and using some
%measurable
	deterministic rule to choose if there exist several matchings of
minimum weight). That is, let $\tau_a(\X)$ be a  matching of $\X$
with $w_a(\tau_a(\X)) = MM_{w_a}(\X)$.   
Let $\tau_{(=a)} (\X) = \{e \in \tau_a(\X): w_{a}(e) =a\}$, 
and  let $N_a (\X) := |\tau_{(=a)}(\X)|$.
 Thus $N_a(\X) $ is the  total number of edges
of weight $a$ in the minimum $w_a$-weight matching of
$\X$.
 Let $\cY_a(\X)$ be the set of vertices in $\X$ that
are incident to edges in $\tau_{(=a)}(\X)$.

Let $\tau'$ be a matching of $\cY_a(\X)$ with
$w(\tau') = MM_w(\cY_a(\X))$.
%n optimal matching, using weight function $w$,
%on $\cY_a(\X)$.
	Then $\tau' \cup (\tau_a(\X)  \setminus \tau_{(=a)}(\X))$ 
	is
	a matching on $\X$ with $w$-weight at most $w_a(\tau_a(\X)) + w(\tau')$.
	Therefore $MM_w(\X) \leq MM_{w_a}(\X) + MM_w(\Y_a(\X) )$, and hence
%using Therefore
	by  the second inequality of (\ref{0316b}) applied to $\cY_a(\X)$,
 $$
MM_w(\X) - MM_{w_a}(\X) \leq
%W^{(a)}(\X) \leq
 C s^p (2 N_a(\X))^{(d-p)/d},   ~~~~~ \X \subset Q_s, \X {\rm ~
finite.} 
$$
	Since $a N_a (\X)  \leq MM_{w_a} ( \X)$,
	this gives us (\ref{0316a}).
\end{proof}

%\begin{proof}[Proof of Theorem \ref{lem:MMunbound}]
	\begin{proof}[Proof of Theorem \ref{lem:MM1}(b)]
		We can follow the same argument as in
	%the proof of Theorem \ref{lem:TSPunbound}, now using
		the proof of Theorem \ref{lem:TSP}(b), now using
	Lemma \ref{LemMMW5} instead of Lemma \ref{LemTSPW5}
	and
	Lemma \ref{lemMMdet} instead of Lemma \ref{lemTSPdet}.
\end{proof}

	\subsection{Minimum-weight bipartite matching}		

Given two disjoint finite sets $\cU,\cV \subset \R^d$,
and given a near-perfect matching $\tau$ of $\cU \cup \cV$, define
$w^*(\tau) := \sum_{e \in \tau} w^*(e)$, with $w^*(xy) := w(xy)$
if $x \in \U, y \in \cV$ or $x \in \cV, y \in \U$, and
$w^*(xy) := \wmax$ if $x,y \in \U$ or $x,y \in \cV$. In other words,
if we say vertices in $\cU$ are of type 1 and
vertices in $\cV$ are of type 2,
we allow edges between vertices of the same type, but give them
a weight which is at least as big as that of any edge between
vertices of different types. We then define 
$BM_w(\cU,\cV) := \min_\tau w^*(\tau)$, where the minimum is
taken over all near-perfect matchings of $\cU \cup \cV$.
When $|\cU|=|\cV|$, clearly the minimum can be attained 
using a {\em bipartite
 matching}, i.e. a  perfect matching
in which each edge is between vertices of different types.
(Hence, when $|\cU|= |\cV|$, $BM_w(\U,\cV) = \min_\tau w(\tau)$
with the minimum taken  over all bipartite matchings.) 
Given $r>0$, we define 
$BM_{w,r}(\cU,\cV):= BM_{w}(r^{-1} \cU,r^{-1} \cV)$.

For $\lambda, s >0$,
let $\H'_\lambda$ denote a homogeneous Poisson process in $\R^d$ of intensity  
$\lambda$, independent of $\H_\lambda$, and set $\H'_{\lambda,s}:=
\H'_\lambda \cap Q_s$.
Let $U_0,V_0,U_1,V_1,U_2,V_2,\ldots$ be independent random
$d$-vectors with common distribution $\mu$.
For all $n \in \N$,
let $\cU_n:=  \{U_1,\ldots,U_n\}, \cV_n := \{V_1,\ldots,V_n\}$.
 We have the following law of large numbers
 for $BM_w(\cU_n,\cV_n)$ in the
	thermodynamic limit, when $\wmax < \infty$. 
	\begin{theo}[LLN for BM]
\label{MBMthm}
	Suppose $w(\cdot)$ satisfies W1--W3 with $\wmax < \infty$.
	Then for all $\lambda >0$,  the limit
 $\rho_{BM_w}(\lambda):= \lim_{s \to \infty} 
	\E [ BM_{w}(\H_{\lambda,s}, \H'_{\lambda,s})/(\lambda s^d)]$
	exists and is finite. Also, if 
	if $n r_n^d \to t \in (0,\infty)$ as $n \to \infty$, then
\bea
n^{-1} BM_{w,r_n}(\cU_n,\cV_n) \toccL \int_{\R^d} \rho_{BM_w}(t f_\mu(x)) f_\mu(x)
dx.  
\label{0227a}
\eea
\end{theo}

{\bf Remarks.}     
(a) The case with $w(x)= \|x\|^p$ for some fixed $p >0$
is considered in \cite{BB} and references therein; 
it is {\em not} covered by
our
Theorem \ref{MBMthm} 
because we require $\wmax < \infty$.
 For this case, results along the lines of 
(\ref{0227a}) have been shown
  \cite{BB} for $p <  d/2$, $d \leq 2$ or $p \leq 1, d \geq 3$
and $f_\mu $ being uniform over a bounded region, or for  $p< d/2$ with
$\mu$ being uniform over a bounded region with no singular part.
In contrast, our result (\ref{0227a}) for bounded $w$
holds for {\em any} choice of $\mu$.

	(b) One may also consider the {\em Bipartite TSP} with
	weight function $w$
	(the case with $w(\cdot) = \|\cdot\|^p$ is considered in 
	\cite{BB,Caracciolo}).
	To define this, let $BTSP_w(\cU,\cV)
	:= \min_{\tau} w^*(\tau)$,  now taking the minimum over all
	{\em tours} $\tau$ of $\cU \cup \cV$; let $BTSP_{w,r}(\cU,\cV)
	:= BTSP_{w}(r^{-1} \cU,r^{-1} \cV)$. If $|\cU| = |\cV|$, then
	$BTSP_w(\cU,\cV)$ is the minimum weight of all tours through
	$\cU \cup \cV$ that alternate between $\cU$ and $\cV$.

	The statement of
	 Theorem \ref{MBMthm} also holds with the functional 
	$BM_w(\cdot)$ replaced by $BTSP_w(\cdot)$ and $BM_{w,r}(\cdot)$
	replaced by $BTSP_{w,r}(\cdot)$. The proof is very similar.

\begin{lemm}
	\label{lemMBMW5}
	Suppose $w(\cdot)$ satisfies W1--W3 and W5. Then
	the conclusion of Theorem \ref{MBMthm} holds.
\end{lemm}
\begin{proof}
	First assume $c_5 =1$ in W5.
 Given finite $\X \subset \R^d$, suppose each point
is assigned a mark in $\{1,2\}$. 
Write $\X^*$ for the resulting marked point set
(a subset  of $\R^d \times \{1,2\}$). 
Define $\zeta^m(\X^*)$ to be 
$BM_w(\cU,\cV)$ where $\cU$ is the set of vertices of $\X$
with mark 1, and $\cV : = \X \setminus \U$. 
Then  define $\zeta(\X):= \E[\zeta^m(\X^*)]$, where
 the expectation is over a random marking of
$\X$ where each point is independently marked 1 or 2 with equal 
 probability. 
Likewise, for $r>0$  define $\zeta^m_r(\X^*) := BM_{w,r}(\cU,\cV)$ and
$\zeta_r(\X):= \E[\zeta^m_r(\X^*)]$.

Given two marked sets $\Y^*$ and $\cZ^*$, just as in the proof 
of Lemma \ref{LemMMW5}  we have
\bea
	\zeta^m(\cY^*) + \zeta^m(\cZ^*) - \wmax |\partial_\cZ\cY|  & \leq &
\zeta^m(\cY^* \cup \cZ^*) 
\nonumber	\\
	&	\leq & 
\zeta^m(\cY^*) + \zeta^m(\cZ^*) +\wmax,
\label{twomarked}
\eea
so for any disjoint finite $\cY, \cZ \subset \R^d$,
	taking expectations over the random markings of $\Y,\cZ$ gives us
$$
\zeta(\cY) + \zeta(\cZ) - \wmax |\partial_\cZ\cY| \leq
\zeta(\cY \cup \cZ) \leq 
\zeta(\cY) + \zeta(\cZ) +\wmax.
$$
Therefore $\zeta$ satisfies P1--P4. 

Next we verify P5 for $\zeta$. Let $\eps >0$. By W2, we can and do
choose $\delta \in (0,\eps^2)$ such
that $w(x) \leq \eps$ for all $x \in B_{2 \delta}(o)$. 
Given $\X \subset B_\delta(o)$ with $|\X| \geq \delta^{-1}$,
 set $n = |\X|$. 
Let $M_n$ be a Binomial($n,1/2$)
variable, representing the number of elements of $\X$ assigned
to type 1 rather than to type 2 in the randomly marked point
set $\X^*$. Then 
$$
\zeta^m(\X^*) \leq   \eps \min(M_n,n-M_n) + |M_n - (n/2)| \wmax, 
$$
so that taking expectations and using Jensen's inequality yields 
\bean
	\zeta(\X) \leq \eps (n/2) + \wmax \E[|M_n- (n/2)|]
	\leq \eps (n/2)  +  \sqrt{ \Var M_n} \wmax
	\\
	= \eps (n/2) + \sqrt{n/4} \wmax,   
\eean
and  therefore $n^{-1} \zeta (\X) 
	\leq (1/2) ( \eps + n^{-1/2} \wmax) \leq \eps(1+\wmax)/2 $,
which yields 
%(\ref{0229a}), i.e.
 P5.

Thus Theorem \ref{th:therm2} is applicable here, so for all $\lambda \in (0,\infty)$
the limit $\rho_1(\lambda):= \lim_{s \to \infty} \E[
	\zeta(\cH_{\lambda,s})/(\lambda s^d)]$
 exists in $\R$, 
and by  (\ref{0302a}) we have as $n \to \infty$ that
\bea
(2n)^{-1} \E[\zeta_{r_n}(\Po_{2n}) ] \to \int_{\R^d} \rho_1(2 t  f_\mu(x)) 
f_\mu(x)dx.
\label{0228a}
\eea

By the Marking theorem \cite{LP}, 
for $\lambda > 0$ the randomly marked point process $\H_{2\lambda}^*$ 
decomposes as the union of two independent copies of $\H_\lambda$.
Therefore
$$
\rho_1(2 \lambda)=  \lim_{s \to \infty}
 \E[ BM_{w} ( \H_{\lambda,s},\H'_{\lambda,s})/(2 \lambda s^d) ]
= (1/2) \rho_{BM_w}(\lambda), 
$$
where $\rho_{BM_w}(\lambda)$ was defined just
after (\ref{0227a}). 

For $n >0$, let $N_n,N'_n$ be independent Poisson$(n)$  random variables,
independent of $(U_1,V_1,U_2,V_2, \ldots)$.
Again by the Marking theorem, and by \cite[Prop. 3.5]{LP},
the left hand side of (\ref{0228a}) is equal to
$(2n)^{-1} \E[BM_{w,r_n}( \cU_{N_n},\cV_{N'_n})]$.
Thus (\ref{0228a}) becomes
\bea
n^{-1} \E[BM_{w,r_n}(\cU_{N_n},\cV_{N'_n}) ] \to \int_{\R^d} \rho_{BM_w}(t f_\mu(x)) 
f_\mu(x) dx.
\label{0228b}
\eea

By taking $|\cY^*|=1$ in  (\ref{twomarked}),
we see that adding one element to $\cU$ (or $ \cV$)
changes  $BM_{w,r_n}(\cU,\cV)$ by at most $\wmax$. Thus by iteration,
$$
|BM_{w,r_n}( \cU_{N_n},\cV_{N'_n}) - BM_{w,r_n}( \cU_{n},\cV_{n}) |
\leq \wmax(|N_n -n| + |N'_n-n|) , 
$$
and hence by (\ref{0228b}), since $\E[|N_n-n|] = o(n)$,
as $n \to \infty$ we have
\bea
n^{-1} \E[BM_{w,r_n}(\cU_{n}, \cV_n) ] \to \int_{\R^d}
 \rho_{BM_w}(t  f_\mu(x)) 
f_\mu(x)dx. 
\label{0228c}
\eea

We now show the complete and $L^2$ convergence. Let  $\eps >0$.
  For $1 \leq i \leq n$,
let 
 %$(\cU_{n,i}, \cV_{n,i}) := \{(U_1, V_1), \ldots,(U_0, V_0)\ldots,(U_n, V_n)\}$
 $\cU_{n,i}  := \{U_1, \ldots,U_0,\ldots,U_n\}$
where $U_0$ appears as the $i$th item in the list, and
let  $\cV_{n,i}  := \{V_1, \ldots,V_0,\ldots, V_n\}$
where $V_0$ appears as the $i$th item in the list.
% and $U_0, V_0$ are independent random $d$-vectors with
% common distribution $\mu$. 
That is,
$(\cU_{n,i}, \cV_{n,i})$ is the pair of point processes $(\cU_n, \cV_n)$ with
$(U_i, V_i)$ replaced by $(U_0, V_0)$. Let $\F_{i}$ be
the $\sigma$-algebra generated by $((U_1, V_1),\ldots, (U_i, V_i))$, and let
$\F_0$ be the trivial $\sigma$-algebra. Then
\bea
BM_{w,r_n}(\cU_n, \cV_n) - \E[ BM_{w,r_n}(\cU_n, \cV_n) ]
= \sum_{i=1}^n
D_{n,i},
\label{martdif2}
\eea
where for each $i \in [n]:= \{1,\ldots,n\}$ we set
\bean
D_{n,i} := 
\E[ BM_{w,r_n}(\cU_n, \cV_n) |\F_i  ]
- \E[ BM_{w,r_n}(\cU_n, \cV_n) |\F_{i-1}   ]
\\
= 
\E[ BM_{w,r_n}(\cU_{n}, \cV_n) - BM_{w,r_n}(\cU_{n,i}, \cV_{n,i}) |\F_i  ].
\eean
By the conditional Jensen inequality, and since replacing
$(U_i, V_i)$ by $(U_0, V_0)$ can change $\zeta_{r_n}$ by at most $4 \wmax$ 
by (\ref{twomarked}),
% (see also smoothness, Lemma~\ref{lemsmooth}),
\bea
|D_{n,i}| \leq 
\E[ (|BM_{w,r_n}(\cU_{n}, \cV_n) - BM_{w,r_n}(\cU_{n,i}, \cV_{n,i})|) |\F_i  ]
\leq 4 \wmax,
\label{200113b}
\eea
almost surely. Also
 $(D_{n,1},\ldots,D_{n,n})$ is a  sequence of martingale differences
with respect to the filtration $(\F_0,\F_1,\ldots,\F_n)$.
Therefore by Azuma's inequality (see e.g. \cite{Yukich} or \cite{RGG}),
\bean
\Pr[| BM_{w,r_n}(\cU_n, \cV_n) - \E[ BM_{w,r_n}(\cU_n, \cV_n) ] |> \eps n ]
\leq 2 \exp \left( \frac{ - (\eps n)^2}{ 32 \wmax^2 n  } \right),
\eean
which is summable in $n$. Combined with 
%~\eqref{0228b} and~
\eqref{0228c},
this demonstrates  the  complete convergence in~\eqref{0227a}.

For the $L^2$ convergence, note that
by (\ref{martdif2}), and  the orthogonality of martingale differences,
and (\ref{200113b}),
$$
\Var [n^{-1} BM_{w,r_n}(\cU_n, \cV_n)] = n^{-2} \sum_{i=1}^n
\E[D_{n,i}^2] \leq 16 \wmax^2 /n,
$$ 
which tends to 0
 as $n \to \infty$,
%and hence $\Var [n^{-1}\zeta_{r_n}(\X_n)] \to 0$,
% completing the
%proof of (\ref{0101c}). 
 so  the $L^2$ convergence in~\eqref{0227a} follows as well.
\end{proof}
\begin{proof}[Proof of Theorem \ref{MBMthm}]
	Assume $w(\cdot)$ satisfies W1--W3 with $\wmax < \infty$.
	Given $a \in (0,\wmax)$, $w_a(\cdot) := \min(w(\cdot),a)$ satisfies
	W5.  By Lemma \ref{lemMBMW5},
	for all $\lambda >0$
	the limit 
	\bea
	\rho_{BM_{w_a}} (\lambda)
	:= \lim_{s \to \infty} 
	\E[BM_{w_a}(\H_{\lambda,s},\H'_{\lambda,s})/(\lambda s^d)]
	\label{0322a}
	\eea
	exists,  and if $n r_n^d \to t \in (0,\infty)$ as $n \to \infty$,
	then 
	\bea
	n^{-1} BM_{w_a}(\U_n,\cV_n) \toccL
	\int \rho_{BM_{w_a}}(t f_\mu(x)) f_\mu(x) dx.
	\label{0322b}
	\eea
	Also $ 0 \leq BM_w(\U,\cV) - BM_{w_a}(\U,\cV)
	\leq (\wmax -a) |\U \cup \cV|$ for any finite $\cU,\cV \subset \R^d$.
	Therefore setting $a(k) := \max(\wmax -1/k,0)$ for each $k \in \N$,
	using the fact that (\ref{0322a}) and (\ref{0322b})  
	hold for $a=a(k) $,
	taking $k \to \infty$ 
	and
	%we can apply Lemma \ref{}(b), taking $\zeta^{(k)} :=
	%so by taking $a \uparrow \wmax$ and  
	arguing similarly to the proof of Lemma \ref{lemlong}(b),
	we can replace $w_a$ by $w$ in (\ref{0322a})
	and (\ref{0322b}), which is the desired conclusion.
\end{proof}

\subsection{Weighted minimum spanning tree (MST) }
\label{subsecMST}
Given  finite $\X \subset \R^d$, 
%Again consider the complete graph on vertex set $\cX$.
 the weighted minimum spanning tree (MST) problem
is to find a spanning tree on vertex set $\X$ of minimum
 weight.
 More formally, for each spanning tree $\tau$
in the complete graph on vertex set $\X$, define $w(\tau) = \sum_{e \in E(\tau)}
w(e)$, where $E(\tau)$ denotes the set of edges of $\tau$.
 Define 
$MST_w(\cX) := \min w(\tau) $, where the minimum is over
all spanning trees in the complete graph on vertex set $\X$.

\begin{theo}\label{lem:MST1}[LLN for weighted MST]
%Suppose  $w(\cdot)$ satisfies conditions W1--W3, with $\wmax < \infty$.  
	Suppose  $w(\cdot)$ satisfies conditions W1--W4, and either
	(a)  $\wmax < \infty$ or (b) $\mu$ has bounded support, 
	$\wmax = \infty$. 
	Then 
when choosing $\zeta(\cX)$ as $MST_w(\cX)$, the limit $\rho(\lambda)$ given by
(\ref{meanconv}) exists and is finite for all $\lambda \in (0,\infty)$, 
and if $n r_n^d \to t \in (0,\infty)$ as $n \to \infty$,
	then (\ref{0101c}) holds.
\end{theo}
In the case where $w (\cdot) = 1- {\bf 1}_{B_1(o)} (\cdot)$, $MST_w(\cX)$ is
the number of components of $G(\X,1)$, minus 1. 
Thus Theorem \ref{th:numcon} is a corollary of Theorem \ref{lem:MST1}(a).

%\begin{theo}\label{lem:MSTunbound}
%Suppose  $w(\cdot)$ satisfies conditions W1--W4, with
%	$\wmax = + \infty$.  Suppose $\mu$ has
%bounded support.
%Then when
% choosing $\zeta(\cX)$ as $MST_w(\cX)$, the limit $\rho(\lambda)$ given by
%(\ref{meanconv}) exists and is finite for all $\lambda \in (0,\infty)$, 
%and if $n r_n^d \to t \in (0,\infty)$ as $n \to \infty$,
%	then (\ref{0101c}) holds.
%\end{theo}
There is some overlap between Theorem \ref{lem:MST1}
%these results
and
 \cite[Theorem 2.3]{PY}. However, 
we here allow $\mu$ to have a singular part, and do not
require $w(x)$ to be a monotonic increasing function of
$\|x\|$, unlike in \cite{PY}.

We work 
towards proving the preceding theorems.
We shall first consider $w(\cdot)$ satisfying an extra condition,
namely that
 {\em there exists $\delta >0$
such that $w(x) =0$ whenever $\|x\| \leq \delta$}; we call this condition
W7. 

Given finite $\X \subset \R^d$ and given a spanning tree $\tau$
on vertex set $\X$,  we shall say, as before, that an edge $uv$ of $\tau$ is {\em short}
if $\|u-v \| \leq 1$, {\em long} if $\|u-v\|>1$.

\begin{lemm}
\label{lemMSTP4}
Suppose $w(\cdot)$ satisfies W1, W5 with $c_5=1$, and W7. 
There is a constant $k$ with the following property.
Given any finite $\cX \subset \R^d$,
 there exists a spanning tree $\tau$ on $\cX$
with $w(\tau) = MST_w(\cX)$, 
such that for each $x \in \cX$ there are at most $k$
short edges
in $\tau$ incident to $x$.
\end{lemm}
\begin{proof}
%Let $\cY, \cZ \subset \R^d$ be finite and disjoint.
	Let $\delta $ be as in W7. 
Partition $\R^d$ into half-open cubes of side $\delta/d$. 
Label  these
 cubes as $C_1,C_2,\ldots$.
% in the partition which intersect $\cX$ as $C_1,\ldots,C_m$.

Let $\cI:= \{ i: |\cX \cap C_i  | \geq 2\}$.
For each $i \in \cI$,
choose a path $\pi_i$ through $\cX \cap C_i$. The edges of this path
all have weight zero. Therefore we can then create a minimum-weight 
spanning  tree $\tau$ on $\X$ 
 that includes the edges of all of the paths $\pi_i, i \in \cI$.
 For any two distinct cubes $C_i,C_j$, $\tau$ has at most
one edge from $\X \cap C_i$ to $\X \cap C_j$ (since it is a tree).

Let $k'$ 
be the number of cubes in the partition
lying within unit Euclidean distance of $C_1$. This is also  
 the number of cubes in the partition
lying within unit distance of $C_i$, for any $i \in \N$.
For any
 $x \in \cX$ 
the  number of short edges in $\tau$  
from $x $  is bounded by $k:= k'+2$, as required.
\end{proof}

\begin{lemm}
\label{LemMSTW5}
Suppose $w(\cdot)$ satisfies W1,  W5 and W7. 
Then when choosing 
$\zeta(\cX)$ as $MST_w(\cX)$, the limit $\rho(\lambda)$ given by
(\ref{meanconv}) exists and is finite for all $\lambda \in (0,\infty)$, 
and if $n r_n^d \to t \in (0,\infty)$ as $n \to \infty$,
	then (\ref{0101c}) holds.
\end{lemm}

\begin{proof}
	We assume first that $c_5=1$ in condition W5.
	Clearly $ \zeta(\cdot) := MST_w(\cdot)$ satisfies
	Properties P1 and P2.
To check P3, 
 let $\cY, \cZ$ be  disjoint finite subsets of $ \R^d$. Let
 $\tau$ be a spanning tree on $\cY$ with weight 
$MST_w(\cY)$,
and let $\tau'$ be a spanning tree on $\cZ$ with weight $MST_w(\cZ)$.
If we combine the edges of $\tau$ and of $\tau'$, along with a single
edge from a vertex of  $\cY$ to a vertex of $\cZ$, we obtain a spanning
tree on $\cY \cup \cZ$ with total weight at most  
$MST_w(\cY)
+ MST_w(\cZ) +\wmax$. Hence,
$
MST_w(\cY \cup \cZ) \le MST_w(\cY) +MST_w(\cZ)+\wmax,
$
so P3 holds with $c_1=\wmax$.
 %since by combining the edges  spanning tree on $ the optimal solution of $\cY \cup \cZ$ is clearly bounded from above by an optimal solution on $\cY$, together with an edge from any point visited in $\cZ$ to $\cZ \setminus \cY$ of weight at most $1$, an optimal solution on $\cZ \setminus \cY$, which is clearly bounded by an optimal solution on $\cZ$, since the weight function is non-negative. Hence, Property P3
%\ref{propsub} holds with $c_1=1$.

%We now
To check Property P4,
 %Let $\cY,\cZ$ be as above.
let $\tau$ be a spanning tree on $\cY \cup \cZ$ with $w(\tau) =
MST_w(\cY \cup \cZ)$, having the property described in Lemma \ref{lemMSTP4}.
Remove from $\tau$ all edges going from $\cY$ to $\cZ$. This leaves us
with a forest, each component of which has either all vertices in $\cY$
or all vertices in $\cZ$. Then add edges to this forest, 
connecting up the components in $\cY$ to make a spanning tree 
$\tau'$ in $\cY$ and 
connecting up the components in $\cZ$ to make a spanning tree 
$\tau''$ in $\cZ$.
The total number of added edges equals the total number of previously
removed edges, minus 1.  

%We thus end up with a spanning tree $\tau'$ on $\cY$ and
% a spanning tree $\tau''$ on $\cZ$, and 
Then $w(\tau') + w(\tau'') - w(\tau)$
equals the total weight of added edges, minus the total
weight of removed edges. This is at most $\wmax$ times the number of
added edges, minus the total weight of removed edges. Therefore
it is at most  $\wmax$ times
	the total number of short removed edges, since
the long removed edges have weight $\wmax$, and there are more
removed edges than added edges. Since $\tau$ has the
property  described in Lemma \ref{lemMSTP4},  for each $y \in \cY$
the number of short removed edges incident to $y$ 
is at most $k$, and is zero if $y \notin \partial_\cZ \cY$, so that
\bean
MST_w(\cY) + MST_w(\cZ) - MST_w(\cY \cup \cZ) & \leq
& w(\tau') + w(\tau'') - w(\tau) \\
& \leq & k \wmax | \partial_{\cZ} \cY|,
\eean
which gives us Property P4 with $c_2 = k \wmax$.
Property P5 follows because 
 $MST_w(\cdot) \leq TSP_w(\cdot)$ (since removing one edge from a tour yields
a spanning tree).
\end{proof}

As usual, for $a \in (0,\wmax)$ we set $w_a(\cdot):= \min(w(\cdot),a)$.

\begin{proof}[Proof of Theorem \ref{lem:MST1}(a)]
	Assume $w(\cdot)$ satisfies W1--W3 with $\wmax < \infty$
(and hence also W4).
	For $k \in \N$, set $a(k) := \max(\wmax - 1/k,0)$ and
	set 
	$$
	w^{(k)}(x):= w_{a(k)}(x)(1 - {\bf 1}_{B_{1/k}(o)}(x)), ~~~~
	x \in \R^d.
	$$
	Then $w^{(k)}(\cdot)$ satisfies W1, W5  and W7, so
	by Lemma \ref{LemMSTW5},
	%by the case already considered,
	for all $\lambda >0$ 
	the limit $\rho_k(\lambda) :=
	\lim_{s \to \infty} \E[MST_{w^{(k)}}(\H_{\lambda,s})/(\lambda s^d)]$
	exists and if $n r_n^d \to t \in (0,\infty) $ as $n \to \infty$,
	then $n^{-1} MST_{w^{(k)}}(r_n^{-1}\X_n) \toccL \int \rho_k(t f_\mu(x))
	f_\mu(x) dx$.

	Now, for all finite $\X \subset \R^d$, set $\zeta^{(k)}(\X):= MST_{w^{(k)}}(\X)$ and $\zeta(\X):= MST_w(\X)$. Then $\zeta^{(k)}(\X) $ is nondecreasing
	in $k $ with $\lim_{k \to \infty} \zeta^{(k)}(\X) = \zeta(\X)$. 
	Moreover, for all $x \in \R^d$, we have
	$$
	0 \leq w(x) - w^{(k)}(x) \leq (\wmax - a(k)) + \sup_{x \in B_{1/k}(o)}
	w(x) =: h(k),
	$$
	and $h(k) \to 0$ as $k \to \infty$ by assumption W2.
	Then $0 \leq \zeta(\X) - \zeta^{(k)}(\X) \leq h(k)|\X|$, and
	so we can apply
	 Lemma \ref{lemlong}(b) with this choice of $\zeta^{(k)}$,
	 to deduce
	 %By that result we obtain
	 that 
	for all $\lambda >0$ 
	the limit $\rho(\lambda) :=
	\lim_{s \to \infty} \E[MST_{w}(\H_{\lambda,s})/(\lambda s^d)]$
	exists and if $n r_n^d \to t \in (0,\infty) $ as $n \to \infty$,
	then $n^{-1} MST_{w}(r_n^{-1}\X_n) \toccL \int \rho(t f_\mu(x))
	f_\mu(x) dx$, as required.
\end{proof}

	Now we %consider unbounded $w$. 
	%We aim to prove Theorem \ref{lem:MSTunbound},
	aim to prove Theorem \ref{lem:MST1}(b),
	so we assume $w(\cdot)$ satisfies
	W1--W4 with $\wmax= +\infty$. Initially we shall also assume W6.
	Again set $w_a(\cdot):= \min(w(\cdot),a)$.

\begin{lemm}
	\label{lemMSTdet}
	Suppose $w(\cdot)$ satisfies W1--W4 and W6.
	There is a constant $C>0$ such that
	for all
	$a \in (0,\infty)$, all
	$s \geq 1$ and all finite $\X \subset Q_s$ we have
\bea
	MST_{w}(\X) \leq C s^p |\X|^{(d-p)/d}
\label{0316b2}
\eea  
	and
	\bea
	0 \leq MST_w(\X) -  MST_{w_a}(\X) \leq C s^p a^{(p/d)-1}
	MST_{w_a}(\X)^{(d-p)/d}.
\label{0316a2}
\eea
\end{lemm}
\begin{proof}
	Let $c_6 $ be as in W6.
	By \cite[(3.6) or (3.7)]{Yukich},
	there exists a finite constant $C >0$ such that 
 for any finite $\cY \subset Q_1$, 
we have
$MST_{c_6\|\cdot\|^p}(\cY) \leq C |\cY |^{(d-p)/d}$.
 Hence by scaling, 
	%for finite $\X \subset Q_s$,
$ MST_w(s \cY ) \leq MST_{c_6\|\cdot\|^p} (s\cY) \leq C s^p 
	|\cY|^{(d-p)/d} $,
which implies (\ref{0316b2}).

Now let $s \geq 1$ and let $\X \subset Q_s$ be finite.
	Let $\tau_a(\X)$ be (the set of edges in)  a
minimum-weight spanning tree on
	$\X$, using weight function $w_a$ (and using some
	deterministic rule to choose if there exist several matchings of
minimum weight). That is, let $\tau_a(\X)$ be a spanning tree of $\X$
with $w_a(\tau_a(\X)) = MST_{w_a}(\X)$. 
Let $\tau_{(=a)} (\X) = \{e \in \tau_a(\X): w_{a}(e) =a\}$, 
and  let $N_a (\X) := |\tau_{(=a)}(\X)|$.
 Thus $N_a(\X) $ is the  total number of edges
of weight $a$ in the minimum $w_a$-weight spanning tree of
$\X$.
 Let $\cY_a(\X)$ be a set consisting of one  vertex from
	each component of the graph $(\X,\tau_a(\X) \setminus \tau_{(=a)}(\X))$,
	chosen in an arbitrary measurable way. Note that 
	$|\cY_a(\X)| = |N_a(\X)|+1$.

	Let $\tau'$ be (the edge-set of) a spanning tree on $\cY_a(\X)$ with
$w(\tau') = MST_w(\cY_a(\X))$.
	Then $\tau' \cup ( \tau_a(\X)  \setminus \tau_{(=a)}(\X))$ 
	is (the edge-set of) a spanning tree on $\X$ with $w$-weight at most $w_a(\tau_a(\X)) + w(\tau')$. 
Therefore $MST_w(\X) \leq MST_{w_a}(\X) + MST_w(\Y_a(\X) )$, and hence
	by (\ref{0316b2}) applied to $\Y_a(\X)$, if $N_a(\X) \geq 1 $ then
 $$
MST_w(\X) - MST_{w_a}(\X) \leq
 C s^p ( N_a(\X)+1)^{(d-p)/d} \le 
	C s^p (2 N_a(\X))^{(d-p)/d}.
$$
	Since $a N_a (\X)  \leq MST_{w_a} ( \X)$,
	this gives us (\ref{0316a2}) (upon changing the
	constant $C$). If $N_a(\X)=0$ then $MST_{w}(\X)=
	MST_{w_a}(\X)  $  so (\ref{0316a2}) holds then too.
\end{proof}

%\begin{proof}[Proof of Theorem \ref{lem:MSTunbound}]
	\begin{proof}[Proof of Theorem \ref{lem:MST1}(b)]
		Assume that $w(\cdot)$ satisfies W1--W4 with $\wmax= +\infty$,
		and that $\mu$ has bounded support.
Given $\delta > 0$, $m \in \N$,
	define the new weight functions $w'_\delta(x) := w(x)
	(1- {\bf 1}_{B_\delta(o)}(x))$ and
	$w'_{\delta,m}(x) := \min(w_\delta(x),m)$.

	Then $w'_{\delta,m}$ satisfies W5 and W7. Hence by 
	Lemma \ref{LemMSTW5}, for all $\lambda >0$ the limit
	$\rho'_{\delta, m} (\lambda) := \lim_{s \to \infty}
	\E[MST_{w'_{\delta,m}} (\H_{\lambda,s})/(\lambda s^d)] $ 
	exists, and if $n r_n^d \to t \in (0,\infty)$ as
	$n \to \infty$, then $n^{-1} MST_{w'_{\delta,m}}(r_n^{-1} \X_n)
	\toccL \int \rho'_{\delta, m}( t f_\mu(x)) f_\mu(x)dx$.

	Since $w(\cdot)$ satisfies W4, the function 
$w'_\delta(\cdot)$ satisfies W6. Therefore by Lemma
	\ref{lemMSTdet}, 
	there is a constant $ C(\delta)$ such that for all $s >0$
	and all finite $\X \in Q_s$, we have $MST_{w'_{\delta,m}}(\X)
	 \leq MST_{w'_{\delta}}(\X)
	\leq C(\delta) s^p |\X|^{(d-p)/d}$ and 
	\bean
	MST_{w'_\delta}(\X) - MST_{w'_{\delta,m}}(\X) \leq
	C(\delta) s^p m^{(p/d)-1} (MST_{w'_{\delta,m}}(\X))^{(d-p)/d}.
	\eean
	Also, setting $\tilh(\ell):= \sup_{x \in B_{1/\ell}(o) } w(x)$, 
	we have
	$
	MST_{w}(\X) - MST_{w'_{1/\ell}}(\X) \leq \tilh(\ell) 
	|\X|.
	$
	Thus the functionals
	$\tzeta^{(\ell,m)}(\cdot):= MST_{w'_{1/\ell,m}}(\cdot)$ 
	satisfy the conditions for Lemma \ref{lemlong2} with
	$h_\ell(m): = C(1/\ell) m^{(p/d)-1} $.  Hence by
	Lemma \ref{lemlong2},
for all $\lambda >0$
the limit 
$
\rho(\lambda) := \lim_{s \to \infty}
\E[MST_{w}(\H_{\lambda,s})/ (\lambda s^d)]
$
exists and is finite, and also if $nr_n^d \to t \in (0,\infty)$ 
as $n \to \infty$, 
then
$
MST_{w}( r_n^{-1} \X_n ) \toccL \int 
	\rho(t f_\mu(x)) f_\mu(x)dx,
	$
	as required.
	%We may now conclude the proof in the same manner as we did for
	%Theorem \ref{lem:TSPunbound}.
\end{proof}

\section{Proof of the general results}
\label{secproofgen}
\allco

In this section we prove the results
stated in Section \ref{secstate}.
%Theorem \ref{th:therm2}. 
Assume throughout this section that the function $\zeta(\cdot)$ has
been specified, taking values in $[0,\infty)$ 
with $\zeta(\emptyset)=0$, satisfying Properties P1--P4.

The following elementary lemma concerns deterministic subadditive functionals of
Borel sets.  As before, let $\cB$ denote the class of bounded
Borel sets in $\R^d$, and let $\cB_1$ denote the class of all Borel sets 
 that are contained in $Q_1$.
Let us say that a set function 
$f: \cB \to \R$ is {\em Borel subadditive} 
if $f(A \cup A' ) \leq f(A) + f(A')$
 whenever $A,A' \in \cB$ with $A \cap A' = \emptyset$.
\begin{lemm}
\label{lemsubadd}
Suppose $f: \cB \to [0,\infty)$ is a Borel subadditive set function 
such that $f(x+A) = f(A)$ for all $A \in \cB$,
$x \in \R^d$, and also such 
that $ \sup\{f(A): A \in \cB_1\} < \infty$.
Then $\overline{f}: = \inf_{s \geq 1} (s^{-d} f(Q_s)) \in [0,\infty) $
and $s^{-d} f(Q_s) \to \overline{f} $ as $s \to \infty$. 
\end{lemm}

More sophisticated stochastic subadditive limit theorems
are available
(see \cite{AK}, or \cite[Theorem 4.9]{Yukich})  but are not
needed here. 
Since we use Lemma
\ref{lemsubadd} repeatedly,
for completeness we include a proof.

\begin{proof}[Proof of Lemma \ref{lemsubadd}]
Set $f^* := \sup\{f(A): A \in \cB_1\} $.
Clearly $0 \leq \overline{f} \leq f^* < \infty$.
 Given $\eps>0$,
choose $s_0 \geq 1$ such that $f(Q_{s_0}) < (\overline{f} + \eps) s_0^d$.   

Given $s >0$, write $s = ns_0+t$ with $t \in [0,s_0)$. Then
there is a constant $c$ such that for all such $s$
we can write $Q_s$ as a disjoint union of $n^d $ translates of
$Q_{s_0}$, together with at most $c n^{d-1}$ further sets in $\cB$ that
are translates of sets in $\cB_1$ (in fact, rectangles with all sides at 
most 1). 
By repeated use of Borel subadditivity,
$$
f(Q_s) \leq n^d f(Q_{s_0}) + c n^{d-1} f^* \leq (ns_0)^d (\overline{f}+ \eps)
+ c n^{d-1} f^*.
$$  
Therefore since $s \geq ns_0$, for large enough $s$ (and hence, large 
enough $n$) we have 
$$
s^{-d} f(Q_s) \leq \overline{f} + \eps + cf^* s_0^{-d} n^{-1}
< \overline{f}+ 2 \eps,
$$
which implies the result.
\end{proof}

Next we show that $\zeta_r$ enjoys a {\em smoothness} property;
the effect on $\zeta_r$ of adding  a single point to an existing
point set is uniformly bounded.
\begin{lemm}
\label{lemsmooth}
%\label{prop3} 
%{\em Smoothness:}
Let $K:= \max(c_1 + \zeta(\{o\}), c_2 - \zeta(\{o\}))$,
where $c_1,c_2$ are as in P3, P4 respectively. Then
%There is a constant $K \in \R$ such that
 $|\zeta_r(\cY) - \zeta_r(\cX)| \le K
% |\cY \Delta \cX|$,  
|\cY \triangle \cX|$
for all  finite $\cX, \cY \subset \R^d$, and all $r > 0$. 
\end{lemm}
\begin{proof}
It suffices to prove the result when $ |\cY \triangle \cX| =1$ and $r=1$.
Let $\cX \subset \R^d$ and $x \in \R^d \setminus \cX$.
By Property P2 (translation invariance), $\zeta(\{x\}) = \zeta(\{o\})$,
and by P3 (almost subadditivity) 
%\ref{propsub})
$$
\zeta(\cX \cup \{x\}) %= \zeta( r^{-1} \cX \cup \{r^{-1} x\} )
 \leq \zeta(\cX) + \zeta(\{o\})
+ c_1.
$$
Similarly,
since $|\partial_{\cX}(\{x\}) | \leq 1$,  by P4  
(superadditivity up to boundary)
%\ref{propsup})
 $$
\zeta(\cX \cup \{x\})  \geq \zeta(\cX) + \zeta(\{o\})
- c_2.
$$
  Combining these inequalities shows that
$|\zeta(\cX \cup \{x\} ) - \zeta(\cX) | \leq K$,
if we take $K = \max(c_1+  \zeta(\{o\}),c_2 - \zeta(\{o\}) )$. 
This gives the result.
\end{proof}

\begin{proof}[Proof of Proposition \ref{propsofrho}]
(i)
By Lemma \ref{lemsmooth},  for all  finite $\cX \subset \R^d$ we have that
$ \zeta(\cX) \leq K|\cX|$. 
Hence for each $A \in \cB$,
 the class of bounded Borel sets in $\R^d$,
 we have
\bea
\E[ \zeta(\cH_\lambda \cap A)]  \leq
 K \E[   |\cH_\lambda \cap A| ] 
=
K \lambda \Leb(A).  
\label{meanbound}
\eea
In particular, this expectation is finite.
 Given $\lambda >0$ and   $A \in \cB$, let 
$$
f_\lambda(A) := \E[\zeta( \cH_\lambda \cap A)]+c_1,
$$
where  $c_1 \in \mathbb{R}$ is the constant given in Property P3.

By Property P3,
for any disjoint  $A,A' \in \cB$
we have $f_\lambda(A \cup A') \leq f_\lambda( A ) + f_\lambda(A')$.
Also
for all $A \in {\cal B}$ and all $x \in  \R^d$
we have $f_\lambda(x+A) = f_\lambda(A)$.
Moreover, by (\ref{meanbound})  $f_\lambda$
has the other properties described in the hypothesis of Lemma
\ref{lemsubadd}. 
Hence by Lemma \ref{lemsubadd}, the limit in (\ref{meanconv}) does exist
in $\R$. We defer the proof of (\ref{meanconvstrong}) until later.

(ii) Suppose  $0 <  \lambda < \lambda' $. 
By the Superposition theorem (see e.g. \cite{LP})
 we can couple $\cH_\lambda$
and $\cH_{\lambda'}$ in such a way that $\cH_\lambda \subset \cH_{\lambda'}$ 
and $|(\H_{\lambda'} \setminus \H_{\lambda}) \cap Q_s|$
  is Poisson with parameter $(\lambda' - \lambda) s^d$, for all 
$s \in (0,\infty)$.
With this coupling, by smoothness (Lemma \ref{lemsmooth}), with
$K$ as in that result
 we have for any $s>0$ that
\bea
\E[|\zeta(\cH_{\lambda', s}) - \zeta(\cH_{\lambda, s})|]
 \leq K (\lambda'-\lambda) s^d,
\label{1230b}
\eea
and therefore  by (\ref{meanconv}),
we obtain that
$|\lambda' \rho(\lambda') -  \lambda \rho(\lambda) | 
 \leq K (\lambda' - \lambda ) $, so $\lambda \rho(\lambda)$ is Lipschitz
 continuous in $\lambda$ with Lipschitz constant at most $K$.
Hence $\rho(\lambda ) $ is also continuous in $\lambda$.

We now check (\ref{meanconvstrong}).
If $s_n \to \infty$ and $\lambda_n \to \lambda$, then given
$\eps >0$, for large enough $n$ we have 
$|s_n^{-d} \E[ \zeta(\cH_{\lambda,s_n})] - \lambda \rho(\lambda)|< \eps$
by (\ref{meanconv}), and
$s_n^{-d} \E |\zeta(\cH_{\lambda_n,s_n}) - \zeta(\cH_{\lambda,s_n})|
< \eps$ by 
(\ref{1230b}).
This gives us
 (\ref{meanconvstrong}).

(iii) By the proof of  Lemma \ref{lemsmooth},
for all finite  $\cX$ we have $\zeta(\cX) \leq
(c_1 + \zeta(\{o\})) |\cX|$.  Hence
for all $s>0$ we have 
$$
\E[ \zeta(\H_{\lambda,s}) ]  \leq  
(c_1+ \zeta(\{o\})) \E[|\H_{\lambda,s} | ] = 
(c_1+ \zeta(\{o\})) \lambda s^d.
$$
Hence by (\ref{meanconv}), 
$\rho(\lambda) \leq c_1 + \zeta(\{o\})  $. 

(iv) Suppose $0 < \lambda < \lambda'$, and couple $\cH_\lambda$
 and $\cH_{\lambda'}$ 
as in the proof of (ii).
Under assumption P6, for all $s >0$ we then have
$\zeta(\cH_{\lambda,s})
 \leq \zeta(\cH_{\lambda',s})$, so
 using
(\ref{meanconv}) we obtain that $\lambda \rho(\lambda ) \leq 
\lambda' \rho(\lambda')$, as required.

(v) 
Given $\lambda>0$, by (\ref{meanconv}) we have 
\bean
\rho(\lambda)
%= \lim_{s \to \infty}
%(\lambda^{1/d} s)^{-d} \E[ \zeta(\cH_\lambda \cap Q_s) ] 
= \lim_{s \to \infty}
(\lambda^{1/d} s)^{-d} \E[ \zeta_{\lambda^{1/d}}(\lambda^{1/d} 
\cH_\lambda \cap Q_{\lambda^{1/d} s } )  ],
\eean
and hence by the Mapping theorem for Poisson processes (see e.g. \cite{LP}),
\bean
\rho(\lambda) = \lim_{s \to \infty}
(\lambda^{1/d} s)^{-d} \E[ \zeta_{\lambda^{1/d}}( 
\cH_1 \cap Q_{\lambda^{1/d} s } )  ]
= \lim_{t \to \infty}
t^{-d} \E[ \zeta_{\lambda^{1/d}}( \cH_1 \cap Q_{t } )  ].
\label{0108b2}
\eean 
Hence by Property P7 (downward monotonicity in $r$),
$\rho(\lambda)$ is nonincreasing in $\lambda$.

(vi)
% Suppose that we can take $c_1=0$ in P3.
Let $\cH_\lambda^0$ be the set of points of $\cH_\lambda$
that are isolated vertices in the graph $G(\cH_\lambda,1)$. Using smoothness
(Lemma \ref{lemsmooth}), and then P4 repeatedly, we have 
\bean
\zeta(\cH_\lambda \cap Q_s) \geq 
\zeta( \cH_\lambda^0 \cap Q_s)  
- K | (\cH_\lambda \setminus \cH_\lambda^0) \cap Q_s|
\\
\geq \zeta(\{o\}) | \cH_\lambda^0 \cap Q_s | 
- K | (\cH_\lambda \setminus \cH_\lambda^0) \cap Q_s|.
\eean
Hence by the Palm-Mecke formula from the theory of Poisson processes
 (see \cite{RGG} or \cite{LP}), 
$$
\E[\zeta(\cH_{\lambda,s}) ] \geq \lambda s^d \zeta(\{o\}) 
e^{-\lambda \pi_d} - K \lambda s^d ( 1- e^{-\lambda \pi_d}) .
% \zeta(\cH_\lambda^0 \cap Q_r) 
%- K | (\cH_\lambda \setminus \cH_lambda^0) \cap Q_r|
$$
Hence by (\ref{meanconv}),  $\rho(\lambda) \geq \zeta(\{o\}) e^{-\lambda \pi_d}
- K (1 - e^{-\lambda \pi_d})$, and hence
$\liminf_{\lambda \downarrow 0} \rho(\lambda) \geq \zeta(\{o\})$.
Combined with part (iii) this gives us the result.

(vii)
Fix $0 < \lambda < \lambda' < \zeta(\{o\}) /(c_2 \pi_d)$.
 Our goal is to show that
 $\lambda'\rho({\lambda')} > \lambda \rho(\lambda)$. 
   Fix $s \in (0, \infty)$.
By the Superposition theorem \cite{LP}, 
 the Poisson point process  $\cH_{\lambda'} $
 may be obtained as the union of the Poisson point process
$\cH_{\lambda}$ 
and another fresh Poisson point process 
 of intensity $\lambda' - \lambda$ added on top of it.
Let $N$ denote the number of points in $Q_s$
of the second Poisson point process to be added, and let $N_0$
be the number of these points that are isolated in 
$G(\cH_{\lambda'} \cap Q_s,1)$,
i.e. have no other point of $\cH_{\lambda',s}$ within unit
Euclidean distance.

For each added point, if that point is isolated then by
P4 it increases the value of $\zeta$ by at least $\zeta(\{o\})$,
and if it is not isolated then by P4 it decreases the value
of $\zeta$ by at most $c_2 - \zeta(\{o\})$. Therefore
adding the points one by one, we obtain
 that
\bea
\zeta( \cH_{\lambda',s}) - \zeta(\cH_{\lambda,s}) 
\geq \zeta(\{o\}) N - c_2(N-N_0).
\label{1206a}
\eea
Now, $\E[N]= (\lambda'-\lambda)s^d$, and by the Palm-Mecke formula
 (see \cite{RGG} or \cite{LP}),
$$
\E[N-N_0] \leq 
 (\lambda'-\lambda)s^d (1 - \exp(- \lambda' \pi_d  )) \leq
 (\lambda'-\lambda)s^d  \lambda' \pi_d . 
$$
Therefore taking expectations in (\ref{1206a}) yields
$$
s^{-d}( \E[\zeta( \cH_{\lambda',s})]  - \E[ \zeta(\cH_{\lambda,s} )
] ) 
\geq
(\lambda' - \lambda)
 (\zeta(\{o\})  - c_2 \pi_d \lambda' ).
$$
By the choice of $\lambda$ and $\lambda'$,
the right hand side above is strictly positive.
Taking the large-$s$ limit and using \eqref{meanconv} shows that 
$\lambda' \rho(\lambda') > \lambda \rho(\lambda)$, as required.
\end{proof}

\begin{proof}[Proof of Proposition \ref{thsubcrit}]
Given a locally finite set $\cX \subset \R^d$, define the
{\em clusters} of $\cX$ to be the vertex sets of the
connected components of $G(\cX,1)$.
For $x \in \cX$ we define $\cC(x,\cX)$ to be
 the cluster of $\cX$ containing $x$. 

By the assumptions P3 (with $c_1=0$) and P4, if 
 $\cX$ and $\cY$ are disjoint finite subsets  of $\R^d$ such that
$\partial_\cZ \cY = \emptyset$, then $\zeta(\cY \cup \cZ) =
\zeta(\cY) + \zeta(\cZ)$. 
Thus for any finite $\cX \subset \R^d$, denoting
the clusters of $\cX$ by $\cC_1,\ldots,\cC_m$ we
have $\zeta(\cX) = \sum_{i=1}^m \zeta(\cC_i)$.
Thus for $s >0$,
$$
	\zeta(\cH_{\lambda,s}) = \sum_{x \in \cH_{\lambda,s}}
\frac{
	\zeta( \cC(x,\cH_{\lambda , s}))
	}{| \cC(x,\cH_{\lambda ,s})|}.
$$
Hence by the Palm-Mecke formula
% (see \cite{RGG} or \cite{LP}),
$$
	\E [ \zeta(\cH_{\lambda,s}) ] = \lambda \int_{Q_s}
\E \left[ \frac{  \zeta( \cC(x, \cH^x_\lambda \cap Q_s)) }{
|   \cC(x, \cH^x_\lambda \cap Q_s) | } \right] dx,
$$
where we write $\cH_\lambda^x$ for $\cH_\lambda \cup \{x\}$.
Taking $y=s^{-1}x$ we have
\bea
\lambda^{-1}
	s^{-d} \E \zeta(\cH_{\lambda,s}) = \int_{Q_1}
\E \left[  \frac{\zeta( \cC({sy}, \cH^{sy}_\lambda \cap Q_s)) }{
|   \cC({sy}, \cH^{sy}_\lambda \cap Q_s) | } \right] dy.
\label{1204a}
\eea

Let $Q_1^o$ denote the interior of $Q_1$, and let $y \in Q_1^o$.
 By P2 (translation invariance) and the stationarity
of $\cH_\lambda$ we have 
\bea
\E \left[  \frac{\zeta( \cC({sy}, \cH^{sy}_\lambda \cap Q_s)) }{
|   \cC({sy}, \cH^{sy}_\lambda \cap Q_s) | } \right] 
= 
\E \left[  \frac{\zeta( \cC(o, \cH^{o}_\lambda \cap (-sy+ Q_s))) }{
|   \cC(o, \cH^{o}_\lambda \cap ( -sy + Q_s)) | } \right]. 
\label{1204b}
\eea
Since  we assume $\lambda < \lambda_c$ the set $\cC(o,\cH_\lambda^o)$ is
almost surely finite, and hence for large enough $s$ this set is
contained in the set $-sy +Q_s$ (which is equal to $s(-y+Q_1)$).
Also, by smoothness (Lemma \ref{lemsmooth}) there
is a constant $K$ such that $|\zeta(\cX)| \leq K |\cX|$ for all finite
$\cX \subset \R^d$. Hence by the Dominated Convergence theorem,
as $s \to \infty$ 
the expression in (\ref{1204b}) converges to
$$
\E \left[  \frac{\zeta( \cC(o, \cH^{o}_\lambda))  }{
|   \cC(o, \cH^{o}_\lambda )  | } \right]. 
$$
By Dominated Convergence again,
the expression (\ref{1204a})  converges to the same limit, and then
the result follows from (\ref{meanconv}).
\end{proof}

It remains, in this section, to prove Theorem \ref{th:therm2}.
Hence, now and for the rest of this section we assume P5 as
well as P1--P4. The existence in $\R$ of the limit $\rho(\lambda)$ 
given by \eqref{meanconv}, for all $\lambda >0$,   
was already proved for Proposition~\ref{propsofrho}(i).

Given  $A \subset \R^d$ and $r >0$ 
let $\partial_{r}A$ denote the set of all points
in $A$ at a Euclidean distance at most  $r$ from $\R^d \setminus A$.
 The following result is 
a consequence of almost subadditivity (Property P3) and
 superadditivity up to boundary
 (Property P4).
\begin{lemm}
\label{lemfromsuper}
Let  $k \in \N$ and let $A_1, \ldots, A_k$ be disjoint subsets of $\R^d$.
 For any $r>0$ and any finite $\X \subset \cup_{i=1}^k A_i$ we have
\bea
\zeta_{r}(\X  ) \le \left( \sum_{i=1}^k 
 \zeta_{r}(\X   \cap A_{i}) \right) 
+ c_1 (k-1),  
\label{fromsub}
\eea
and
\bea
\zeta_{r}(\X  ) \ge \left( \sum_{i=1}^k 
%\left(
 \zeta_{r}(\X   \cap A_{i}) \right) 
- c_2 \sum_{i=2}^k  |\X \cap \partial_{r}  A_i|,  
%- ck |\X \cap \partial_{r}  A_i)|\right),  
\label{fromsuper}
\eea
where $c_1,$ $c_2$ are the constants in Properties P3, P4 respectively,
%\ref{propsup},
 and the second
sum in (\ref{fromsuper}) is interpreted as zero if $k =1$.
\end{lemm}
\begin{proof}
Since we assume Property P3
%\ref{propsub}
 (almost subadditivity of
$\zeta$), an analogous property also holds for $\zeta_r$. 
We then obtain (\ref{fromsub}) by a straightforward induction in $k$.

To prove (\ref{fromsuper}),
first suppose $k=2$. For $i=1,2$ set $\cY_i:= r^{-1}(\cX \cap A_i)$.
Then by P4,
\bean
\zeta_r(\cX) = \zeta(\cY_1 \cup \cY_2)
\geq \zeta(\cY_1 ) + \zeta (\cY_2) - c_2 |\partial_{\cY_1}\cY_2|.
\eean
Also, since $A_1$ and $A_2$ are disjoint,
$$
| \partial_{\cY_1}\cY_2 | \leq 
| r^{-1} \X \cap \partial_1(r^{-1} A_2)| 
= 
|  \X \cap \partial_r A_2|. 
$$
Together these two inequalities give us the result for $k=2$, and
the  general result follows by a straightforward induction on $k$.
\end{proof}

Following \cite{Yukich}, we say the distribution $\mu$ is {\em blocked}
if for some $\bigblock \in \N$ and some
 $m \in \N$ a power of $2$, the density function $f_\mu$
of the absolutely continuous part of $\mu$
% (in its Lebesgue decomposition)
 is constant on
each of the $(\bigblock m)^d$  cubes in
 the subdivision of $Q_\bigblock$ into 
%$m^d$
  half-open cubes of side $1/m$, which we call a {\em block-partition} of 
$Q_\bigblock$,
 and $f_\mu \equiv 0$ outside $Q_\bigblock$.
Note that if $\mu$ is both blocked and absolutely continuous, then
$\mu$ is supported by $Q_M$.

\begin{lemm}
\label{lemacblock}
Suppose that $\mu$ is both blocked and absolutely continuous, and that
 $nr_n^d \to \lambda$ as $n \to \infty$, for some
$\lambda \in (0,\infty)$.  Let $t >0$. Then  
\bea
n^{-1} \E[ \zeta_{r_n}( \Po_{nt}) ] \to \int_{\R^d} \rho(\lambda  t f_\mu(x))
 t f_\mu(x) dx
~~~{\rm as}~~ n \to \infty.
\label{1231f0}
\eea
\end{lemm}
\begin{proof}
Let $\eps >0$.
Let 
 $\bigblock$ and $m$  be as in the description of the blocked distribution
 $\mu$.  Fix a cube
$C$ of side $1/m$ in the block-partition of $Q_\bigblock$, and let the constant
value taken by $f_\mu$ on this cube be denoted $b$.
Assume $b  >0$.  Then
\bean
\E [ \zeta_{r_n}(\Po_{nt} \cap C) ] = 
\E [ \zeta (r_n^{-1} (\Po_{nt} \cap C) )] 
= \E [ \zeta (\H_{nr_n^d b  t} \cap  r_n^{-1} C )] , 
\eean
because the restriction of  $\Po_{nt}$ to $C$ is a
homogeneous Poisson process on $C$ with intensity $n b  t$, so 
by the Mapping theorem the restriction of $r_n^{-1} \Po_{nt}$  
to $r_n^{-1}C$ is a 
homogeneous Poisson process on $r_n^{-1} C$ with intensity $n b  t r_n^d$.
Since $nr_n^{d} \to \lambda$,
  by (\ref{meanconvstrong}) from Proposition \ref{propsofrho}(i)
 we have 
as $n \to \infty$ that
\bean
m^d r_n^d \E [ \zeta_{r_n} (\Po_{nt} \cap C) ] 
\to \lambda b  t \rho(\lambda b  t),
 \eean
so that
\bea
n^{-1} \E [ \zeta_{r_n}(\Po_{nt} \cap C) ]  \to  b  t  
 \rho(\lambda b  t) m^{-d} = \rho (\lambda b  t) b  t \Leb(C).  
\label{1231c0}
\eea
Let the cubes in the block-partition be denoted $C_1,\ldots,
C_{(\bigblock m)^d}$. Then 
 for some $b _1,\ldots,b _{(\bigblock m)^d} \in [0,\infty)$
we have $f_\mu = \sum_{i=1}^{( \bigblock m)^d} b _i {\bf 1}_{C_i}$.
 Assume the cubes
are enumerated
 in such a way that for some $k \leq (\bigblock m)^d$
we have that $b _i >0$ for $i=1,\ldots,k$  
and  $b _i =0$ for $i=k+1,\ldots,(Mm)^d$.

By (\ref{fromsub}) from Lemma \ref{lemfromsuper}, we have
$$
\zeta_{r_n}(\Po_{nt}) 
\leq \left( \sum_{i=1}^k \zeta_{r_n}(\Po_{nt} \cap C_i) \right)  
+c_1k,
%\label{1231e0}
$$ 
where $c_1$ is the constant from Property P3.
%\ref{propsub}.}
%
Therefore using (\ref{1231c0}),
%and  (\ref{1231e0}),
 we obtain that 
\bea
\limsup_{n \to \infty} \E[ n^{-1} 
\zeta_{r_n}(\Po_{nt} )  ] \leq
 \sum_{i=1}^k
\rho(\lambda b _{i}t) b _{i}t \Leb (C_i) 
\label{1207a} \\
= \int_{\R^d} \rho(\lambda t f_\mu(x) ) t f_\mu(x) dx .
\label{0103a0}
\eea

For a lower bound, note that for each $n$,
since $\Po_{nt} \subset \cup_{i=1}^k C_i$, by (\ref{fromsuper}) from
Lemma \ref{lemfromsuper}
we have 
$$
\zeta_{r_n} (\Po_{nt} ) \geq \left(
 \sum_{i=1}^k  \zeta_{r_n}(\Po_{nt} \cap C_i) \right)
- c_2 \sum_{i=2}^k | \Po_{nt} \cap \partial_{r_n} C_i| .
$$
Since $\mu$ is absolutely continuous and 
% $n r_n^d \to \lambda$, we have 
$r_n \to 0$, 
 as $n \to \infty$ we have
\bean
 n^{-1} \E[| \Po_{nt}  \cap \partial_{r_n} C_i|] 
= t \mu(\partial_{r_n}C_i) \to 0, ~~~
1  \leq i \leq k.
%\label{0830a0}
\eean
Thus, using (\ref{1231c0}) we obtain that
\bean
\liminf_{n \to \infty} \E[ n^{-1} \zeta_{r_n}(\Po_{nt} )]
 \geq \sum_{i=1}^k
\rho(\lambda b _i t) b _i t \Leb (C_i).
%\nonumber \\
% = \int_{\R^d} \rho(\lambda t f_\mu(x) ) t f_\mu(x) dx.
\eean
Combined with (\ref{1207a}) and (\ref{0103a0}) this gives us (\ref{1231f0}).
\end{proof}

{\bf Remark.}  The next lemma is the only place in the proof
of Theorem \ref{th:therm2}
where we use Property P5. If we assume only P1--P4, but also assume
$\mu$ is absolutely continuous, then we can use 
 Lemma \ref{lemacblock} instead of Lemma \ref{lemblock}, and still
obtain the conclusion of Theorem \ref{th:therm2}.

\begin{lemm}
\label{lemblock}
Suppose $\mu$ is blocked, and $nr_n^d \to \lambda$ as $n \to \infty$,
$\lambda \in (0,\infty)$.  Then  
\bea
\lim_{n \to \infty} \E[ n^{-1} \zeta_{r_n}(\Po_n) ] 
=
\int_{\R^d} \rho(\lambda f_\mu(x) ) f_\mu(x) dx.
\label{1231d}
\eea 
\end{lemm}
\begin{proof}
Let $\eps >0$.
By P5, we can and do choose $\delta \in (0,1)$  
 such that for any
 $\X \subset B_{\delta d}(o)$
with $\delta^{-1} \leq  |\X | < \infty $ we have $\zeta(\X) \leq  \eps |\X|$.

	Let $\mu'$ and (as before) $\mus$ denote the continuous and singular parts of $\mu$,
respectively 
(so $\mu'$ is a measure with density function $f_\mu$). 
Assume that $\bigblock$ (in the description
of the blocked distribution $\mu$) is chosen large enough so that
$\mus(\R^d \setminus Q_\bigblock) < \eps/2$.

As in the proof of \cite[eqn. (7.3)]{Yukich}, provided $m$ (in the definition of a blocked distribution) is
chosen large enough we can assume that
the  support of $\mus$ is contained in the union of two disjoint sets $A$ and
$D$, where $\Leb(A) =0$ and $\mu(A) \leq \eps$, while $D$   is
a union of cubes $C_i$ in the block-partition of $Q_\bigblock$
 with total Lebesgue
measure at most $\eps \delta^{d+1}\lambda/2$.
This is because the support of the restriction of $\mus$ to $Q_M$ has zero Lebesgue outer measure, so we can find a countable collection of dyadic cubes
contained in $Q_M$ with total volume less than $\eps \delta^{d+1} \lambda/2$
 that contains this supporting set. Then we can take
 $D$ to be the union of a sufficiently large finite
 subcollection of these dyadic cubes.

By the superposition theorem for Poisson processes (see for example \cite{LP}),
\bea
\E[\zeta_{r_n}(\Po_n) ]  = \E [\zeta_{r_n} 
( \Po_n^{(1)} \cup \Po_n^{(2)} \cup \Po_n^{(3)} ) ],
\label{0106b}
\eea
where $\Po_n^{(1)}$, $\Po_n^{(2)}$ and
$\Po_n^{(3)}$ are independent Poisson processes with 
intensities $n \mu'|_{D^c}$, $n \mus|_A$ and $n \mu|_D$ respectively.
Here we set $D^c := \R^d \setminus D$ and  for any Borel measure $\nu$ on 
$\R^d$ and Borel $E \subset \R^d$ we write
$\nu|_E$ for the restriction of $\nu$ to $E$.

Set $t= \mu'(D^c)$. If $t >0$, then
 $t^{-1} \mu'|_{D^c}$ is  an absolutely
continuous and blocked probability distribution with
density $t^{-1} f_\mu \cdot{\bf 1}_{D^c}$.
Hence by Lemma \ref{lemacblock} we have 
\bea
\lim_{n \to \infty}( n^{-1} 
\E[ \zeta_{r_n}( \Po_n^{(1)}) ] ) =
  \int_{D^c } 
\rho(\lambda f_\mu(x) ) f_\mu(x) dx.
\label{0924a}
\eea  
If $t=0$, then $\Po_n^{(1)} = \emptyset$, and
(\ref{0924a}) still holds because both sides are zero.

Let the cubes in the block-partition be denoted $C_1,\ldots,
C_{(\bigblock m)^d}$. Then 
 for some $b _1,\ldots,b _{(\bigblock m)^d} \in [0,\infty)$
we have $f_\mu= \sum_{i=1}^{( \bigblock m)^d} b _i {\bf 1}_{C_i}$.
 Assume the cubes are  enumerated 
 in such a way that for some $k \leq k' \leq \ell
 \leq (\bigblock m)^d$
we have that $b _i >0$ for $i=1,\ldots,k'$  
and  $b _i =0$ for $i=k'+1,\ldots,(Mm)^d$,
 while   $D= \cup_{i=k+1 }^{\ell} C_i$.

%Recall that the cubes in the block-partition have side length denoted $1/m$.
Now, for each  $n \in \N$, divide $D$ into smaller cubes (boxes)
 of equal
 side length $ s_n$ where $(1/(ms_n)) \in \N$,
 and $s_n \sim \delta r_n$ as $n \to \infty$, and $s_n \leq \delta r_n$
for all large enough $n$, say for $n \geq n_0$. Denote these boxes
 $B_{n,1},\ldots, B_{n,m_n}$. Then the volume of each one
of these boxes
is asymptotic  to  $\delta^d r_n^d$, and hence to $\delta^d \lambda/n$,
so that $m_n \sim \Leb(D) n/(\lambda \delta^d)$ as $n \to \infty$.
 Since $\Leb(D) \leq 
\eps \delta^{d+1} \lambda/2$, therefore $m_n \leq  \eps \delta n$ for all large enough $n$.

By Lemma \ref{lemfromsuper}, writing 
$\Po_n^{(1,3)}$ for $\Po_n^{(1)} \cup \Po_n^{(3)}$
and noting that $\Po_n^{(3)} \subset D$ while $\Po_n^{(1)} \subset 
\cup_{i=1}^k C_i$, 
we have
\begin{align}
\zeta_{r_n}(\Po_n^{(1,3)})  \leq
 \zeta_{r_n}(\Po^{(1)}_n ) + 
\left( \sum_{i=1}^{m_n}
 \zeta_{r_n}(\Po^{(3)}_n   \cap B_{n,i}) \right) +c_1m_n,
\label{1231e}
\end{align}
 where $c_1$ is the constant from Property P3.
%\ref{propsub}.}
We now show that the expected value of the sum on the right hand side
of (\ref{1231e}) is small. 
For each $i \in \{1,\ldots,m_n\}$, the
 box $B_{n,i}$ is contained in a ball of radius $\delta d r_n$,
so by translation invariance (Property P2) and our choice of $\delta$,
 for any finite
 $\X \subset B_{n,i}$ we have
$
\zeta_{r_n}(\X) \leq  \eps |\X| 
$ if $|\X| \geq \delta^{-1}$.
Moreover, 
 with $K$ as given in  Lemma \ref{lemsmooth}, by that result
 $\zeta_{r_n}(\X) \leq K /\delta$ if $|\X| \leq \delta^{-1}$. Thus
 in all cases
$ \zeta_{r_n}(\X) \leq \eps |\X| + K /\delta$.  
Hence for all large enough $n$,
\bean
\sum_{i=1}^{m_n} \E [ \zeta_{r_n}(\Po_n^{(3)} \cap B_{n,i} ) ] 
\leq 
\sum_{i =1}^{m_n} ( K/\delta + \eps \E [ |\Po_n^{(3)} \cap B_{n,i} | ] ) 
\leq ( K + 1) \eps n.
\eean
%
%and (\ref{gamstardef}),
%$$
% \sum_{i=k+1}^\ell
% \zeta_{r_n}(\Po_n^{(3)}  \cap C_i)  
%=\sum_{i=k+1}^{\ell} \zeta(r_n^{-1}(\Po_n^{(3)}  \cap C_i))  
%\leq (\ell -k) \zeta^* (Q_{1/(r_n m)}).
%$$
%By (\ref{wellknown2}), $m^d r_n^d \zeta^*(Q_{1/(r_n m)}) \to
%  \overline{\zeta}
%$ as  $n \to \infty$, and
%since $ nr_n^d \to  \lambda  $,
%it follows that
%\bean
%\limsup_{n \to \infty} n^{-1} 
%\sum_{i=k+1}^{\ell}
%\E[ \zeta_{r_n}(\Po_n^{(3)}  \cap C_i) ]  
%\leq 
%\lambda^{-1} (\ell -k)
%%\sum_{i=k+1}^{\ell} 
%m^{-d} \overline{\zeta}
%%\\
%\leq  \lambda^{-1} \overline{\zeta} \eps,
%\eean
%because $(\ell -k)m^{-d} $ is the Lebesgue measure of $D$.

Also  $\rho(r) \geq 0$ for all $r \geq 0$, by (\ref{meanconv}).
Therefore using (\ref{0924a}) and  
(\ref{1231e}), and the assumption $\delta <1$, we obtain that 
\bea
\limsup_{n \to \infty} \E[ n^{-1} 
\zeta_{r_n}(\Po_n^{(1,3)} )  ] \leq
\int_{Q_M} \rho(\lambda f_\mu(x) ) f_\mu(x) dx
 + ( K+1 +  c_1) \eps.
\label{0103a}
\eea

For a lower bound, we start with $\Po_n^{(1)}$. 
By (\ref{meanconv}) and  Lemma \ref{lemsmooth}  we have
$\rho(a) \leq K$ for all $a >0$.  Since $\mu$ is blocked
we have that $f_{{\rm max}} := \sup_{x \in \R^d} f_\mu(x) <\infty$. 
Therefore
 $ \int_{D} \rho(\lambda f_\mu(x) ) \lambda f_\mu(x) dx
\leq K \lambda f_{{\rm max} } \Leb(D) \leq 
K \lambda^2 f_{{\rm max}} \eps$.  
%By (\ref{0922a}) and the fact that  
%since $\Leb(D) < \eps$,
% we have
% $ \int_{D} \rho(\lambda f_\mu(x) ) \lambda f_\mu(x) dx
%\leq   \ozeta \, \eps$.
 Hence 
by (\ref{0924a}), since $f_\mu$ is supported by $Q_\bigblock$,
\bea
\liminf_{n \to \infty} \E[ n^{-1} \zeta_{r_n}
(\Po_n^{(1)}  ) ]
\geq
  \int_{Q_\bigblock } \rho(\lambda f_\mu(x)) f_\mu(x) dx
  - K \lambda^2 f_{{\rm max}} \eps.
 %\lambda^{-1} \overline{\zeta}
 %\eps.
\label{0103b2}  
\eea 

Also since $\Po_n^{(3)} \subset D$ and
 $\Po_n^{(1)} \subset \cup_{i=1}^k C_i$ which is disjoint from $D$,
 by (\ref{fromsuper}) from Lemma \ref{lemfromsuper}, and the non-negativity
 of $\zeta$,
\bean
\zeta_{r_n}( \Po_n^{(1,3)}) - \zeta_{r_n}(\Po_n^{(1)})  
\geq \zeta_{r_n}(\Po_n^{(3)}) - c_2 | \Po_n^{(1)} \cap 
 \partial_{r_n}
( \cup_{i=1}^k
 C_i ) |
\\
\geq
 %\zeta_{r_n}(\Po_n^{(3)} ) 
- c_2 \sum_{i=1}^k |\Po_n^{(1)} \cap \partial_{r_n} C_i |.
\eean
Therefore taking expectations we obtain that
$$
\liminf_{n \to \infty} n^{-1} \E[
\zeta_{r_n}( \Po_n^{(1,3)}) - \zeta_{r_n}(\Po_n^{(1)}) ] 
\geq -c_2 \lim_{n \to \infty} \sum_{i=1}^k \mu'( \partial_{r_n}  C_i) =0
$$  
Hence by (\ref{0103b2}),
\bea
\liminf_{n \to \infty} \E[ n^{-1} \zeta_{r_n}
(\Po_n^{(1,3)}  ) ]
\geq
  \int_{Q_\bigblock } \rho(\lambda f_\mu(x)) f_\mu(x) dx
  -
% \lambda^{-1} \overline{\zeta}
K \lambda^2 f_{{\rm max}} 
 \eps.
\label{0103b}  
\eea

Also by (\ref{0106b}) and Lemma \ref{lemsmooth}, 
for all $n$ we have
\bean
%n^{-1} 
|\E[ \zeta_{r_n}( \Po_n )] - \E[\zeta_{r_n}( \Po_n^{(1,3)} )
% \cup \Po_n^{(3)}))  
]|   \leq 
K \E[ | \Po_n^{(2)}|]
\leq K n \eps.
\eean
Therefore using (\ref{0103a}) and (\ref{0103b}),
since  $\limsup |a_n| = \max (\limsup(a_n),-\liminf(a_n))$
for any real-valued sequence $(a_n)$,
we obtain that
\bean
\limsup_{n \to \infty} \left| n^{-1} \E [
 \zeta_{r_n}( \Po_n)]  - 
  \int_{Q_\bigblock} \rho(\lambda f_\mu(x)) f_\mu(x) dx
\right| 
 \leq ( 
2K +1  + c_1 + K \lambda^2 f_{{\max}} )
\eps, 
\eean
and since $\eps >0$ is arbitrary this yields
\bean
\lim_{n \to \infty} \E[ n^{-1} \zeta_{r_n}(\Po_n) ] 
=
\int_{Q_\bigblock} \rho(\lambda f_\mu(x) ) f_\mu(x) dx
=
\int_{\R^d} \rho(\lambda f_\mu(x) ) f_\mu(x) dx,
\eean
which is (\ref{1231d}), as required.
\end{proof}

We now drop the restriction, in the last lemma, to blocked
density functions.
\begin{lemm}
\label{lemunblock}
  Suppose 
$nr_n^d \to \lambda$
as $n \to \infty$, for some $\lambda \in (0,\infty)$. Then
 (\ref{1231d}) holds.
\end{lemm}
\begin{proof}
Let $\eps >0$. 
As in the proof of Lemma 7.3 of \cite{Yukich}, we can find a blocked
 distribution $\nu$ with the same singular part as $\mu$,
such that 
\bean
\int |f_{\nu}(x) -f_{\mu}(x)|dx < \eps,
\eean
where we now write $\int$ for $\int_{\R^d}$. 
Also as in that proof, we can find a pair of coupled random variables
$(X,Y)$ such that $X$ has distribution $\mu$, $Y$ has distribution
$\nu$, and $\Pr[X \neq Y] \leq  \eps$. Then taking 
$(X_i,Y_i)_{i=1,2,3,\ldots}$
 to
be a sequence of i.i.d. coupled pairs with the distribution
of $(X,Y)$, we may consider coupled Poisson
point processes $\Po_n := \{X_1,\ldots,X_{N_n}\}$
and $\Q_n:= \{Y_1,\ldots,Y_{N_n}\}$. 
By smoothness (Lemma \ref{lemsmooth}),
\bea
\E [ n^{-1} |\zeta_{r_n}(\Po_n) - \zeta_{r_n}(\Q_n) | ]
\leq K \E[n^{-1}|\Po_n \triangle \Q_n|] \leq 2K \eps.
\label{0101e}
\eea 
By Lemma \ref{lemblock}, for large enough $n$ we have
\bea
\left| \E [n^{-1} \zeta_{r_n}(\Q_n) ] -
 \int \rho(\lambda f_\nu(x) ) f_\nu(x) dx \right| < \eps.
\label{0101f}
\eea 
Moreover, by the Lipschitz continuity in $\lambda$ 
of $\lambda \rho(\lambda)$ (see Proposition \ref{propsofrho}(ii)),
 \bea
\left| \int \rho(\lambda f_\nu(x) ) f_\nu(x) dx  -
\int \rho(\lambda f_\mu(x) ) f_\mu(x) dx   \right|
\nonumber \\
 =
\lambda^{-1}  \left| \int (\lambda f_\nu(x) \rho (\lambda f_\nu(x) )
 - \lambda f_\mu(x) \rho(\lambda f_\mu(x))  ) dx  
\right|
\nonumber \\
  \leq  \lambda^{-1} K \int |\lambda f_\nu(x) -\lambda f_\mu(x)|dx
\leq   K\eps.
\label{0101d}
 \eea
Combining (\ref{0101e}), (\ref{0101f}), and (\ref{0101d}), we obtain
for all large enough $n$ that 
$$
\left| \E [n^{-1} \zeta_{r_n}(\Po_n) ] - \int  \rho (\lambda f_\mu(x))  f_\mu(x) dx
\right|  \leq
\eps( 1 + 3K ),
$$
and since $\eps$ is arbitrarily small, (\ref{1231d}) follows.
\end{proof}

\begin{proof}[Proof of Theorem \ref{th:therm2}]
Suppose $n r_n^d \to t \in (0,\infty)$. By 
Lemma \ref{lemunblock}, $n^{-1} \E[\zeta_{r_n}(\Po_n)]$ converges
to
$\int_{\R^d} \rho(t f_\mu(x)) f_\mu(x) dx$, as $n \to \infty$.

Recalling from Section \ref{secprelim} our
 coupling of $\Po_n$ and $\X_n$,  and
that $N_n:= |\Po_n|$, by 
Lemma \ref{lemsmooth} we have 
$$
\E[ |\zeta_{r_n}(\Po_n)  - \zeta_{r_n}(\X_n)| ] \leq
K\E[|\Po_n \triangle \X_n|]  = K \E[|N_n -n |]= o(n),
$$ 
and therefore 
\bea
n^{-1} \E[ \zeta_{r_n}( \X_n) ] \to 
\int_{\R^d} \rho(t f_\mu(x)) f_\mu(x) dx
~~~{\rm as}~~ n \to \infty.
\label{1231f}
\eea

The proof of (\ref{0101c}) starting from (\ref{1231f})  
is completed by martingale difference methods as in the proof of
 Theorem~\ref{0227a} starting from (\ref{0228c}),
 with $(\cU_n, \cV_n)$ therein replaced by $\cX_n$ and
 $(\cU_{n,i}, \cV_{n,i})$  replaced by $\cX_{n,i}$, which
 is obtained from $\cX_n$ by replacing $X_i$ by an independent copy $X_0$.
Due to the similarity to  
 the corresponding part of the proof of Theorem~\ref{0227a}, we do
not give the details this time.

By Lemma \ref{lemsmooth}, for all $\eps >0, n \in \N$ we have
$$
\Pr[|\zeta_{r_n}(\Po_n) - \zeta_{r_n}(\X_n) | > \eps n]
\leq \Pr[| N_n- n | > \eps n/K ], 
$$ 
which is summable in $n$ by a Chernoff bound (e.g. \cite[Lemma 1.2]{RGG}).

Using this, and the complete convergence in (\ref{0101c}),
gives us  the complete convergence in (\ref{0302a}).
Also $\E[(n^{-1} (\zeta_{r_n}(\Po_n) - \zeta_{r_n}(\X_n)) )^2] \leq
\E[n^{-2} K^2 |N_n-n|^2]$, which tends to zero, so
using the $L^2$ convergence in (\ref{0101c}),
we have the $L^2$ convergence in (\ref{0302a}).
\end{proof}

\begin{proof}[Proof of Lemma \ref{lemdet}]
By (\ref{gamstardef}) and P3, $\zeta^* (A \cup A') +c_1
\leq (\zeta^*(A) +c_1) + ( \zeta^*(A') + c_1)$ for any 
$A, A' \subset \R^d$ (not necessarily disjoint).
Also  $\zeta^*(A) \leq \zeta^*(A')$ whenever $A \subset A'$.

Setting $k:= \kappa(Q_2)$, we can cover $Q_1$ by $k$ balls of
radius $1/2$, and repeatedly using the above
subadditivity and monotonicity properties we obtain 
for all $A \subset Q_1$ that $\zeta^*(A) +c_1 \leq 
k(\zeta^*(B_{1/2}(o)) +c_1)$, which is finite by  Property P$5'$.
Hence taking $f(A) = \zeta^*(A) + c_1$
gives a functional satisfying all the conditions of Lemma \ref{lemsubadd},
and that result gives us
\eqref{wellknown2}.

	By (\ref{gamstardef}),
	$\E[\zeta( \cH_{\lambda,s})] \le \zeta^*(Q_{s})$,
%Given 
	$s >0$. 
%Therefore 
	Then by (\ref{meanconv}) and (\ref{wellknown2}), 
we have (\ref{0922a}).
\end{proof}

\begin{proof}[Proof of Theorem \ref{th:lambig}]
Since we assume $\zeta(\cdot) \geq 0$,
we have that
 $\ozeta \geq 0$ by 
(\ref{wellknown2}), and $\rho(\lambda) \geq 0$ for all $\lambda \in \R_+$
by (\ref{meanconv}). If $\ozeta =0$, then by (\ref{0922a}) $\rho(\lambda)=0$
for all $\lambda$ so the result (\ref{eq:lambig}) holds. 
Therefore we may assume without loss of generality that $\ozeta >0$.

Let $\eps \in (0,1)$. Choose $ s_0 \in (0,\infty)$ with 
$(s_0/(s_0+2))^d > 1-\eps$ and $s_0^{-d} \zeta^*(Q_{s_0}) > (1-\eps)
\ozeta$. Then choose a finite set $\cX \subset Q_{s_0}$
with $ \zeta(\cX) > (1- \eps) \zeta^*(s_0)$.

Choose $\delta  \in (0,1)$ such that 
 $\inf[\{\|x-y\| -1 : x,y \in \cX\} \cap (0,\infty)] > \delta$ 
(using the convention $\inf[\emptyset]:=  +\infty$ if needed); this
can be done because
 the infimum is over a finite set of strictly
positive numbers.   Define the re-scaled
point set $\cX' := (1+\delta)^{-1}\cX$; then $G(\cX',1)$
is isomorphic to $G(\cX,1)$ so $\zeta(\cX') = \zeta(\cX)  $ by assumption
P8, and $\cX'$ has no point on the boundary of $Q_{s_0}$, and 
$\|x-y\| \neq 1$ for all $x,y \in \cX'$.

Enumerate the elements of
 $\cX'$ as $x_1,\ldots,x_k$, say. By the preceding conclusion there
exists $\delta' >0$ such that the balls $B_i := B_{\delta'}(x_i)$
$1 \leq i \leq k$, are disjoint, are all contained in $Q_i$, and and have
 the property that for any $\cY= \{y_1,\ldots,y_k \}$ with 
$y_i \in B_i$ for each $i$ we have 
that $G(\cY,1)$ isomorphic
to $G(\cX',1)$ and hence $\zeta(\cY) = \zeta(\cX' ) = \zeta(\cX)$ by P8.
By P6 (upwards monotonicity in  $\cX$),
\bean
%\liminf_{\lambda \to \infty} 
\Pr[ \zeta(\cH_\lambda \cap Q_{s_0}) > (1-\eps)
\zeta^*(Q_{s_0}) ]
\geq \Pr \left[ \cap_{i=1}^k \{\cH_\lambda \cap B_i \neq \emptyset \} \right]
\\
 \to 1 {\rm ~~~~ as ~} \lambda \to \infty.
\eean
Hence, using also the nonnegativity 
of $\zeta$, we have that
$$
	\liminf_{\lambda \to \infty} \E[ \zeta(\cH_{\lambda,s_0}) ]  > 
(1-\eps)^2 \zeta^*(Q_{s_0}) .
$$
Pick $\lambda_0 \in (0,\infty )$ such that
\bea
\E[ \zeta(\cH_{\lambda ,s_0}) ]  > 
(1-\eps)^2 \zeta^*(Q_{s_0}) , ~~~~ \forall ~ \lambda \geq \lambda_0.
\label{1230c}
\eea

Given $s>0$, write $s= n(s_0+2) +t$ with $0 \leq t < s_0+2$.
Let
$Q'_1, \ldots, Q'_{n^d}$ be a collection of cubes of side $s_0$,
that are disjoint and distant at least 2 from each other, and
contained in $Q_s$. Then using P6 followed repeatedly by P4, we have
$$
\zeta( \cH_\lambda \cap Q_s) \geq \zeta 
\left( \cH_\lambda \cap \left( \cup_{i=1}^{n^d} Q'_i \right) \right)
\geq \sum_{i=1}^{n^d} \zeta ( \cH_\lambda \cap Q'_i ).
$$
Therefore by (\ref{1230c}) for fixed $\lambda \geq \lambda_0$, provided $s$ (and hence $n$) is large enough
\bean
s^{-d} \E[ \zeta(\cH_\lambda \cap Q_s) ] \geq \frac{ n^d \E[\zeta(\cH_\lambda \cap
Q_{s_0} ) ] }{ (n+1)^d (s_0+2)^d} 
\\
\geq \left(\frac{n}{n+1} \right)^d \left( \frac{s_0}{s_0+2} \right)^d
\left( \frac{(1- \eps)^2 \zeta^*(Q_{s_0}) }{s_0^d} \right)
> (1-\eps)^{5}\,\ozeta.
\eean
Using (\ref{meanconv}), 
this shows that $\lambda \rho(\lambda) \geq (1-\eps)^5 \, \ozeta $
for $\lambda \geq \lambda_0$,
so that $\liminf_{\lambda \to \infty} (\lambda \rho(\lambda)) \geq \ozeta$.
Combined with (\ref{0922a}), this gives us the result.
\end{proof}

\begin{proof}[Proof of Theorem \ref{th:densegen}]
Assume that $r_n \to 0$ and $n r_n^d \to \infty$ as $n \to \infty$, and
that
% $\mu = \mu_U$
$\mu$ is absolutely continuous with density $f_\mu$ chosen in
such a way that the set $A:= f_\mu^{-1}((0,\infty))$ is
Riemann measurable.
 Then,
given $\eps >0$, we can find $a >0$ and $k \in \N$, such that we can find a covering of the set $A$
 by disjoint cubes $C_1,\ldots,C_k$
of side $a$, with $k a^d < \Leb(A) + \eps$. 
%such a way that $f_\mu^{-1}((0,\infty))$ is
Then by (\ref{fromsub}) and (\ref{gamstardef}),  almost surely
\bean
%r_n^d
 \zeta_{r_n} (\cX_n) \leq
\left( \sum_{i=1}^k 
%r_n^d
 \zeta_{r_n} (\X_n \cap C_i) \right) + c_1 k 
\
\leq
k 
% r_n^d
 \zeta^*(Q_{r_n^{-1}a}) + c_1 k,
\eean
so that by (\ref{wellknown2}),
\bea
\limsup_{n \to \infty} (r_n^{d}  \zeta_{r_n}(\X_n)) \leq k a^d \, \ozeta
\leq (\Leb(A) + \eps) \,  \ozeta.
\label{1209a}
\eea
For an inequality the other way, let $\lambda \in (0,\infty)$ and let
 $(m(n), n \in \N)$ be an $\N$-valued
sequence with $m(n) \leq n$ for all $n$ and with $m(n) r_n^d \to \lambda$
as $n \to \infty$. Then
by P6, almost surely 
$$
r_n^d \zeta_{r_n}( \cX_n) 
\geq r_n^d\zeta_{r_n} (\X_{m(n)}) =
(1+ o(1)) (\lambda/m(n)) \zeta_{r_n}( \cX_{m(n) }),  
$$
so by Theorem \ref{th:therm2}, almost surely
\bea
\liminf_{n \to \infty} ( r_n^d \zeta_{r_n}( \cX_n) )
 \geq \lambda \int_{\R^d} \rho(\lambda f_\mu(x) ) 
f_\mu(x)dx.
\label{1228}
\eea  
By Theorem \ref{th:lambig} and monotone convergence 
 (which is applicable by Proposition \ref{propsofrho} (iv)),
we have 
$ \int_A \lambda f_\mu(x) \rho(\lambda f_\mu(x)) dx \to \ozeta \, \Leb(A)$
 as $\lambda \uparrow \infty$.
Since (\ref{1228}) holds for arbitrarily large $\lambda \in (0,\infty)$, 
we obtain that
$
\liminf_{n \to \infty} ( r_n^d \zeta_{r_n}( \cX_n) ) \geq \ozeta \, \Leb(A)$,
and combined with (\ref{1209a}) this gives us the result.
\end{proof}

\section{Proofs for domination number at high density}\label{secdom}

%{\color{Cyan}
%\section{Proofs for the dense limit}
\label{secsubd}
\allco

In this section we shall prove Theorems \ref{th:onemore} and
\ref{thsubcond}.
Throughout  this section we take $\zeta(\cX) = \gamma(G(\X,1))$,
the domination number of $G(\X,1)$, for all finite 
$\X \subset \R^d$. Therefore
$\zeta_r(\X)$ is equal to $\gamma(G(\X,r))$.
% Also we take $B= B_1(o)$, the Euclidean unit ball in $d$ dimensions. 

Throughout  this section we assume that
 a sequence $(r_n)_{n \geq 1}$ is given, taking values in $(0,\infty)$ and
 satisfying $r_n \to 0  $ and $n r_n^d \to \infty$ as $n \to \infty$, but otherwise arbitrary. 
This is part of the hypothesis of Theorem \ref{thsubcond} but we

\subsection{Deterministic preliminaries}

\begin{lemm}[Smoothness of the domination number]
\label{lemA1C}
There is a constant $K$ such that if
 $\X, \cY$ are finite subsets of $\R^d$, and $r > 0$,
% and $x \in \R^d \setminus \X$.
 then 
\bea
|\zeta_r(\cY) - \zeta_r(\X) | \leq K |\cY \triangle \cX|.
\label{1221j}
\eea  
\end{lemm}
\begin{proof}
This follows from Lemma \ref{lemsmooth}; we already checked in
Section \ref{subsecdom} that our current choice of $\zeta$
satisfies P1--P4.
\end{proof}

For $s >0$ set $\oQ_s := [-s/2,s/2]^d$,
the closure of $Q_s$. The sequence $(r_n)_{n \in \N}$ was specified already.  
Also we fix a constant $\delta \in (0,1/4)$. Recall the definitions
(\ref{kappadef}) and (\ref{psidef}) of $\kappa(\cdot)$ 
and $\okappa$, respectively.

\begin{lemm}
\label{lemcov}
%Let $\delta \in (0,1/4)$.
There exists a sequence $(\LL_n)_{n \in \N}$ 
of finite sets with $\LL_n \subset \oQ_{r_n^{-1}}$ 
and $Q_{r_n^{-1} } \subset \cup_{x \in \LL_n} B_{1- \delta}(x)$
for each $n$, such
that $k_n := |\LL_n|$ satisfies
\bea
(1-  \delta)^{-d} \okappa
\leq
\liminf_{n \to \infty} (r_n^d k_n) \leq 
\limsup_{n \to \infty} (r_n^d k_n) \leq 
(1- 4 \delta)^{-d} \okappa,
% ~~~{\rm as}~ n \to \infty.
\label{0102a}
\eea
and moreover with $\|x -y\| > 3 \delta$ for all $x,y \in \LL_n$ such that
$x \neq y$.
\end{lemm}
\begin{proof}
Set
$
k^0_n := \kappa(Q_{(1-4 \delta)^{-1} r_n^{-1}}), 
$
and let $\LL_n^0 \subset \R^d$ 
with $Q_{r_n^{-1}} \subset \cup_{x \in \LL_n^0} B_{1-4 \delta}(x)$ and
 $|\LL_n^0 | = k^0_n$.
We may assume without loss of generality that all elements of $\LL_n^0$ lie
in $\oQ_{r_n^{-1}}$, since if $x \in \LL_n^0 \setminus \oQ_{r_n^{-1}}$,  
and $y$ is the closest  point of $\oQ_{r_n^{-1}} $ to $x$, then
$\|z-y\| \leq \|z-x \|$ for all $z \in Q_{r_n^{-1}}$, so we could replace the
point $x$  in $\LL^0_n$ by $y$.

Let $\LL_n \subset \LL_n^0$ be a 
maximum independent set of $G(\LL_n^0, 3 \delta)$, i.e.
an  independent set of $G(\LL_n^0, 3 \delta)$
with $|\LL_n| = \alpha(G( \LL_n^0,3 \delta))$, where $\alpha(\cdot)$
was defined in Section \ref{secfurther}.
Then $\LL_n$ is a dominating set for 
$G( \LL_n^0,3 \delta)$; otherwise, a further element
of $\LL_n^0 \setminus \LL_n$ could be added to $\LL_n$ without
it ceasing to be an independent set, contradicting the 
maximality of $\LL_n$.

Then $Q_{r_n^{-1} } \subset \cup_{x \in \LL_n} B_{1- \delta}(x)$. Indeed,
for each $x \in \LL^0_n \setminus \LL_n$ there exists $y \in \LL_n$
with $\|x-y\| \leq 3 \delta$, and therefore by the triangle
inequality $B_{1-4 \delta}(x) \subset B_{1-  \delta}(y)$. 

Set $k_n := |\LL_n|$. 
By (\ref{psidef}),
$
r_n^d k^0_n \to (1- 4 \delta)^{-d} \okappa ~{\rm as}~ n \to \infty.
$
Since $k_n \leq k_n^0$ this gives the last inequality of
(\ref{0102a}). Since 
$Q_{r_n^{-1} } \subset \cup_{x \in \LL_n} B_{1- \delta}(x)$,
by (\ref{kappadef}) 
we have $k_n \geq \kappa(Q_{(1-\delta)^{-1} r_n^{-1} }) $,
and hence by (\ref{psidef}), the first inequality of (\ref{0102a}).
\end{proof}

Given $n \in \N$,
enumerate $\LL_n$ as $\LL_n = \{x_{n,1},\ldots,x_{n,k_n}\}$.
Choose $n_0 \in \N$ such that $r_n \leq 1 $ for all $n \geq n_0$. Then
 for all $n \geq n_0$ and each $i \in [k_n]$, since
$x_{n,i} \in \overline{Q}_{r_n^{-1}}$, the ball $B_\delta (x_{n,i}) 
\cap Q_{r_n^{-1}}$ contains a cube of side $\delta/d$ with a corner
at $ x_{n,i}$, and hence
\bea
\Leb (B_{\delta}(x_{n,i}) \cap Q_{r_n^{-1}} ) \geq (\delta/d)^d,
~~~~ \forall ~ n \geq n_0, i \in [k_n]. 
\label{0101a}
\eea
% For $1 \leq i \leq k_n$
%define the ball $B_{n,i} := B_{\delta}(x_{n,i})$.
Given  
 finite $\X \subset Q_{r_n^{-1}}$,
define the set of `good' indices 
\bea
I_n (\X) :=\{i \in [ k_n]: \X \cap B_\delta(x_{n,i}) \neq \emptyset \}.
\label{Indef}
\eea
Set $K_0 := \kappa(B_2(o))$.

\begin{lemm}
\label{lemub}
For all $n \in \N$ and all finite $\X \subset Q_{r_n^{-1}}$,
 \bea
\zeta(\X)  \leq k_n + K_0|[k_n] \setminus I_n(\X)|.
\label{1220b}
\eea
\end{lemm}
\begin{proof}
Let $n \in \N$. For each $i \in [k_n]$,
cover the ball $B_{1}(x_{n,i})$ by
balls $B_{n,i,1}, \ldots,B_{n,i,K_0}$ of radius $1/2$.
 We now define a set $\cA \subset \X$ as follows. 

For each `good' index $i \in I_n(\X)$, let $y_{n,i}$ be
the first element of $\X \cap B_\delta(x_{n,i})$
in the lexicographic ordering.  If
 $i \in [k_n] \setminus I_n(\cX)$  
(so $i$ is a `bad' index), then
 for each $j \in [K_0]$ such that 
 $\cX \cap B_{n,i,j} \neq \emptyset$,
let  $z_{n,i,j}$ be the
first element
 of $\cX \cap B_{n,i,j}$ in the lexicographic ordering.
Set
$$
\cA := \{y_{n,i}: i \in I_n(\X)\} \cup \{z_{n,i,j}:i \in [k_n] \setminus
I_n(\X), j \in [K_0], \X \cap B_{n,i,j} \neq \emptyset \}.
$$

We assert that $\cA$ is a dominating set for $G(\X,1)$.
Indeed, if $i$ is a good index then the set $B_{1}(y_{n,i})$ 
contains the whole of the ball $B_{1-\delta}(x_{n,i})$.
If $i$ is a bad index, then for each 
 $x \in \cX \cap B_{1-\delta}(x_{n,i})$,
 choosing $j \in [K_0]$ such that $x \in B_{n,i,j}$ we have
$x \in B_{1}(z_{n,i,j})$. Since 
$\X \subset Q_{r_n^{-1}}$, every point of $\X$ lies in
at least one of the balls $B_{1- \delta}(x_{n,i})$, $i \in [k_n]$.
Therefore $\cA$ is a dominating set as asserted. 
Moreover, $|\cA| \leq k_n + K_0|[k_n] \setminus I_n(\X)|$, 
and (\ref{1220b}) follows.
\end{proof}

Now we derive a bound the other way. Given $n \in \N$,
partition  $Q_{r_n^{-1}}$ into cubes of
side $\delta_n := r_n^{-1}/\lceil r_n^{-1}/\delta  \rceil$, denoted 
$S_{n,1},S_{n,2},\ldots,S_{n,\ell_n}$. Note that 
$\delta_n \leq \delta$ and  $\delta_n \to \delta$
 as $n \to \infty$, so that
\bea
\ell_n \sim (r_n^{-1}/\delta )^d ~~~~~{\rm as~~} n \to \infty.
\label{1230d}
\eea 
%
%Also, set $K_1:= \kappa(Q_1)$.
For $ i \in [\ell_n]$,
%let $y_{n,i,1},\ldots,y_{n,i,K_1} \in \R^d$ 
let $v_{n,i}$ denote the centre of the cube $S_{n,i}$. 
%such that $S_{n,i} \subset \cup_{j=1}^{K_1} B_1(y_{n,i,j})$. 
% with
%the centre of cube $S_{n,i}$ denoted $y_{n,i}$
%for $i \in [\ell_n]$. 
Given finite $\X \subset Q_{r_n^{-1} } $, set
$$
J_n(\X) : = \{ i \in [\ell_n]: \X \cap S_{n,i} \neq \emptyset \}. 
$$
\begin{lemm}
\label{lemlow1}
Let $n \in \N$ and let $\X \subset Q_{r_n^{-1}} $ be finite. Then
\bea
\zeta(\X) \geq r_n^{-d} (1+d \delta)^{-d} \okappa 
- |[\ell_n] \setminus J_n(\cX)|. 
\label{1221e}
\eea
\end{lemm}
\begin{proof}
Suppose $\cA \subset \X$
is a dominating set for $G(\X,1)$.
Then for each $i \in J_n(\X)$ we have
for some $x \in \cA$ and $y \in \X \cap S_{n,i}$   that 
$y \in B_1(x)$, and hence $S_{n,i} \subset B_{1+ d \delta}(x)$.
Therefore
$$
\cup_{i \in J_n(\cX)} S_{n,i} \subset \cup_{x \in \cA}  B_{1+ d \delta}(x).
$$
%Also, clearly $S_{n,j} \subset B_{1+2 d \delta}(y_{n,j})$
%for each $j \in [\ell_n]$.  
Also, since the cube $S_{n,i}$ has 
side at most $\delta$, we have $S_{n,i} \subset B_{1+ d \delta}(v_{n,i})$,
for $1 \leq i \leq \ell_n$.
Hence the union of balls of radius
 $1+ d \delta$ centred on the points of the set
 $$
\cA \cup \{ v_{n,i} : i \in [\ell_n] \setminus J_n(\X)
%, j \in [K_1]
  \}
$$
contains the whole of $Q_{r_n^{-1}}$. Therefore by the definition
(\ref{psidef}) of $\okappa$,
we have
$$
(1+ d \delta)^d r_n^d (|\cA | +  |[\ell_n] \setminus J_n(\X)|)  
\geq
 \okappa,
$$
for all choices of $\cA$. Hence
% for large enough $n$, if also 
%$(1+\delta)^d |[\ell_n] \setminus J_n| \leq \delta r_n^{-d}$ 
%we have
 (\ref{1221e}) holds. 
\end{proof}

\subsection{Proof of Theorem \ref{th:onemore}}

\begin{proof}[Proof of Theorem \ref{th:onemore}]

Let $(r_n)_{n \geq 1}$ and $\delta \in (0,1/4)$ be as specified
already in this section.
Let $\lambda >0$.  
Let $\LL_n := \{x_{n,1},\ldots,x_{n, k_n}\}$
 be as described
in Lemma \ref{lemcov}.
 For finite $\X \subset Q_{r_n^{-1}}$,
define the set of `good' indices $I_n (\cX)$ by (\ref{Indef}). 
Then by Lemma \ref{lemub},
\bean
\zeta(\H_\lambda \cap Q_{r_n^{-1}}) 
\leq k_n + K_0 |[k_n] \setminus I_n (  \H_\lambda \cap Q_{r_n^{-1}})|,
\eean
 and $\limsup_{n \to \infty} (r_n^{d} k_n) \leq (1- 4 \delta)^{-d} \okappa$
by (\ref{0102a}).
Then by (\ref{0101a}), 
$$
\E[|[k_n] \setminus I_n( \H_\lambda \cap Q_{r_n^{-1}}) | ]
\leq
% (k_n - k'_n) +
 k_n  \exp(- \lambda  (\delta/d)^d).
$$ 
Therefore, using also (\ref{meanconv}), 
 we have
\bean
\lambda \rho(\lambda) = \lim_{n \to \infty}(r_n^{d} 
\E [\zeta(\H_\lambda \cap Q_{r_n^{-1}} ) ])
 \leq
(1- 4 \delta)^{-d} \okappa ( 1 +  \exp(-\lambda  (\delta/d)^d)),
\eean
so that
\bea
\limsup_{\lambda \to \infty} (\lambda \rho(\lambda)) \leq (1- 4 \delta)^{-d} \okappa. 
\label{0107a}
\eea

By Lemma \ref{lemlow1} we have
\bean
\zeta(\H_\lambda \cap Q_{r_n^{-1}}) 
%= \gamma_{r_n}(r_n \H_\lambda \cap Q_1)  
\geq (1+  d \delta)^{-d} r_n^{-d} \okappa -
  | [\ell_n ] \setminus J_n( \H_\lambda \cap Q_{r_n^{-1}} )|,  
\eean
and $\ell_n \sim (r_n^{-1} /\delta)^d$ by (\ref{1230d}).
For fixed $\lambda$ we have as $n \to \infty$ that
$$
 \E[ | [\ell_n ] \setminus J_n( \H_\lambda \cap Q_{r_n^{-1}})| ] 
\sim (r_n^{-1}/\delta)^d \exp(-\lambda \delta^d ).
$$  
Therefore, 
$$
\lambda \rho(\lambda) = \lim_{n \to \infty}
%(r_n^{d} \E [\zeta(\H_\lambda \cap Q_{r_n^{-1}}) ]) \geq 
	(r_n^{d} \E [\zeta(\H_{\lambda ,r_n^{-1}}) ]) \geq 
(1+ d \delta)^{-d} \okappa -  \delta^{-d} \exp(- \lambda \delta^d ),
$$
so that
$$
\liminf_{\lambda \to \infty} (\lambda \rho(\lambda)) \geq 
(1+ d \delta)^{-d} \okappa. 
$$
Combining this with (\ref{0107a}), since $\delta  \in (0,1/4)$ is arbitrary, 
gives us Theorem~\ref{th:onemore}.
\end{proof}

\subsection{Proof of Theorem \ref{thsubcond}}

Assume in this subsection that $\mu = \mu_U$.
We shall prove Theorem \ref{thsubcond} using separate
arguments according 
to whether $n r_n^d $ grows faster or more slowly than $n^{1/8}$.
Set
\bea
\cN: = \{ n \in \N:r_n^d \geq n^{-7/8}\}. 
\label{cNdef}
\eea
Let $\delta \in (0,1/4)$, and
let $\LL_n$ and $k_n$ be as described in Lemma \ref{lemcov}.
\begin{lemm}
\label{lemlimsup1}
Let $\eps >0$.
% and now assume $  n^{-1/d} \ll r_n \ll 1$. 
Then
\bea
\sum_{n \in \cN}
 \Pr [  r_n^d \zeta_{r_n}(\X_n) > \okappa + \eps] < \infty .
\label{eqlimsup1}
\eea
\end{lemm}
\begin{proof}
Assume $\cN$ is infinite (otherwise (\ref{eqlimsup1}) is trivial).
%Note that for any $a \geq 1$ and $x \in \oQ_a$ we can find a cube
%of side $\delta/d$ with one corner at $x$ that is contained in 
%$ Q_a \cap B_\delta(x)$.
For all $n \in \cN$  and $i \in [k_n]$, recalling the definition
of $I_n(\X)$ at (\ref{Indef}), and using (\ref{0101a}), if $n \geq n_0$ we have 
$$
\Pr[ i \notin I_n(r_n^{-1} \X_n)] \leq (1 - (\delta/d)^d r_n^d)^n \leq
 \exp(-n 
(\delta/d)^d r_n^d) \leq \exp ( - (\delta/d)^d  n^{1/8}).
$$
Since $r_n^{-d} \ll n$ so $k_n \ll n$ by (\ref{0102a}),
 by the union bound we have 
for all large enough $n \in \cN$ that
\bea
\Pr[ [k_n] \setminus I_n(r_n^{-1} \X_n) \neq \emptyset ] \leq 
k_n \exp(- (\delta/d)^d 
n^{1/8}) \leq   n \exp(- (\delta/d)^d  n^{1/8}).
\label{0102b}
\eea
Also, using (\ref{1220b}) and then (\ref{0102a}),
 for all large enough $n \in \cN$ we have 
% and (\ref{1220a}) 
%we obtain
%that almost surely, for large enough $n \in \cN$ we have
\bean
%\Pr[ r_n^d \zeta_{r_n}(\X_n)  > (1-4 \delta)^{-d}\okappa +  \delta ]  = 
& & \Pr[ r_n^d \zeta(r_n^{-1} \X_n)  > (1-4 \delta)^{-d}\okappa +  \delta ]  
\nonumber \\
&  \leq & \Pr[ r_n^d (k_n + K_0  |[k_n ] \setminus I_n(r_n^{-1} \X_n) | 
) > (1-4\delta)^{-d} \okappa + \delta ]  
\nonumber \\
& \leq & \Pr[ [k_n] \setminus I_n(r_n^{-1} \X_n) \neq \emptyset ],
\eean
which is summable in $n$ 
 by (\ref{0102b}).
Since $\delta >0$ is arbitrarily small, (\ref{eqlimsup1}) follows.
\end{proof}

\begin{lemm}
\label{lemlimsup2}
Let $\eps > 0$.
% Assume $n^{-1/d} \ll r_n \ll 1$. 
Then
\bea
\sum_{n \in \N \setminus \cN}
\Pr[ r_n^d \zeta_{r_n}(\X_n) > \okappa + \eps]  < \infty.
\label{eqlimsup1a}
\eea
\end{lemm}
\begin{proof}
We first consider the point process 
  $\Po_n$ defined earlier; since $\mu= \mu_U$ here,
 $\Po_n$ is a homogeneous Poisson process in $Q_1$ of intensity $n$, so 
by the Mapping theorem 
 $r_n^{-1} \Po_n$ is a homogeneous Poisson process in $Q_{r_n^{-1}} $
 of intensity $nr_n^d$. By (\ref{0101a}),  
we  have for all 
$n \geq n_0$ 
and
$i \in [k_n]$ that
$ \Pr[i \notin I_n(r_n^{-1} \Po_n)] \leq p_n$, where we set
$$
p_n := \exp (- n  (\delta/d)^d r_n^d) = o(1),
$$  
because we assume $n r_n^d \to \infty$.

Also the balls $B_\delta(x_{n,i}), 1 \leq i \leq k_n$, are disjoint, and setting
 $M_n := |[k_n] \setminus I_n(r_n^{-1} \Po_n) | $, we have that
$M_n$ is stochastically dominated  by a 
Binomial with parameters $k_n$ and $p_n$.
 Hence,  given  any fixed $\delta' >0$,
by (1.7) of \cite{RGG}, 
\bean
%\Pr[ | [k'_n] \setminus I_n(r_n^{-1} \Po_n) | \geq \delta' k_n ]
%&\leq &
 \Pr[  M_n  \geq \delta' k_n  ]
%\nonumber
% \\
%&
 \leq \exp \left(
 %\left( - \left(\frac{\delta' k_n}{2} \right) 
  - (\delta' k_n/2 ) 
\log ( \delta'/ p_n) \right) 
=  \exp (- \omega( k_n )).
% = \exp( -\omega(r_n^{-d} ) ).
\eean
Therefore by (\ref{0102a}),  
for large enough $n \in \N \setminus \cN$ we have 
\bea
\Pr[ M_n
%| [k'_n] \setminus I_n(r_n^{-1} \Po_n) |
 \geq \delta' k_n ]
%\leq
%\leq  \exp (- \omega( k_n )) 
%= \exp( -\omega(r_n^{-d} ) )
%\nonumber \\
\leq \exp(- n^{7/8}).
\label{1221b}
\eea
%which is summable in $n$.
By 
(\ref{1220b}) and (\ref{0102a}),
% and Lemma \ref{lemknkn}, 
 for all large enough $n$
% $n \in \N \setminus \cN$
 we have
\bean
\Pr \left[
r_n^d \zeta(r_n^{-1} \Po_n)  \geq (1-4\delta)^{-d}\okappa +  
 2 \delta \right]
~~~~~~~~~~~~~~~~~~~~~~~~~~~~~~~~~~~~~~~~~~~~~~~~~~~
\\
\leq \Pr[r_n^d( k_n + K_0 |[k_n] \setminus I_n(r_n^{-1} \Po_n) | )
> (1-4 \delta)^{-d} \okappa + 2 \delta ]
%\\
 \leq  
\Pr[  r_n^d K_0 M_n
% |[k'_n] \setminus I_n(r_n^{-1}\Po_n) | 
\geq \delta ].
% k_n ].
%\label{1220c}
\eean
Therefore since
$\zeta(r_n^{-1} \Po_n) =
\zeta_{r_n}(\Po_n) $,
for large enough $n \in \N \setminus \cN$ we have
 by (\ref{0102a}) and (\ref{1221b}) that 
\bea
\Pr \left[ r_n^d \zeta_{r_n}( \Po_n)  \geq (1-4\delta)^{-d}\okappa +  
 2 \delta \right] \leq  
\exp(- n^{7/8}).
\label{1220c}
\eea
%which is summable in $n$.

With our coupling of $\X_n$ and $\Po_n$, 
 by \cite[Lemma 1.2]{RGG} and Taylor expansion about 1 of
the function $h(x) := 1-x + x \log x$ we have for large enough $n$ that
\bean
%\Pr[ r_n^d |\gamma_{r_n}(\X_n) - \gamma_{r_n}(\Po_n) | 
%> K r_n^d n^{3/4} ] 
%\nonumber \\
\Pr[ |\X_n \triangle \Po_n| > n^{3/4} ] 
\leq 
\exp \left( -n h \left( \frac{n+ n^{3/4}}{n} \right) \right)   
+
\exp \left( -n h \left( \frac{n- n^{3/4}}{n} \right) \right)   
\nonumber \\
\leq 2 \exp( - n^{1/2 } /3). 
\eean
 By Lemma \ref{lemA1C},
if $ |\X_n \triangle \Po_n| \leq n^{3/4} $
then $  |\zeta_{r_n}(\X_n) - \zeta_{r_n}(\Po_n) | 
\leq K n^{3/4}$. Hence
\bea
\Pr[ r_n^d |\zeta_{r_n}(\X_n) - \zeta_{r_n}(\Po_n) | 
> K r_n^d n^{3/4} ] \leq
 2 \exp( - n^{1/2 } /3). 
%\Pr[ |\X_n \triangle \Po_n| > n^{3/4} ], 
\label{1228a}
\eea
%which is summable in $n$.
Since $K r_n^d n^{3/4} \leq K n^{-1/8} < \delta $ for all large enough
$n \in \N \setminus \cN$,
 combining (\ref{1228a}) with
%(\ref{1221b}) and
 (\ref{1220c}) shows that for all large enough $n \in \N \setminus \cN$,
$$
\Pr \left[ r_n^d \zeta_{r_n}( \X_n)  \geq (1-4\delta)^{-d}\okappa +  
 3 \delta \right] \leq  
 \exp(- n^{7/8})  + 2 \exp(- n^{1/2}/3), 
$$
which is summable in $n$. Since $\delta $ can be taken
 arbitrarily small, this
 gives us the result.
\end{proof}

\begin{proof}[Proof of Theorem \ref{thsubcond}]
By Lemma \ref{lemlow1}, for all $n \in \N$ we have 
\bea
\Pr[ r_n^{d} \zeta_{r_n}(\X_n) < (1+  d \delta)^{-d} \okappa  ]  
\leq  
\Pr[ 
%r_n^{d} |
[\ell_n] \setminus J_n(r_n^{-1} \X_n)
%| > \delta 
\neq \emptyset ].
\label{0108a}
\eea
For all sufficiently large $n \in \cN$ and all $i \in [\ell_n]$,
\bean
\Pr[ i \notin J_n(r_n^{-1} \X_n)] \leq \left( 1 - \left( \frac{ \delta r_n}{2 } 
\right)^d \right)^n 
\leq 
\exp \left[ - \left(\frac{\delta}{2}\right)^d n r_n^d \right] \leq 
\exp \left[ - \left(\frac{\delta}{2}\right)^d n^{1/8} \right].
\eean

Thus
by the union bound, since $\ell_n = O(r_n^{-d})=
o(n)$ by (\ref{1230d}), 
%there is a constant $c$ such that
for all sufficiently large $n \in \cN$ we have 
\bea 
\Pr[ [\ell_n] \setminus J_n(r_n^{-1} \X_n) \neq \emptyset ] \leq  n 
\exp\left( - (\delta/2)^d  n^{1/8} \right).
\label{1221d}
\eea

Suppose  $n \in \N \setminus \cN$, so that $nr_n^d \leq n^{1/8}$. 
Then $|[\ell_n] \setminus J_n(r_n^{-1} \Po_n)|$ is binomial with parameters $\ell_n$
and $q_n$, where we set $q_n:= \exp(- n \delta_n^dr_n^d)
$
so $q_n \to 0$ as  $n \to \infty$. Therefore, similarly to
(\ref{1221b})  and using 
(\ref{1230d}), 
%given $\delta' >0$
 we have for all 
large enough 
$n \in \N \setminus \cN $ that
\bea
\Pr[ | [\ell_n] \setminus J_n(r_n^{-1} \Po_n) | \geq \delta r_n^{-d} ]
& \leq &  \exp ( - (\delta r_n^{-d}/2) \log [\delta r_n^{-d}/(\ell_n q_n) ] ) 
\nonumber \\
 & \leq  &   \exp( - n^{7/8}  ).
\label{0101b}
\eea
%so by the Borel-Cantelli lemma, almost surely we have 
%$r_n^{d} |[\ell_n] \setminus J_n(\Po_n)| \to 0$ as $n \to \infty$
%through $\N \setminus \cN$ (if $\N \setminus \cN$ is finite).
By Lemma \ref{lemlow1},
we have for all large enough $n \in \N \setminus \cN$ that
\bea
\Pr[ r_n^{d} \zeta(r_n^{-1} \Po_n)
 < (1+ d \delta)^{-d} \okappa -  \delta ]
\leq  \Pr [ | [\ell_n] \setminus J_n(r_n^{-1} \Po_n) | \geq \delta r_n^{-d} ].
\label{1221g}
\eea
As before, we have
(\ref{1228a}) for all 
 $n $, and $K r_n^d n^{3/4} \leq K n^{-1/8} \leq \delta $
for all large enough $n \in \N \setminus \cN$. Hence
by (\ref{0101b}) and (\ref{1221g}), 
for all  large enough  $n \in \N \setminus \cN$ we have
$$
\Pr[r_n^d \zeta_{r_n}(\X_n) < (1+  d \delta)^{-d} \okappa - 2 \delta]
\leq \exp(-n^{7/8}) + 2 \exp (- n^{1/2}/3) .
$$
Together with (\ref{0108a}) and (\ref{1221d}),
this shows that 
$$
\sum_{n=1}^\infty \Pr[ r_n^d \zeta_{r_n}(\X_n) < (1+  d \delta)^{-d} \okappa
- 2\delta] < \infty.
$$
Combined with Lemmas \ref{lemlimsup1} and \ref{lemlimsup2} and
using the Borel-Cantelli lemma, since
$\delta >0$ can be arbitrarily small this 
shows that $r_n^d \zeta_{r_n}(\X_n) \to \okappa$ almost surely, which
gives us (\ref{1221h}) as required.
%$$
%r_n^{d} \gamma(G(\X_n,r_n))
% \geq (1+\delta)^{-d}( (1-\delta)(\psi(B)/\pi_d) -  \delta). 
%$$
\end{proof}
%}

%Mathew D. Penrose, Department of Mathematical Sciences,
%University of Bath, Bath BA2 7AY United Kingdom:
%{\texttt mathew.penrose@durham.ac.uk}

\end{document}